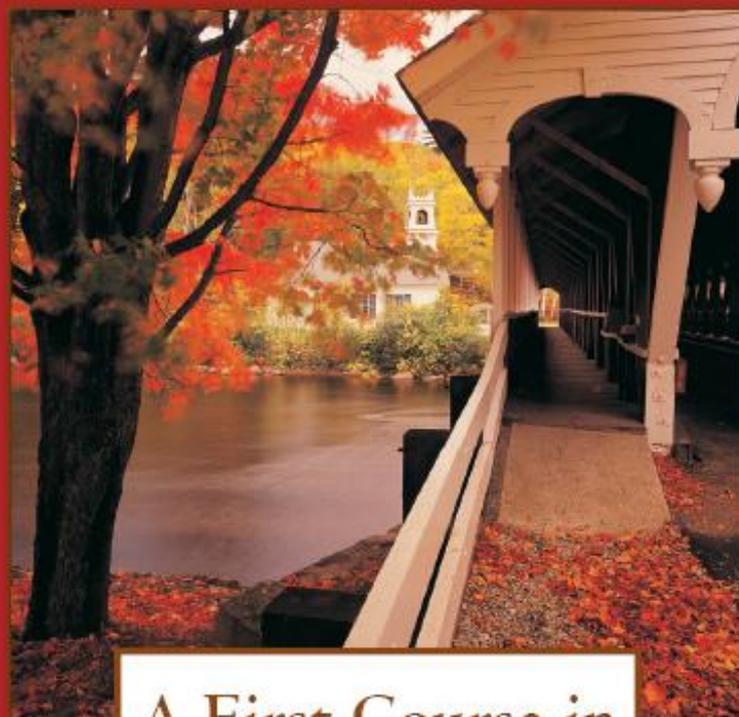

# A First Course in Linear Algebra

By

Mohammed K A Kaabar

# About the Author

**Mohammed Kaabar** is a math tutor at the Math Learning Center (MLC) at Washington State University, Pullman, and he is interested in linear algebra, scientific computing, numerical analysis, differential equations, and several programming languages such as SQL, C#, Scala, C++, C, JavaScript, Python, HTML 5 and MATLAB. He is a member of of Institute of Electrical and Electronics Engineers (IEEE), IEEE Antennas and Propagation Society, IEEE Consultants Network, IEEE Smart Grid Community, IEEE Technical Committee on RFID, IEEE Life Sciences Community, IEEE Green. ICT Community, IEEE Cloud Computing Community, IEEE Internet of Things Community, IEEE Committee on Earth Observations, IEEE Electric Vehicles Community, IEEE Electron Devices Society, IEEE Communications Society, and IEEE Computer Society. He participated in several competitions, conferences, research papers and projects. He is an online instructor of numerical analysis at Udemy Inc, San Francisco, CA. In addition, he is also a Technical Program Committee (TPC) member, reviewer and presenter at CCA-2014, WSMEAP 2014, EECSI 2014, JIEEEC 2013 and WCEEENG 2012. He worked as electrical engineering intern at Al-Arabia for Safety and Security L.L.C. He also received several educational awards and certificates from accredited institutions. For more information about the author and his free online courses, please visit his personal website: http://www.mohammed-kaabar.net.

# Table of Contents

# Introduction

In this book, I wrote five chapters: Systems of Linear Equations, Vector Spaces, Homogeneous Systems, Characteristic Equation of Matrix, and Matrix Dot Product. I also added exercises at the end of each chapter above to let students practice additional sets of problems other than examples, and they can also check their solutions to some of these exercises by looking at "Answers to Odd-Numbered Exercises" section at the end of this book. This book is very useful for college students who studied Calculus I, and other students who want to review some linear algebra concepts before studying a second course in linear algebra. According to my experience as a math tutor at Math Learning Center, I have noticed that some students have difficulty to understand some linear algebra concepts in general, and vector spaces concepts in particular. Therefore, my purpose is to provide students with an interactive method to explain the concept, and then provide different sets of examples related to that concept. If you have any comments related to the contents of this book, please email your comments to mohammed.kaabar@email.wsu.edu.

I wish to express my gratitude and appreciation to my father, my mother, and my brother. I would also like to give a special thanks to my mathematics professor Dr. Ayman Rateb Badawi for his brilliant efforts in revising the content of this book, and I would also like to thank all my mathematics professors who taught me math courses at Washington State University, Pullman and American University of Sharjah. Ultimately, I would appreciate to consider this book as a milestone for devolving more math books that can serve our mathematical society.

# Chapter 1
# Systems of Linear Equations

In this chapter, we discuss how to solve $n \times m$ systems of linear equations using row operations method. Then, we give an introduction to basic algebra of matrix including matrix addition, matrix subtraction and matrix multiplication. We cover in the remaining sections some important concepts of linear algebra such as linear combinations, determinants, square matrix, inverse square matrix, transpose matrix, Cramer's Rule and Adjoint Method.

## 1.1 Row Operations Method

First of all, let's start with a simple example about $n \times m$ systems of linear equation.

**Example 1.1.1** Solve for x and y for the following $2 \times 2$ system of linear equations:
$$\begin{cases} 3x + 2y = 5 \\ -2x + y = -6 \end{cases}$$

**Solution:** Let's start analyzing this $2 \times 2$ system.

First of all, each variable in this system is to the power 1 which means that this system is a linear system. As given in the question itself, $2 \times 2$ implies that the number of equations is 2, and the number of unknown variables is also 2.

$$2 \times 2$$

(**Number of Equations**) × (**Number of Unkown Variables**)

Then, the unknown variables in this question are x and y. To solve for x and y, we need to multiply the first equation $3x + 2y = 5$ by 2, and we also need to multiply the second equation $-2x + y = -6$ by 3. Hence, we obtain the following:

$$\begin{cases} 6x + 4y = 10 \ldots\ldots\ldots\ldots\ldots..1 \\ -6x + 3y = -18 \ldots\ldots\ldots\ldots.2 \end{cases}$$

By adding equations 1 and 2, we get the following:

$$7y = -8$$

Therefore, $y = \frac{-8}{7}$

Now, we need to substitute the value of y in one of the two original equations. Let's substitute y in the first equation $3x + 2y = 5$ as follows:

$3x + 2\left(\frac{-8}{7}\right) = 5$  is equivalent to

$3x = 5 - 2\left(\frac{-8}{7}\right) = 5 + \left(\frac{16}{7}\right) = \frac{51}{7}$  Then,  $x = \frac{17}{7}$

Therefore, $x = \frac{17}{7}$ and $y = \frac{-8}{7}$

Hence, we solve for x and y.

**Example 1.1.2** Describe the following system:
$$\begin{cases} 2x + 3y - 2z = 10 \\ 5x + 2y + z = 13 \end{cases}$$

**Solution:** Let's start asking ourselves the following questions.

Question 1: How many equations do we have?
Question 2: How many unknown variables do we have?
Question 2: What is the power of each unknown variable?

If we can answer the three questions above, then we can discuss the above systems.

Let's start now answering the three questions above.

Answer to Question 1: We have 2 equations: $2x + 3y - 2z = 10$ and $5x + 2y + z = 13$.

Answer to Question 2: We have 3 unknown variables: x, y and z.

Answer to Question 3: Each unknown variable is to the power 1.

After answering the above three question, we have now an idea about the above system. As we remember from example 1.1.1 that the system of equations has the following form:

(**Number of Equations**) × (**Number of Unkown Variables**)

Since the number of equations in this example is 2 and the number of unknown variables is 3, then using the above form, the above system is $2 \times 3$ system. Recall that each unknown variable is to the power 1. This means that it is a $2 \times 3$ linear system.

Solving for unknown variables in $2 \times 2$ or $2 \times 3$ systems of linear equations using some arithmetic operations such as additions, subtractions and multiplications of linear equations is very easy. But if we have for example $4 \times 4$ system of linear equations, using some arithmetic operations in this case will be very difficult and it will take long time to solve it. Therefore, we are going to use a method called Row Operation Method to solve complicated $n \times m$ systems of linear equations. Now, let's start with an example discussing each step of Row Operation Method for solving $n \times m$ system of linear equations.

**Example 1.1.3** Solve for $x_1$, $x_2$ and $x_3$ for the following $3 \times 3$ system of linear equations:

$$\begin{cases} x_1 + x_2 - x_3 = 2 \\ 2x_1 - x_2 + x_3 = 2 \\ -x_2 + 2x_3 = 1 \end{cases}$$

**Solution:** In the above system, we have 3 equations and 3 unknown variables $x_1$, $x_2$ and $x_3$. To solve this $3 \times 3$ system, we need to use Row Operation Method. The following steps will describe this method.

Step 1: Rewrite the system as a matrix, and this matrix is called Augmented Matrix. Each row in this matrix is a linear equation.

$$\begin{matrix} x_1 & x_2 & x_3 & c \end{matrix}$$
$$\begin{pmatrix} 1 & 1 & -1 & | & 2 \\ 2 & -1 & 1 & | & 2 \\ 0 & -1 & 2 & | & 1 \end{pmatrix}$$

As we can see from the above augmented matrix, the first column represents the coefficients of $x_1$ in the three linear equations. The second column represents the coefficients of $x_2$ in the three linear equations. The third column represents the coefficients of $x_3$ in the three linear equations. The fourth column is not an actual column but it just represents the constants because each linear equation equals to a constant.

Step 2: Start with the first row in matrix and make the first non-zero number equals to 1. This 1 is called a leader number.

$$\begin{array}{cccc} x_1 & x_2 & x_3 & c \end{array}$$
$$\begin{pmatrix} \boxed{1} & 1 & -1 & 2 \\ 2 & -1 & 1 & 2 \\ 0 & -1 & 2 & 1 \end{pmatrix}$$

Since the leader number is already 1 in this matrix, then we do not need to do anything with this number. Otherwise, we need to make this leader number equals to 1 by multiplying the whole row with a non-zero number.

Step 3: Eliminate all numbers that are exactly below the leader number in other words use the leader number to eliminate all things below it.

Step 4: Move to the second row and repeat what has been done in step 2.

$$-2R_1+R_2 \dashrightarrow R_2$$

(This means that we multiply the first row by -2 and we add it to the second row and the change will be only in the second row and no change in the first row)

Hence, we obtain:

$$\begin{pmatrix} 1 & 1 & -1 & | & 2 \\ 0 & -3 & 3 & | & -2 \\ 0 & -1 & 2 & | & 1 \end{pmatrix}$$

$$\frac{-1}{3} R_2$$

(This means that we multiply the second row by $\frac{-1}{3}$ and the change will be only in the second row and no change in other rows)

Hence, we obtain:

$$\begin{pmatrix} 1 & 1 & -1 & | & 2 \\ 0 & \boxed{1} & -1 & | & \frac{2}{3} \\ 0 & -1 & 2 & | & 1 \end{pmatrix}$$

$$R_2+R_3 \dashrightarrow R_3$$

(This means that we add the second row to the third row and the change will be only in the third row and no change in the second row)

$$-R_2+R_1 \dashrightarrow R_1$$

(This means that we multiply the second row by -1, and we add it to the first row and the change will be only in the first row and no change in the second row)

Hence, we obtain:

$$\begin{pmatrix} 1 & 0 & 0 & \bigg| & \frac{4}{3} \\ 0 & 1 & -1 & \bigg| & \frac{2}{3} \\ 0 & 0 & 1 & \bigg| & \frac{5}{3} \end{pmatrix}$$

Step 5: Move to the third row and do the same as we did in the first row and the second row of matrix.

$$R_3+R_2 \dashrightarrow R_2$$

(This means that we add the third row to the second row and the change will be only in the second row and no change in the third row)

Hence, we obtain:

$$\begin{pmatrix} 1 & 0 & 0 & \bigg| & \frac{4}{3} \\ 0 & 1 & 0 & \bigg| & \frac{7}{3} \\ 0 & 0 & 1 & \bigg| & \frac{5}{3} \end{pmatrix}$$

Therefore, $x_1 = \frac{4}{3}$ , $x_2 = \frac{7}{3}$ and $x_3 = \frac{5}{3}$

These are one solution only.

Hence, the system has a unique solution (one solution).

After this example, we will get introduced to two new definitions.

**Definition 1.1.1** $n \times m$ system of linear equations is consistent if it has a solution; otherwise, it is called inconsistent.

**Definition 1.1.2** $n \times m$ system of linear equations has a unique solution if each variable has exactly one value.

Now, let's apply what we have learned from example 1.1.3 for the following example 1.1.4.

**Example 1.1.4** Solve for $x_1$, $x_2$, $x_3$, $x_4$ and $x_5$ for the following $3 \times 5$ system of linear equations:

$$\begin{cases} x_2 - x_3 + x_4 - x_5 = 1 \\ -2x_1 + x_3 - x_5 = 0 \\ -x_2 + x_3 + 2x_4 - 10x_5 = 12 \end{cases}$$

**Solution:** In the above system, we have 3 equations and 5 unknown variables $x_1$, $x_2$, $x_3$, $x_4$ and $x_5$. To solve this $3 \times 5$ system, we need to use Row Operation Method. First, we need to construct the augmented matrix.

$$\begin{array}{ccccccc} x_1 & x_2 & x_3 & x_4 & x_5 & c \end{array}$$
$$\begin{pmatrix} 0 & 1 & -1 & 1 & -1 & | & 1 \\ -2 & 0 & 1 & 0 & -1 & | & 0 \\ 0 & -1 & 1 & 2 & -10 & | & 12 \end{pmatrix}$$

The leader number here is 1.

$$\begin{pmatrix} 0 & \boxed{1} & -1 & 1 & -1 & | & 1 \\ -2 & 0 & 1 & 0 & -1 & | & 0 \\ 0 & -1 & 1 & 2 & -10 & | & 12 \end{pmatrix}$$

By doing the same steps in example 1.1.3, we obtain the following:

$$R_1+R_3 \dashrightarrow R_3$$

$$\begin{array}{cccccc} x_1 & x_2 & x_3 & x_4 & x_5 & c \end{array}$$
$$\begin{pmatrix} 0 & 1 & -1 & 1 & -1 & | & 1 \\ -2 & 0 & 1 & 0 & -1 & | & 0 \\ 0 & 0 & 0 & 3 & -11 & | & 13 \end{pmatrix}$$

$$\frac{-1}{2} R_2$$

$$\begin{array}{cccccc} x_1 & x_2 & x_3 & x_4 & x_5 & c \end{array}$$
$$\begin{pmatrix} 0 & 1 & -1 & 1 & -1 & | & 1 \\ 1 & 0 & \frac{-1}{2} & 0 & \frac{1}{2} & | & 0 \\ 0 & 0 & 0 & 3 & -11 & | & 13 \end{pmatrix}$$

$$\frac{1}{3} R_3$$

$$\begin{array}{cccccc} x_1 & x_2 & x_3 & x_4 & x_5 & c \end{array}$$
$$\begin{pmatrix} 0 & 1 & -1 & 1 & -1 & | & 1 \\ 1 & 0 & -1/2 & 0 & 1/2 & | & 0 \\ 0 & 0 & 0 & 1 & -11/3 & | & \frac{13}{3} \end{pmatrix}$$

$$-R_3+R_1 \dashrightarrow R_1$$

$$\begin{array}{cccccc} x_1 & x_2 & x_3 & x_4 & x_5 & c \end{array}$$
$$\begin{pmatrix} 0 & \boxed{1} & -1 & 0 & 8/3 & | & -10/3 \\ \boxed{1} & 0 & -1/2 & 0 & 1/2 & | & 0 \\ 0 & 0 & 0 & \boxed{1} & -11/3 & | & \frac{13}{3} \end{pmatrix}$$

Hence, we obtain from the above matrix:

$$x_2 - x_3 + \frac{8}{3}x_5 = -\frac{10}{3}$$
$$x_1 - \frac{1}{2}x_3 + \frac{1}{2}x_5 = 0$$
$$x_4 - \frac{11}{3}x_5 = \frac{13}{3}$$

Since we have 3 leader numbers (numbers equal to 1) in $x_1$, $x_2$ and $x_4$ columns, then $x_1$, $x_2$ and $x_4$ are called leading variables. All other variables are called free variables such that $x_3, x_5 \in \mathbb{R}$.

Since we have leading and free variables, we need to write the leading variables in terms of free variables.

$$x_1 = \frac{1}{2}x_3 - \frac{1}{2}x_5$$
$$x_2 = x_3 - \frac{8}{3}x_5 - \frac{10}{3}$$
$$x_4 = \frac{11}{3}x_5 + \frac{13}{3}$$

Now, let $x_3 = 0$ and $x_5 = 0$ (You can choose any value for $x_3$ and $x_5$ because both of them are free variables). Therefore, we obtain

$$x_1 = \frac{1}{2}(0) - \frac{1}{2}(0) = 0$$
$$x_2 = 0 - \frac{8}{3}(0) - \frac{10}{3} = -\frac{10}{3}$$
$$x_4 = \frac{11}{3}(0) + \frac{13}{3} = \frac{13}{3}$$

Hence, $x_1 = 0, x_2 = -\frac{10}{3}$ and $x_4 = \frac{13}{3}$

Another possible solution for example 1.1.4:

Now, let $x_3 = 1$ and $x_5 = 0$ (You can choose any value for $x_3$ and $x_5$ because both of them are free variables). Therefore, we obtain

$$x_1 = \frac{1}{2}(1) - \frac{1}{2}(0) = \frac{1}{2}$$
$$x_2 = 1 - \frac{8}{3}(0) - \frac{10}{3} = 1 - \frac{10}{3} = -\frac{7}{3}$$
$$x_4 = \frac{11}{3}(0) + \frac{13}{3} = \frac{13}{3}$$

Hence, we obtain

$$x_1 = \frac{1}{2}, x_2 = -\frac{7}{3} \text{ and } x_4 = \frac{13}{3}$$ This is another solution.

**Summary of Row Operations Method in $n \times m$ Systems of Linear Equations:**

Suppose $\propto$ is a non-zero constant, and $i$ and $k$ are row numbers in the augmented matrix.

* $\propto R_i$ , $\propto \neq 0$ (Multiply a row with a non-zero constant $\propto$).

* $R_i \leftrightarrow R_k$ (Interchange two rows).

* $\propto R_i + R_k \dashrightarrow R_k$ (Multiply a row with a non-zero constant $\propto$, and add it to another row).

Note: The change is in $R_k$ and no change is in $R_i$.

Let's start with other examples that will give us some important results in linear algebra.

**Example 1.1.5** Solve for $x_1$, $x_2$ and $x_3$ for the following $2 \times 3$ system of linear equations:
$$\begin{cases} x_1 + 2x_2 - x_3 = 2 \\ 2x_1 + 4x_2 - 2x_3 = 6 \end{cases}$$

**Solution:** In the above system, we have 2 equations and 3 unknown variables $x_1$, $x_2$ and $x_3$. To solve this $2 \times 3$ system, we need to use Row Operation Method. First, we need to construct the augmented matrix.

$$\begin{array}{cccc} x_1 & x_2 & x_3 & c \end{array}$$
$$\begin{pmatrix} \boxed{1} & 2 & -1 & | & 2 \\ 2 & 4 & -2 & | & 6 \end{pmatrix}$$

$-2R_1 + R_2 \dashrightarrow R_2$

$$\begin{array}{cccc} x_1 & x_2 & x_3 & c \end{array}$$
$$\begin{pmatrix} 1 & 2 & -1 & | & 2 \\ 0 & 0 & 0 & | & 2 \end{pmatrix}$$

We stop here and read the above matrix as follows:
$$x_1 + 2x_2 - x_3 = 2$$
$$0 = 2$$

Hence, we get introduced to a new result.

**Result 1.1.1** The system is inconsistent if and only if after reducing one of the equations, we have 0 equals to a non-zero number.

**Example 1.1.6** Solve for $x_1$, $x_2$ and $x_3$ for the following $2 \times 3$ system of linear equations:
$$\begin{cases} x_1 + 2x_2 - x_3 = 2 \\ 2x_1 + 4x_2 - 2x_3 = 4 \end{cases}$$

**Solution:** In the above system, we have 2 equations and 3 unknown variables $x_1$, $x_2$ and $x_3$. To solve this $2 \times 3$ system, we need to use Row Operation Method. First, we need to construct the augmented matrix as we did in the previous example.

$$\begin{array}{cccc} x_1 & x_2 & x_3 & c \end{array}$$
$$\begin{pmatrix} \boxed{1} & 2 & -1 & | & 2 \\ 2 & 4 & -2 & | & 4 \end{pmatrix}$$

$-2R_1+R_2 \dashrightarrow R_2$

$$\begin{array}{cccc} x_1 & x_2 & x_3 & c \end{array}$$
$$\begin{pmatrix} 1 & 2 & -1 & | & 2 \\ 0 & 0 & 0 & | & 0 \end{pmatrix}$$

We stop here and read the above matrix as follows:
$$x1 + 2x2 - x3 = 2$$
$$0 = 0$$

The solution is $x_1 = -2 x_2 + x_3 + 2$

$x_1$ is a leading variable, and $x_2, x_3 \in \mathbb{R}$ free variables.

Now, let $x_2 = 1$ and $x_3 = 0$ (You can choose any value for $x_2$ and $x_3$ because both of them are free variables).

Therefore, we obtain $x1 = -2(0) + (0) + 2 = 2$ and
$$x_2 = 0 \text{ and } x_3 = 0.$$

Hence, we get introduced to new results.

**Result 1.1.2** The $n \times m$ system is consistent if and only if it has a unique solution (no free variable).

**Result 1.1.3** The $n \times m$ system is consistent if and only if it has infinitely many solutions (free variables).

**Result 1.1.4** Assume the $n \times m$ system is consistent, and we have more variables than equations (m>n). Then, the system has infinitely many solutions.

**Example 1.1.7** Given an augmented matrix of a system: $\begin{pmatrix} 2 & 5 & 3 & | & -1 \\ -2 & 4 & b & | & 3 \\ -2 & -5 & -3 & | & d \end{pmatrix}$

For what values of b and d will the system be consistent?

**Solution:** We need to multiply the leader number (2) by a non-zero constant to make it equal to 1 instead of 2 as follows:

$$\tfrac{1}{2} R_1$$

$$\begin{pmatrix} 1 & \tfrac{5}{2} & \tfrac{3}{2} & \big| & -\tfrac{1}{2} \\ -2 & 4 & b & \big| & 3 \\ -2 & -5 & -3 & \big| & d \end{pmatrix}$$

$$2R_1 + R_2 \dashrightarrow R_2$$
$$2R_1 + R_3 \dashrightarrow R_3$$

$$\begin{pmatrix} 1 & \tfrac{5}{2} & \tfrac{3}{2} & \big| & -\tfrac{1}{2} \\ 0 & 9 & b+3 & \big| & 2 \\ 0 & 0 & 0 & \big| & d-1 \end{pmatrix}$$

Therefore, $b \in \mathbb{R}$ and d=1

When the system is consistent, then we must have infinitely many solutions. In this case, b is a free variable.

**Example 1.1.8** Given an augmented matrix of a system: $\begin{pmatrix} 1 & -1 & 2 & \big| & a \\ -1 & 2 & -2 & \big| & b \\ -2 & 2 & -4 & \big| & c \end{pmatrix}$

For what values of a, b and c will the system be consistent?

**Solution:** We need to use the Row Operation Method to find a, b and c as follows:

$$R_1 + R_2 \dashrightarrow R_2$$

$$2R_1+R_3 \dashrightarrow R_3$$

$$\begin{pmatrix} 1 & -1 & 2 & | & a \\ 0 & 1 & 0 & | & a+b \\ 0 & 0 & 0 & | & 2a+c \end{pmatrix}$$

$$R_2+R_1 \dashrightarrow R_1$$

$$\begin{pmatrix} 1 & 0 & 2 & | & 2a+b \\ 0 & 1 & 0 & | & a+b \\ 0 & 0 & 0 & | & 2a+c \end{pmatrix}$$

We stop here and read the above matrix. Hence, we obtain:

$x_1 + 2x_3 = 2a + b$

$x_2 = a + b$

$0 = 2a + c$

Therefore, $2a + c = 0$ which means that $c = -2a$

Hence, the solution is $a, b \in \mathbb{R}$ (free variables), and $c = -2a$.

Let's now discuss some new definitions and results of linear algebra.

**Definition 1.1.3** Suppose we have the following $3 \times 3$ system of linear equations:

$$\begin{cases} 2x_1 - 3x_2 + x_3 = 0 \\ x_1 + x_2 - 3x_3 = 0 \\ -x_1 + 3x_2 - x_3 = 0 \end{cases}$$

This system is called homogeneous system of linear equations because all constants on the right-side are zeros.

**Result 1.1.5** Every $n \times m$ homogeneous system is consistent.

**Result 1.1.6** Since every $n \times m$ homogeneous system is consistent, then $x_1=0$, $x_2=0$, $x_3=0$, …, $x_m=0$ is a solution. This solution is called Trivial Solution.

**Result 1.1.7** The solution is called Non-Trivial if it has at least one value that is not zero.

**Result 1.1.8** In a homogeneous system if number of variables is greater than number of equations, then the system has infinitely many solutions.

**Result 1.1.9** Given the following $3 \times 6$ matrix:

$$\begin{bmatrix} 0 & \boxed{1} & 0 & 0 & 1 & 1 \\ \boxed{1} & 0 & 3 & 2 & 5 & 6 \\ 0 & 0 & 0 & 0 & 0 & \boxed{1} \end{bmatrix}$$ This matrix is called Reduced Matrix because all the numbers exactly below each leader number are zeros.

**Result 1.1.10** Given the following $3 \times 6$ matrix:

$$\begin{bmatrix} 0 & 0 & \boxed{1} & 2 & 0 \\ \boxed{1} & 5 & 0 & 2 & 0 \\ 0 & 0 & 0 & 0 & \boxed{1} \end{bmatrix}$$ This matrix is called Completely-Reduced Matrix because all the numbers exactly below and above each leader number are zeros.

**Result 1.1.11** The matrix is called Echelon Matrix if it is a combination of Reduced Matrix and Leader Numbers that are arranged so that the leader number in i$^{\text{th}}$ row is to the right of the leader number in the previous row, and the rows that are entirely zeros come at the end of the matrix.

$$\begin{bmatrix} 0 & \boxed{1} & 2 & 3 & 7 \\ 0 & 0 & \boxed{1} & 5 & 6 \\ 0 & 0 & 0 & \boxed{1} & 3 \\ 0 & 0 & 0 & 0 & \boxed{1} \\ 0 & 0 & 0 & 0 & 0 \\ 0 & 0 & 0 & 0 & 0 \end{bmatrix}$$ This is an Echelon Matrix.

$$\begin{bmatrix} 0 & 1 & 0 & 0 & 1 & 1 \\ 1 & 0 & 3 & 2 & 5 & 6 \\ 0 & 0 & 0 & 0 & 0 & 1 \end{bmatrix}$$ This is NOT Echelon Matrix.

**Example 1.1.9** Given the following reduced matrix:
$$\begin{bmatrix} 0 & 1 & 0 & 0 & 1 & 1 \\ 1 & 0 & 3 & 2 & 5 & 6 \\ 0 & 0 & 0 & 0 & 0 & 1 \end{bmatrix}$$
Covert this reduced matrix to echelon matrix.
**Solution:** We just need to interchange two rows as follows:
$$R_1 \longleftrightarrow R_2$$
$$\begin{bmatrix} \boxed{1} & 0 & 3 & 2 & 5 & 6 \\ 0 & \boxed{1} & 0 & 0 & 1 & 1 \\ 0 & 0 & 0 & 0 & 0 & 1 \end{bmatrix}$$
Now, the matrix is an echelon matrix.

**Result 1.1.12** To convert Reduced Matrix to Echelon Matrix, we just need to use interchange rows.

**Result 1.1.13** To convert Completely-Reduced Matrix to Reduced-Echelon Matrix, we just need to use interchange rows.

# 1.2 Basic Algebra of Matrix

In this section, we discuss an example of basic algebra of matrix including matrix addition, subtraction and multiplication.

**Example 1.2.1** Given the following three matrices:
$$A = \begin{bmatrix} 2 & 1 & 3 \\ 0 & 1 & 5 \end{bmatrix} \quad, \quad B = \begin{bmatrix} 1 & 2 \\ 6 & 1 \\ 3 & 1 \end{bmatrix}, \quad C = \begin{bmatrix} -3 & 2 & 5 \\ 1 & 7 & 3 \end{bmatrix}$$

a) Find 3A.
b) Find A+B.
c) Find A-B.
d) Find 3A-2C.

**Solution: Part a:** We just need to multiply 3 by matrix A as follows:

$3A = 3 * \begin{bmatrix} 2 & 1 & 3 \\ 0 & 1 & 5 \end{bmatrix} = \begin{bmatrix} 3*2 & 3*1 & 3*3 \\ 3*0 & 3*1 & 3*5 \end{bmatrix} = \begin{bmatrix} 6 & 3 & 9 \\ 0 & 3 & 15 \end{bmatrix}$

**Part b:** In matrix addition we can only add matrices of the same order only.

The order of matrix A is 2 × 3.

The order of matrix B is 3 × 2.

Since matrix A and matrix B have different orders, then we cannot find A+B.

Hence, there is no answer for part b.

**Part c:** In matrix subtraction we can only subtract matrices of the same order only.

The order of matrix A is 2 × 3.

The order of matrix B is 3 × 2.

Since matrix A and matrix B have different orders, then we cannot find A-B.

Hence, there is also no answer for part c.

**Part d:** First, we need to multiply 2 by matrix C as follows:

$2C = 2 * \begin{bmatrix} -3 & 2 & 5 \\ 1 & 7 & 3 \end{bmatrix} = \begin{bmatrix} 2*-3 & 2*2 & 2*5 \\ 2*1 & 2*7 & 2*3 \end{bmatrix} = \begin{bmatrix} -6 & 4 & 10 \\ 2 & 14 & 6 \end{bmatrix}$

Since we are already have 3A from part a, we just need to find 3A-2C. As we know, in matrix subtraction we can only subtract matrices of the same order only.

The order of matrix 3A is 2 × 3.

The order of matrix 2C is 2 × 3.

Both of them have the same order. Hence, we can find 3A-2C as follows:

$3A\text{-}2C = \begin{bmatrix} 6 & 3 & 9 \\ 0 & 3 & 15 \end{bmatrix} - \begin{bmatrix} -6 & 4 & 10 \\ 2 & 14 & 6 \end{bmatrix} = \begin{bmatrix} 12 & -1 & -1 \\ -2 & -11 & 9 \end{bmatrix}$

# 1.3 Linear Combinations

In this section, we discuss the concept of linear combinations of either columns or rows of a certain matrix.

**Example 1.3.1** Suppose we have 5 apples, 7 oranges and 12 bananas. Represent them as a linear combination of fruits.

**Solution:** To represent them as a linear combination, we need to do the following steps:

**Step 1:** Put each fruit individually:

        **Apples    Oranges    Bananas**

**Step 2:** Separate each one of them by addition sign.

        **Apples  +  Oranges  +  Bananas**

**Step 3:** Put 1 box in front of each one.

      ☐ **Apples**  +  ☐ **Oranges**  +  ☐ **Bananas**

**Step 4:** Write the number of each fruit in each box.

      [5] **Apples**  +  [7] **Oranges**  +  [12] **Bananas**

This representation is called a Linear Combination of Fruits.

Now, we can apply what we have learned in example 1.3.1 on rows and columns of a certain matrix.

**Example 1.3.2** Given the following matrix A:

$$A = \begin{bmatrix} 1 & 2 & 3 \\ 0 & 5 & 1 \\ 2 & 1 & 2 \end{bmatrix}$$

  a) Find a linear combination of the columns of matrix A.
  b) Find a linear combination of the rows of matrix A.

**Solution: Part a:** To represent the columns of matrix A as a linear combination, we need to do the same steps as we did in example 1.3.1.

**Step 1:** Put each column of matrix A individually:

$$\begin{bmatrix} 1 \\ 0 \\ 2 \end{bmatrix} \quad \begin{bmatrix} 2 \\ 5 \\ 1 \end{bmatrix} \quad \begin{bmatrix} 3 \\ 1 \\ 2 \end{bmatrix}$$

**Step 2:** Separate each one of them by addition sign.

$$\begin{bmatrix} 1 \\ 0 \\ 2 \end{bmatrix} + \begin{bmatrix} 2 \\ 5 \\ 1 \end{bmatrix} + \begin{bmatrix} 3 \\ 1 \\ 2 \end{bmatrix}$$

**Step 3:** Put 1 box in front of each one.

$$\square \begin{bmatrix} 1 \\ 0 \\ 2 \end{bmatrix} + \square \begin{bmatrix} 2 \\ 5 \\ 1 \end{bmatrix} + \square \begin{bmatrix} 3 \\ 1 \\ 2 \end{bmatrix}$$

**Step 4:** Write random number for each column in each box because it is not mentioned any number for any column in the question.

$$\boxed{5} \begin{bmatrix} 1 \\ 0 \\ 2 \end{bmatrix} + \boxed{-3} \begin{bmatrix} 2 \\ 5 \\ 1 \end{bmatrix} + \boxed{2} \begin{bmatrix} 3 \\ 1 \\ 2 \end{bmatrix}$$

This representation is called a Linear Combination of the Columns of Matrix A.

**Part b:** To represent the rows of matrix A as a linear combination, we need to do the same steps as we did in part a.

**Step 1:** Put each row of matrix A individually:

$$[1 \quad 2 \quad 3] \quad [0 \quad 5 \quad 1] \quad [2 \quad 1 \quad 2]$$

**Step 2:** Separate each one of them by addition sign.

$$[1 \quad 2 \quad 3] \; + \; [0 \quad 5 \quad 1] \; + \; [2 \quad 1 \quad 2]$$

**Step 3:** Put 1 box in front of each one.

$$\boxed{\phantom{0}}[1 \quad 2 \quad 3] + \boxed{\phantom{0}}[0 \quad 5 \quad 1] + \boxed{\phantom{0}}[2 \quad 1 \quad 2]$$

**Step 4:** Write random number for each row in each box because it is not mentioned any number for any row in the question.

$$\boxed{4}[1 \quad 2 \quad 3] + \boxed{-8}[0 \quad 5 \quad 1] + \boxed{11}[2 \quad 1 \quad 2]$$

This representation is called a Linear Combination of the Rows of Matrix A.

**Definition 1.3.1** Suppose matrix A has an order (size) $n \times m$, and matrix B has an order (size) $m \times k$, then the number of columns of matrix A equals to the number of rows of matrix B. If we assume that the usual multiplication of matrix A by matrix B equals to a new matrix C such that $A \cdot B = C$, then the matrix C has an order (size) $n \times k$. (i.e. If A has an order (size) $3 \times 7$, and B has an order(size) $7 \times 10$, then C has an order(size) $3 \times 10$).

22                                                                M. Kaabar**Example 1.3.3** Given the following matrix A and matrix B:

$$A = \begin{bmatrix} 1 & 2 \\ 3 & 6 \\ 1 & 1 \end{bmatrix}, \quad B = \begin{bmatrix} 1 & 1 & 1 \\ 0 & 2 & 1 \end{bmatrix}$$

   a) Find AB such that each column of AB is a linear combination of the columns of A.
   b) Find AB such that each row of AB is a linear combination of the rows of B.
   c) Find AB using the usual matrix multiplication.
   d) If AB = C, then find $c_{23}$.
   e) Find BA such that each column of BA is a linear combination of the columns of B.
   f) Find BA such that each row of BA is a linear combination of the rows of A.

**Solution: Part a:** Since A has an order (size) $3 \times 2$, and B has an order(size) $2 \times 3$, then AB will have an order(size) $3 \times 3$ according to Definition 1.3.1.

**Step 1:** 1st column of AB:

$$1 * \begin{bmatrix} 1 \\ 3 \\ 1 \end{bmatrix} + 0 * \begin{bmatrix} 2 \\ 6 \\ 1 \end{bmatrix} = \begin{bmatrix} 1 \\ 3 \\ 1 \end{bmatrix}$$

**Step 2:** 2nd column of AB:

$$1 * \begin{bmatrix} 1 \\ 3 \\ 1 \end{bmatrix} + 2 * \begin{bmatrix} 2 \\ 6 \\ 1 \end{bmatrix} = \begin{bmatrix} 5 \\ 15 \\ 3 \end{bmatrix}$$

**Step 3:** 3rd column of AB:

$$1 * \begin{bmatrix} 1 \\ 3 \\ 1 \end{bmatrix} + 1 * \begin{bmatrix} 2 \\ 6 \\ 1 \end{bmatrix} = \begin{bmatrix} 3 \\ 9 \\ 2 \end{bmatrix}$$

Hence, AB = $\begin{bmatrix} 1 & 5 & 3 \\ 3 & 15 & 9 \\ 1 & 3 & 2 \end{bmatrix}$

**Part b:** As we did in part a but the difference here is rows instead of columns.

**Step 1:** 1st row of AB:

$$1 * [1 \quad 1 \quad 1] + 2 * [0 \quad 2 \quad 1] = [1 \quad 5 \quad 3]$$

**Step 2:** 2nd row of AB:

$$3 * [1 \quad 1 \quad 1] + 6 * [0 \quad 2 \quad 1] = [3 \quad 15 \quad 9]$$

**Step 3:** 3rd row of AB:

$$1 * [1 \quad 1 \quad 1] + 1 * [0 \quad 2 \quad 1] = [1 \quad 3 \quad 2]$$

Hence, AB = $\begin{bmatrix} 1 & 5 & 3 \\ 3 & 15 & 9 \\ 1 & 3 & 2 \end{bmatrix}$

**Part c:** Here we use the usual matrix multiplication.

AB = $\begin{bmatrix} 1 & 2 \\ 3 & 6 \\ 1 & 1 \end{bmatrix} * \begin{bmatrix} 1 & 1 & 1 \\ 0 & 2 & 1 \end{bmatrix}$

$$AB = \begin{bmatrix} 1*1+2*0 & 1*1+2*2 & 1*1+2*1 \\ 3*1+6*0 & 3*1+6*2 & 3*1+6*1 \\ 1*1+1*0 & 1*1+1*2 & 1*1+1*1 \end{bmatrix}$$

Hence, $AB = \begin{bmatrix} 1 & 5 & 3 \\ 3 & 15 & 9 \\ 1 & 3 & 2 \end{bmatrix}$

**Part d:** Since AB = C, then $c_{23}$ means that we need to find the number that is located in 2nd row and 3rd column.

$$C = \begin{bmatrix} c_{11} & c_{12} & c_{13} \\ c_{21} & c_{22} & c_{23} \\ c_{31} & c_{32} & c_{33} \end{bmatrix} = \begin{bmatrix} 1 & 5 & 3 \\ 3 & 15 & 9 \\ 1 & 3 & 2 \end{bmatrix}$$

Now, let's explain each element of matrix C.

$c_{11}$ means that the number that is located in 1st row and 1st column. $c_{11} = 1$.

$c_{12}$ means that the number that is located in 1st row and 2nd column. $c_{12} = 5$.

$c_{13}$ means that the number that is located in 1st row and 3rd column. $c_{13} = 3$.

$c_{21}$ means that the number that is located in 2nd row and 1st column. $c_{21} = 3$.

$c_{22}$ means that the number that is located in 2nd row and 2nd column. $c_{22} = 15$.

$c_{23}$ means that the number that is located in 2nd row and 3rd column. $c_{23} = 9$.

$c_{31}$ means that the number that is located in 3rd row and 1st column. $c_{31} = 1$.

$c_{32}$ means that the number that is located in 3rd row and 2nd column. $c_{32} = 3$.

$c_{33}$ means that the number that is located in 3rd row and 3rd column. $c_{33} = 2$.

Hence, $c_{23} = 9$.

**Part e:** As we did in part a. Since B has an order(size) $2 \times 3$, and A has an order (size) $3 \times 2$, then BA will have an order(size) $2 \times 2$ according to Definition 1.3.1.

**Step 1:** $1^{st}$ column of BA:

$$1 * \begin{bmatrix} 1 \\ 0 \end{bmatrix} + 3 * \begin{bmatrix} 1 \\ 2 \end{bmatrix} + 1 * \begin{bmatrix} 1 \\ 1 \end{bmatrix} = \begin{bmatrix} 5 \\ 7 \end{bmatrix}$$

**Step 2:** $2^{nd}$ column of BA:

$$2 * \begin{bmatrix} 1 \\ 0 \end{bmatrix} + 6 * \begin{bmatrix} 1 \\ 2 \end{bmatrix} + 1 * \begin{bmatrix} 1 \\ 1 \end{bmatrix} = \begin{bmatrix} 9 \\ 13 \end{bmatrix}$$

Hence, $BA = \begin{bmatrix} 5 & 9 \\ 7 & 13 \end{bmatrix}$

**Part f:** As we did in part b.

**Step 1:** $1^{st}$ row of AB:

$$1 * [1 \quad 2] + 1 * [3 \quad 6] + 1 * [1 \quad 1] = [5 \quad 9]$$

**Step 2:** $2^{nd}$ row of AB:

$$0 * [1 \quad 2] + 2 * [3 \quad 6] + 1 * [1 \quad 1] = [7 \quad 13]$$

Hence, $BA = \begin{bmatrix} 5 & 9 \\ 7 & 13 \end{bmatrix}$

**Result 1.3.1** In general, matrix multiplication is not commutative (i.e. AB is not necessarily equal to BA).

**Example 1.3.4** Given the following matrix A and matrix B:

$A = \begin{bmatrix} 1 & 0 \\ 0 & 0 \end{bmatrix}$ , $B = \begin{bmatrix} 0 & 0 \\ 1 & 0 \end{bmatrix}$

a) Find AB.
b) Find BA.

**Solution: Part a:** Using the usual matrix multiplication, we obtain:

$AB = \begin{bmatrix} 1 & 0 \\ 0 & 0 \end{bmatrix} \begin{bmatrix} 0 & 0 \\ 1 & 0 \end{bmatrix} = \begin{bmatrix} 0 & 0 \\ 0 & 0 \end{bmatrix}$

**Part b:** Using the usual matrix multiplication, we obtain:

$BA = \begin{bmatrix} 0 & 0 \\ 1 & 0 \end{bmatrix} \begin{bmatrix} 1 & 0 \\ 0 & 0 \end{bmatrix} = \begin{bmatrix} 0 & 0 \\ 1 & 0 \end{bmatrix}$

**Result 1.3.2** It is possible that the product of two non-zero matrices is a zero matrix. However, this is not true for product of real numbers.

**Example 1.3.5** Given the following matrix A, matrix B and matrix D:

$A = \begin{bmatrix} 1 & 0 \\ 0 & 0 \end{bmatrix}$ , $B = \begin{bmatrix} 0 & 0 \\ 1 & 0 \end{bmatrix}$ , $D = \begin{bmatrix} 0 & 0 \\ 0 & 1 \end{bmatrix}$

a) Find AB.
b) Find AD.

**Solution: Part a:** Using the usual matrix multiplication, we obtain:

$AB = \begin{bmatrix} 1 & 0 \\ 0 & 0 \end{bmatrix} \begin{bmatrix} 0 & 0 \\ 1 & 0 \end{bmatrix} = \begin{bmatrix} 0 & 0 \\ 0 & 0 \end{bmatrix}$ zero matrix.

**Part b:** Using the usual matrix multiplication, we obtain:

$AD = \begin{bmatrix} 1 & 0 \\ 0 & 0 \end{bmatrix} \begin{bmatrix} 0 & 0 \\ 0 & 1 \end{bmatrix} = \begin{bmatrix} 0 & 0 \\ 0 & 0 \end{bmatrix}$ zero matrix.

**Result 1.3.3** In general, If $AB = AD$ and A is not a zero matrix, then it is possible that $B \neq D$.

**Definition 1.3.2** Suppose matrix A has an order (size) $n \times m$, A is a zero matrix if each number of A is a zero number.

**Example 1.3.5** Given the following system of linear equations:
$$\begin{cases} x_1 + x_2 + x_3 = 3 \\ -x_1 - x_2 + 2x_3 = 0 \\ x_2 + 4x_3 = 5 \end{cases}$$
Write the above system in matrix-form.

**Solution:** We write the above system in matrix-form as follows:

(Coefficient Matrix) · (Variable Column) = (Constant Column)

CX = A where C is a Coefficient Matrix, X is a Variable Column, and A is a Constant Column.

$$\begin{bmatrix} 1 & 1 & 1 \\ -1 & -1 & 2 \\ 0 & 1 & 4 \end{bmatrix} \begin{bmatrix} x_1 \\ x_2 \\ x_3 \end{bmatrix} = \begin{bmatrix} 3 \\ 0 \\ 5 \end{bmatrix}$$

The above matrix-form means the following:

$$x_1 * \begin{bmatrix} 1 \\ -1 \\ 0 \end{bmatrix} + x_2 * \begin{bmatrix} 1 \\ -1 \\ 1 \end{bmatrix} + x_3 * \begin{bmatrix} 1 \\ 2 \\ 4 \end{bmatrix}$$

Hence, we obtain:

$$\begin{bmatrix} x_1 + x_2 + x_3 = 3 \\ -x_1 - x_2 + 2x_3 = 0 \\ x_2 + 4x_3 = 5 \end{bmatrix} = \begin{bmatrix} 3 \\ 0 \\ 5 \end{bmatrix}$$

Now, let $x_1 = 1, x_2 = 1,$ and $x_3 = 1$.
$$\begin{bmatrix} 1 & 1 & 1 \\ -1 & -1 & 2 \\ 0 & 1 & 4 \end{bmatrix} \begin{bmatrix} 1 \\ 1 \\ 1 \end{bmatrix} = \begin{bmatrix} 3 \\ 0 \\ 5 \end{bmatrix}$$

**Result 1.3.4** Let CX = A be a system of linear equations
$$C \begin{bmatrix} x_1 \\ x_2 \\ \cdot \\ \cdot \\ x_n \end{bmatrix} = \begin{bmatrix} a_1 \\ a_2 \\ \cdot \\ \cdot \\ a_n \end{bmatrix}$$
Then, $x_1 = b_1, x_2 = b_2, \ldots, x_n = b_n$ is a solution to the system if and only if $b_1 C_1 + b_2 C_2 + \cdots + b_n C_m = A$, $C_1$ is the 1st column of C, $C_2$ is the 2nd column of C, …. , $C_m$ is $n^{th}$ column of C.

**Example 1.3.6** Given the following system of linear equations:
$$\begin{cases} x_1 + x_2 = 2 \\ 2x_1 + 2x_2 = 4 \end{cases}$$
Write the above system in matrix-form.
**Solution:** We write the above system in matrix-form CX = A as follows:
$$C = \begin{bmatrix} 1 & 1 \\ 2 & 2 \end{bmatrix}, \quad X = \begin{bmatrix} x1 \\ x2 \end{bmatrix}, \quad A = \begin{bmatrix} 2 \\ 4 \end{bmatrix}$$
Hence, CX = A:
$$\begin{bmatrix} 1 & 1 \\ 2 & 2 \end{bmatrix} \begin{bmatrix} x_1 \\ x_2 \end{bmatrix} = \begin{bmatrix} 2 \\ 4 \end{bmatrix}$$
$$x_1 * \begin{bmatrix} 1 \\ 2 \end{bmatrix} + x_2 * \begin{bmatrix} 1 \\ 2 \end{bmatrix} = \begin{bmatrix} 2 \\ 4 \end{bmatrix}$$
Now, we need to choose values for $x_1$ and $x_2$ such that these values satisfy that above matrix-form.
First, let's try $x_1 = 3$ and $x_2 = 4$.
$$3 * \begin{bmatrix} 1 \\ 2 \end{bmatrix} + 4 * \begin{bmatrix} 1 \\ 2 \end{bmatrix} = \begin{bmatrix} 7 \\ 14 \end{bmatrix} \neq \begin{bmatrix} 2 \\ 4 \end{bmatrix}$$

Therefore, $x_1 = 3$, $x_2 = 4$ is not solution, and our assumption is wrong.
Now, let's try $x_1 = 0$ and $x_2 = 2$.
$0 * \begin{bmatrix} 1 \\ 2 \end{bmatrix} + 2 * \begin{bmatrix} 1 \\ 2 \end{bmatrix} = \begin{bmatrix} 2 \\ 4 \end{bmatrix} = \begin{bmatrix} 2 \\ 4 \end{bmatrix}$
Therefore, $x_1 = 0$, $x_2 = 2$ is solution.

**Result 1.3.5** Let $CX = A$ be a system of linear equations. The constant column A can be written uniquely as a linear combination of the columns of C.

# 1.4 Square Matrix

In this section, we discuss some examples and facts of square matrix and identity matrix.

**Definition 1.4.1** In general, let A is a matrix, and n be a positive integer such that $A^n$ is defined only if number of rows of A = number of columns of A. (i.e. Let A is a matrix with a size (order) $\times m$, $A^2 = A \cdot A$ is defined only if n = m). This matrix is called a Square-Matrix.
The following is some facts regarding square-matrix and set of real numbers $\mathbb{R}$.

**Fact 1.4.1** Let A is a matrix, $A^{-n}$ is not equal to $\frac{1}{A^n}$.

**Fact 1.4.2** Let A is a matrix, $A^{-n}$ has a specific meaning if A is a square-matrix and A is invertible (nonsingular). We will talk about it in section 1.5.

**Fact 1.4.3** A set of real numbers $\mathbb{R}$ has a multiplicative identity equals to 1. (i.e. $a \cdot 1 = 1 \cdot a = a$).

**Fact 1.4.4** A set of real numbers $\mathbb{R}$ has an additive identity equals to 0. (i.e. $a + 0 = 0 + a = a$).

**Fact 1.4.5** $\mathbb{R}^{2\times 2} = M_2(\mathbb{R})$ = set of all $2 \times 2$ matrices.

Now, we give some helpful notations:
$\mathbb{Z}$ : Set of all integers.
$\mathbb{Q}$ : Set of all rational numbers.
$\mathbb{R}$ : Set of all real number.
$\mathbb{N}$ : Set of all natural numbers.

The following figure 1.4.1 represents three main sets: $\mathbb{Z}, \mathbb{Q}$ and $\mathbb{R}$.

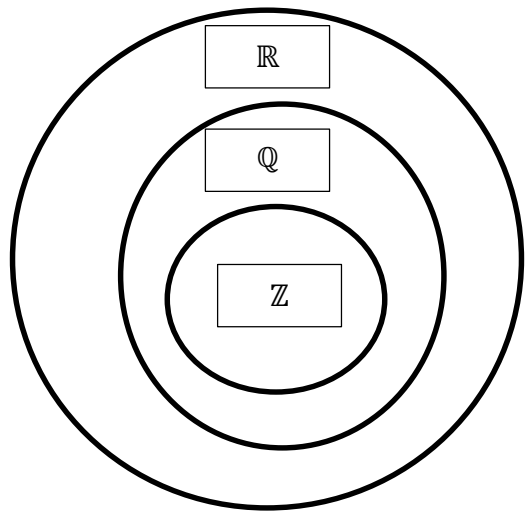

Figure 1.4.1: Representation of Three Sets of Numbers

The above figure shows that the largest set is $\mathbb{R}$ , and the smallest set is $\mathbb{Z}$. In addition, we know that a rational number is an integer divided by a non-zero integer.

**Fact 1.4.6** $3 \in \mathbb{R}$ means that $3$ is an element of the set $\mathbb{R}$. ("$\in$" is a mathematical symbol that means "belong to").

**Fact 1.4.7** $\sqrt{2} \in \mathbb{R}$ means that $\sqrt{2}$ is an element of the set $\mathbb{R}$. ("$\in$" is a mathematical symbol that means "belong to").

**Fact 1.4.8** $\frac{1}{2} \notin \mathbb{Z}$ means that $\frac{1}{2}$ is not an element of the set $\mathbb{Z}$. ("$\notin$" is a mathematical symbol that means "does not belong to").

**Fact 1.4.9** $A \in M_2(\mathbb{R})$ means that A is a $2 \times 2$ matrix.

**Fact 1.4.10** $\mathbb{R}^{n \times n} = M_n(\mathbb{R}) =$ set of all $n \times n$ matrices. (i.e. $\mathbb{R}^{3 \times 3} = M_3(\mathbb{R}) =$ set of all $3 \times 3$ matrices).

After learning all above facts, let's test our knowledge by asking the following question.

**Question 1.4.1** Does $M_n(\mathbb{R})$ has a multiplicative identity?

**Solution:** Yes; there is a matrix called $I_n$ where $I_n \in M_n(\mathbb{R})$ such that $AI_n = I_n A = A$ for every $A \in M_n(\mathbb{R})$.

**Example 1.4.1** Given $M_2(\mathbb{R}) = \mathbb{R}^{2 \times 2}$. What is $I_2$?

**Solution:** We need to find the multiplicative identity for all $2 \times 2$ matrices. This means that $AI_2 = I_2 A = A$ for every $A \in M_2(\mathbb{R})$.

$I_2 = \begin{bmatrix} 1 & 0 \\ 0 & 1 \end{bmatrix}$ and $A = \begin{bmatrix} a_{11} & a_{12} \\ a_{21} & a_{22} \end{bmatrix}$

Now, let's check if the multiplicative identity is right:

$AI_2 = \begin{bmatrix} a_{11} & a_{12} \\ a_{21} & a_{22} \end{bmatrix} \begin{bmatrix} 1 & 0 \\ 0 & 1 \end{bmatrix} = \begin{bmatrix} a_{11} & a_{12} \\ a_{21} & a_{22} \end{bmatrix} = A$

$I_2 A = \begin{bmatrix} 1 & 0 \\ 0 & 1 \end{bmatrix} \begin{bmatrix} a_{11} & a_{12} \\ a_{21} & a_{22} \end{bmatrix} = \begin{bmatrix} a_{11} & a_{12} \\ a_{21} & a_{22} \end{bmatrix} = A$

Hence, the multiplicative identity for all $2 \times 2$ matrices is $I_2 = \begin{bmatrix} 1 & 0 \\ 0 & 1 \end{bmatrix}$.

**Example 1.4.2** Given $M_3(\mathbb{R}) = \mathbb{R}^{3 \times 3}$. What is $I_3$?

**Solution:** We need to find the multiplicative identity for all $3 \times 3$ matrices. This means that $AI_3 = I_3A = A$ for every $A \in M_3(\mathbb{R})$.

$$I_3 = \begin{bmatrix} 1 & 0 & 0 \\ 0 & 1 & 0 \\ 0 & 0 & 1 \end{bmatrix} \text{ and } A = \begin{bmatrix} a_{11} & a_{12} & a_{13} \\ a_{21} & a_{22} & a_{23} \\ a_{31} & a_{32} & a_{33} \end{bmatrix}$$

Now, let's check if the multiplicative identity is right:

$$AI_3 = \begin{bmatrix} a_{11} & a_{12} & a_{13} \\ a_{21} & a_{22} & a_{23} \\ a_{31} & a_{32} & a_{33} \end{bmatrix} \begin{bmatrix} 1 & 0 & 0 \\ 0 & 1 & 0 \\ 0 & 0 & 1 \end{bmatrix} = \begin{bmatrix} a_{11} & a_{12} & a_{13} \\ a_{21} & a_{22} & a_{23} \\ a_{31} & a_{32} & a_{33} \end{bmatrix} = A$$

$$I_3 A = \begin{bmatrix} 1 & 0 & 0 \\ 0 & 1 & 0 \\ 0 & 0 & 1 \end{bmatrix} \begin{bmatrix} a_{11} & a_{12} & a_{13} \\ a_{21} & a_{22} & a_{23} \\ a_{31} & a_{32} & a_{33} \end{bmatrix} = \begin{bmatrix} a_{11} & a_{12} & a_{13} \\ a_{21} & a_{22} & a_{23} \\ a_{31} & a_{32} & a_{33} \end{bmatrix} = A$$

Hence, the multiplicative identity for all $3 \times 3$ matrices is $I_3 = \begin{bmatrix} 1 & 0 & 0 \\ 0 & 1 & 0 \\ 0 & 0 & 1 \end{bmatrix}$.

**Result 1.4.1** Assume we have a square-matrix with $n \times n$ size and the main-diagonal is $a_{11}, a_{22}, a_{33}, \dots, a_{nn}$

$$\begin{bmatrix} a_{11} & \cdots & & a_{1n} \\ \vdots & a_{22} & \cdots & \vdots \\ & \cdots & \cdots & \\ a_{n1} & & \cdots & a_{nn} \end{bmatrix}$$

Then, the multiplicative identity for the above $n \times n$ matrix is as follows:

$$I_n = \begin{bmatrix} 1 & \cdots & & 0 \\ \vdots & 1 & 0 & \vdots \\ & 0 & 1 & \\ 0 & & \cdots & 1 \end{bmatrix}$$

All ones on the main-diagonal and zeros elsewhere.

## 1.5 Inverse Square Matrix

In this section, we discuss the concept of inverse square matrix, and we give some examples of elementary matrix.

**Example 1.5.1** Given the following matrices with the following row-operations steps:

$$Z = \begin{bmatrix} 1 & 2 & 3 & 4 \\ 0 & 1 & -1 & 2 \\ 0 & 1 & 1 & 3 \end{bmatrix}, \; D = \begin{bmatrix} -2 & -4 & -6 & -8 \\ 0 & 1 & -1 & 2 \\ 0 & 1 & 1 & 3 \end{bmatrix}$$

$$\begin{bmatrix} 1 & 2 & 3 & 4 \\ 0 & 1 & -1 & 2 \\ 0 & 1 & 1 & 3 \end{bmatrix}$$

$-2R_1 \rightarrow$

$$\begin{bmatrix} -2 & -4 & -6 & -8 \\ 0 & 1 & -1 & 2 \\ 0 & 1 & 1 & 3 \end{bmatrix}$$

$2R_1 + R_3 \dashrightarrow R_3$

$-3R_1 + R_2 \dashrightarrow R_2$

Matrix B    $R_3 \leftrightarrow R_1$    Matrix C

a) Find a matrix D such that DZ=W.
b) Find a matrix K such that KW=Z.
c) Find Elementary Matrices, $F_1$, $F_2$, $F_3$, ..., $F_m$ such that $F_1 F_2 \ldots F_m Z = C$.

   d) Find Elementary Matrices, $L_1, L_2, \ldots, L_n$ such that $L_1 L_2 \ldots L_n B = Z$.

   e) Find a matrix S such that SC=B.

   f) Find a matrix X such that XB=C.

**Solution: Part a:** Since Z is $3 \times 4$ matrix and we multiply from left by -2, then D must be a square matrix. Hence, D must be $3 \times 3$.
We use the multiplicative identity, and we multiply it by -2 from left as follows:

$$I_3 = \begin{bmatrix} 1 & 0 & 0 \\ 0 & 1 & 0 \\ 0 & 0 & 1 \end{bmatrix} \text{-2}R_1 \rightarrow \begin{bmatrix} -2 & 0 & 0 \\ 0 & 1 & 0 \\ 0 & 0 & 1 \end{bmatrix} = D$$

D is called Elementary Matrix of Type I.

Thus, $\begin{bmatrix} -2 & 0 & 0 \\ 0 & 1 & 0 \\ 0 & 0 & 1 \end{bmatrix} Z = W$

**Part b:** Since W is $3 \times 4$ matrix and we multiply from right by $-\frac{1}{2}$ because the inverse of -2 is $-\frac{1}{2}$, then K must be a square matrix. Hence, K must be $3 \times 3$.
We use the multiplicative identity, and we multiply it by $-\frac{1}{2}$ from left as follows:

$$I_3 = \begin{bmatrix} 1 & 0 & 0 \\ 0 & 1 & 0 \\ 0 & 0 & 1 \end{bmatrix} -\frac{1}{2}R_1 \rightarrow \begin{bmatrix} -\frac{1}{2} & 0 & 0 \\ 0 & 1 & 0 \\ 0 & 0 & 1 \end{bmatrix} = K$$

K is also called Elementary Matrix of Type I.

Thus, $\begin{bmatrix} -\frac{1}{2} & 0 & 0 \\ 0 & 1 & 0 \\ 0 & 0 & 1 \end{bmatrix} W = Z$

**Result 1.5.1** Each row-operation corresponds to one and only one elementary matrix.

**Part c:** According to the given four row-operations steps, we need to find four elementary matrices from Z to C which means m=4 because we want to go 4 steps forward from Z to C. From part a, we already have the following:

$$I_3 = \begin{bmatrix} 1 & 0 & 0 \\ 0 & 1 & 0 \\ 0 & 0 & 1 \end{bmatrix} \text{-2R}_1 \rightarrow \begin{bmatrix} -2 & 0 & 0 \\ 0 & 1 & 0 \\ 0 & 0 & 1 \end{bmatrix} \text{ (Step One)}$$

Hence, $\begin{bmatrix} -2 & 0 & 0 \\ 0 & 1 & 0 \\ 0 & 0 & 1 \end{bmatrix} Z = C$

$$I_3 = \begin{bmatrix} 1 & 0 & 0 \\ 0 & 1 & 0 \\ 0 & 0 & 1 \end{bmatrix} 2R_1 + R_3 \dashrightarrow R_3 \begin{bmatrix} 1 & 0 & 0 \\ 0 & 1 & 0 \\ 2 & 0 & 1 \end{bmatrix} \text{ (Step Two)}$$

Hence, $\begin{bmatrix} 1 & 0 & 0 \\ 0 & 1 & 0 \\ 2 & 0 & 1 \end{bmatrix} \begin{bmatrix} -2 & 0 & 0 \\ 0 & 1 & 0 \\ 0 & 0 & 1 \end{bmatrix} Z = C$

$$I_3 = \begin{bmatrix} 1 & 0 & 0 \\ 0 & 1 & 0 \\ 0 & 0 & 1 \end{bmatrix} \text{-3R}_1 + R_2 \dashrightarrow R_2 \begin{bmatrix} 1 & 0 & 0 \\ -3 & 1 & 0 \\ 0 & 0 & 1 \end{bmatrix} \text{ (Step Three)}$$

Hence, $\begin{bmatrix} 1 & 0 & 0 \\ -3 & 1 & 0 \\ 0 & 0 & 1 \end{bmatrix} \begin{bmatrix} 1 & 0 & 0 \\ 0 & 1 & 0 \\ 2 & 0 & 1 \end{bmatrix} \begin{bmatrix} -2 & 0 & 0 \\ 0 & 1 & 0 \\ 0 & 0 & 1 \end{bmatrix} Z = C$

$$I_3 = \begin{bmatrix} 1 & 0 & 0 \\ 0 & 1 & 0 \\ 0 & 0 & 1 \end{bmatrix} R_3 \leftrightarrow R_1 \begin{bmatrix} 0 & 0 & 1 \\ 0 & 1 & 0 \\ 1 & 0 & 0 \end{bmatrix} \text{(Step Four)}$$

Hence, we obtain the following four elementary matrices:

$$\begin{bmatrix} 0 & 0 & 1 \\ 0 & 1 & 0 \\ 1 & 0 & 0 \end{bmatrix} \begin{bmatrix} 1 & 0 & 0 \\ -3 & 1 & 0 \\ 0 & 0 & 1 \end{bmatrix} \begin{bmatrix} 1 & 0 & 0 \\ 0 & 1 & 0 \\ 2 & 0 & 1 \end{bmatrix} \begin{bmatrix} -2 & 0 & 0 \\ 0 & 1 & 0 \\ 0 & 0 & 1 \end{bmatrix} Z = C$$

**Part d:** According to the given three row-operations steps from B to Z, we need to find three elementary matrices which means n=3 because we want to go 3 steps backward from B to Z. Backward steps mean that we need to do inverse steps (i.e. The inverse of -2$R_1$ is $-\frac{1}{2} R_1$ because it is row-multiplication step). We start from the third step as follows:

The inverse of -3$R_1$+$R_2$--→ $R_2$ is 3$R_1$+$R_2$--→ $R_2$ because it is row-addition step.

$$I_3 = \begin{bmatrix} 1 & 0 & 0 \\ 0 & 1 & 0 \\ 0 & 0 & 1 \end{bmatrix} 3R_1+R_2 \dashrightarrow R_2 \begin{bmatrix} 1 & 0 & 0 \\ 3 & 1 & 0 \\ 0 & 0 & 1 \end{bmatrix} \text{(Step Three)}$$

Hence, $\begin{bmatrix} 1 & 0 & 0 \\ 3 & 1 & 0 \\ 0 & 0 & 1 \end{bmatrix} B = Z$

The inverse of 2$R_1$+$R_3$--→ $R_3$ is -2$R_1$+$R_3$--→ $R_3$ because it is row-addition step.

$$I_3 = \begin{bmatrix} 1 & 0 & 0 \\ 0 & 1 & 0 \\ 0 & 0 & 1 \end{bmatrix} \text{-2}R_1+R_3 \dashrightarrow R_3 \begin{bmatrix} 1 & 0 & 0 \\ 0 & 1 & 0 \\ -2 & 0 & 1 \end{bmatrix} \text{(Step Two)}$$

Hence, $\begin{bmatrix} 1 & 0 & 0 \\ 0 & 1 & 0 \\ -2 & 0 & 1 \end{bmatrix} \begin{bmatrix} 1 & 0 & 0 \\ 3 & 1 & 0 \\ 0 & 0 & 1 \end{bmatrix} B = Z$

The inverse of $-2R_1$ is $-\frac{1}{2} R_1$ because it is row-multiplication step

$I_3 = \begin{bmatrix} 1 & 0 & 0 \\ 0 & 1 & 0 \\ 0 & 0 & 1 \end{bmatrix} -\frac{1}{2}R_1 \rightarrow \begin{bmatrix} -\frac{1}{2} & 0 & 0 \\ 0 & 1 & 0 \\ 0 & 0 & 1 \end{bmatrix}$ (Step One)

Hence, we obtain the following three elementary matrices:

$\begin{bmatrix} -\frac{1}{2} & 0 & 0 \\ 0 & 1 & 0 \\ 0 & 0 & 1 \end{bmatrix} \begin{bmatrix} 1 & 0 & 0 \\ 0 & 1 & 0 \\ -2 & 0 & 1 \end{bmatrix} \begin{bmatrix} 1 & 0 & 0 \\ 3 & 1 & 0 \\ 0 & 0 & 1 \end{bmatrix} B = Z$

**Part e:** According to the given last row-operations step from C to B, we need to find one elementary matrix because we want to go 1 step backward from C to B. The inverse of $R_3 \leftrightarrow R_1$ is $R_1 \leftrightarrow R_3$ which is the same as $R_3 \leftrightarrow R_1$.

$I_3 = \begin{bmatrix} 1 & 0 & 0 \\ 0 & 1 & 0 \\ 0 & 0 & 1 \end{bmatrix} R_3 \leftrightarrow R_1 \begin{bmatrix} 0 & 0 & 1 \\ 0 & 1 & 0 \\ 1 & 0 & 0 \end{bmatrix}$

Hence, we obtain: $\begin{bmatrix} 0 & 0 & 1 \\ 0 & 1 & 0 \\ 1 & 0 & 0 \end{bmatrix} C = B$. This means that $S = \begin{bmatrix} 0 & 0 & 1 \\ 0 & 1 & 0 \\ 1 & 0 & 0 \end{bmatrix}$.

**Part f:** According to the given last row-operations step from B to C, we need to find one elementary matrix because we want to go 1 step forward from B to C. $I_3 =$
$$\begin{bmatrix} 1 & 0 & 0 \\ 0 & 1 & 0 \\ 0 & 0 & 1 \end{bmatrix} R_3 \leftrightarrow R_1 \begin{bmatrix} 0 & 0 & 1 \\ 0 & 1 & 0 \\ 1 & 0 & 0 \end{bmatrix}$$

Hence, we obtain: $\begin{bmatrix} 0 & 0 & 1 \\ 0 & 1 & 0 \\ 1 & 0 & 0 \end{bmatrix} B = C$. This means that

$X = \begin{bmatrix} 0 & 0 & 1 \\ 0 & 1 & 0 \\ 1 & 0 & 0 \end{bmatrix}$.

**Example 1.5.2** Given the following matrix:
$$A = \begin{bmatrix} 2 & 1 \\ 4 & 0 \end{bmatrix}$$
Find $A^{-1}$ if possible.

**Solution:** Finding $A^{-1}$ if possible means that we need to find a possible inverse matrix called $A^{-1}$ such that $AA^{-1} = I_2 = A^{-1}A$. To find this possible matrix $A^{-1}$, we need to do the following steps:

**Step 1:** Write A and $I_2$ in the following form:
$$\left( \text{Matrix A} \,\middle|\, \text{Identity Matrix } I_2 \right)$$

$$\begin{pmatrix} 2 & 1 & 1 & 0 \\ 4 & 0 & 0 & 1 \end{pmatrix}$$

**Step 2:** Do some Row-Operations until you get the following:
$$\left( \text{Do you see } I_2? \text{ If yes, then } X = A^{-1} \,\middle|\, X \right)$$
$$\left( \text{Do you see } I_2? \text{ If no, then A has no inverse matrix} \,\middle|\, X \right)$$

Now, let's do the above step until we get the Completely-Reduced-Echelon matrix.

$$\begin{pmatrix} 2 & 1 & | & 1 & 0 \\ 4 & 0 & | & 0 & 1 \end{pmatrix} \xrightarrow{-\frac{1}{2}R_1} \begin{pmatrix} 1 & \frac{1}{2} & | & \frac{1}{2} & 0 \\ 4 & 0 & | & 0 & 1 \end{pmatrix}$$

$$\begin{pmatrix} 1 & \frac{1}{2} & | & \frac{1}{2} & 0 \\ 4 & 0 & | & 0 & 1 \end{pmatrix} \xrightarrow{-4R_1+R_2 \dashrightarrow R_2} \begin{pmatrix} 1 & \frac{1}{2} & | & \frac{1}{2} & 0 \\ 0 & -2 & | & -2 & 1 \end{pmatrix}$$

$$\begin{pmatrix} 1 & \frac{1}{2} & | & \frac{1}{2} & 0 \\ 0 & -2 & | & -2 & 1 \end{pmatrix} \xrightarrow{-\frac{1}{2}R_2} \begin{pmatrix} 1 & \frac{1}{2} & | & \frac{1}{2} & 0 \\ 0 & 1 & | & 1 & -\frac{1}{2} \end{pmatrix}$$

$$\begin{pmatrix} 1 & \frac{1}{2} & | & \frac{1}{2} & 0 \\ 0 & 1 & | & 1 & -\frac{1}{2} \end{pmatrix} \xrightarrow{-\frac{1}{2}R_2+R_1 \dashrightarrow R_1} \begin{pmatrix} 1 & 0 & | & 0 & \frac{1}{4} \\ 0 & 1 & | & 1 & -\frac{1}{2} \end{pmatrix}$$

Since we got the Completely-Reduced-Echelon matrix which is the identity matrix $I_2$, then A has an inverse matrix which is $A^{-1}$.

Hence, $A^{-1} = \begin{bmatrix} 0 & \frac{1}{4} \\ 1 & -\frac{1}{2} \end{bmatrix}$

**Example 1.5.3** Given the following matrix:
$$A = \begin{bmatrix} 1 & 0 & 2 \\ 0 & 1 & 0 \\ 0 & -1 & 1 \end{bmatrix}$$
Find $A^{-1}$ if possible.

**Solution:** Finding $A^{-1}$ if possible means that we need to find a possible inverse matrix called $A^{-1}$ such that

$AA^{-1} = I_3 = A^{-1}A$. To find this possible matrix $A^{-1}$, we need to do the following steps:

**Step 1:** Write A and $I_3$ in the following form:

$$(\text{Matrix A} | \text{Identity Matrix } I_3)$$

$$\begin{pmatrix} 1 & 0 & 2 & | & 1 & 0 & 0 \\ 0 & 1 & 0 & | & 0 & 1 & 0 \\ 0 & -1 & 1 & | & 0 & 0 & 1 \end{pmatrix}$$

**Step 2:** Do some Row-Operations until you get the following:

$$(\text{Do you see } I_3 \text{? If yes, then } X = A^{-1} | X)$$

$$(\text{Do you see } I_3 \text{? If no, then A has no inverse matrix} | X)$$

Now, let's do the above step until we get the Completely-Reduced-Echelon matrix.

$$\begin{pmatrix} 1 & 0 & 2 & | & 1 & 0 & 0 \\ 0 & 1 & 0 & | & 0 & 1 & 0 \\ 0 & -1 & 1 & | & 0 & 0 & 1 \end{pmatrix} R_2 + R_3 \dashrightarrow R_3$$

$$\begin{pmatrix} 1 & 0 & 2 & | & 1 & 0 & 0 \\ 0 & 1 & 0 & | & 0 & 1 & 0 \\ 0 & 0 & 1 & | & 0 & 1 & 1 \end{pmatrix} \text{-}2R_3 + R_1 \dashrightarrow R_1$$

$$\begin{pmatrix} 1 & 0 & 0 & | & 1 & -2 & -2 \\ 0 & 1 & 0 & | & 0 & 1 & 0 \\ 0 & 0 & 1 & | & 0 & 1 & 1 \end{pmatrix}$$

Since we got the Completely-Reduced-Echelon matrix which is the identity matrix $I_3$ , then A has an inverse matrix which is A$^{-1}$.

Hence, A$^{-1}$ = $\begin{bmatrix} 1 & -2 & -2 \\ 0 & 1 & 0 \\ 0 & 1 & 1 \end{bmatrix}$

**Result 1.5.2** Given $n \times m$ matrix A and identity matrix $I_n$, and with some row-operations, we got $A'$ and $I_n'$ as follows: $(A | I_n)$ Row-Operations-$\rightarrow$ $(A' | I_n')$

Then, $I_n' A = A'$.

**Result 1.5.3** Given $n \times n$ matrix A and identity matrix $I_n$, and with some row-operations, we got $A^{-1}$ and $I_n$ as follows: $(A | I_n)$ Row-Operations-$\rightarrow$ $(I_n | A^{-1})$

Then, AA$^{-1}$ = $I_n$ = A$^{-1}$A.

**Fact 1.5.1** Given $n \times m$ matrix A and identity matrix $I_n$, $I_n$ is called a left identity for all $n \times m$ matrices such that $I_n A = A$. (i.e. Given $3 \times 5$ matrix A, then $I_3 A = A$).

**Fact 1.5.2** Given $n \times n$ matrix A and identity matrix $I_m$, $I_m$ is called a right identity for all $n \times n$ matrices such that $A I_m = A$. (i.e. Given $3 \times 5$ matrix A, then $A I_5 = A$).

**Result 1.5.4** Given $n \times n$ matrix A, A$^{-4}$ has a meaning if and only if A$^{-1}$ exists.

**Result 1.5.5** Given $n \times n$ matrix A, If A$^{-1}$ exists, then A$^{-4}$ = A$^{-1}$ × A$^{-1}$ × A$^{-1}$ × A$^{-1}$.

## 1.6 Transpose Matrix

In this section, we introduce the concept of transpose matrix, and we discuss some examples of symmetric and skew-symmetric matrices.

**Definition 1.6.1** Given $n \times m$ matrix A, $A^T$ is called A transpose and it is $m \times n$ matrix. We obtain $A^T$ from A by making columns of A rows or by making rows of A columns.

**Example 1.6.1** Given the following matrix:
$$A = \begin{bmatrix} 2 & 3 & 1 & 0 \\ 5 & 1 & 2 & 3 \\ 1 & 1 & 0 & 5 \end{bmatrix}$$
Find $A^T$.

**Solution:** According to definition 1.6.1, A is $3 \times 4$ matrix. Thus, $A^T$ should be $4 \times 3$. By making columns of A rows, we obtain the following:
$$A^T = \begin{bmatrix} 2 & 5 & 1 \\ 3 & 1 & 1 \\ 1 & 2 & 0 \\ 0 & 3 & 5 \end{bmatrix}$$

**Definition 1.6.2** Given $n \times m$ matrix A and $m \times n$ matrix $A^T$, then $AA^T$ is always defined, and it is $n \times n$ matrix. (i.e. Let $AA^T = D$, then $A_{n \times m} A^T_{m \times n} = D_{n \times n}$).

**Definition 1.6.3** Given $n \times m$ matrix A and $m \times n$ matrix $A^T$, then $A^T A$ is always defined, and it is $m \times m$ matrix. (i.e. Let $A^T A = F$, then $A^T_{m \times n} A_{n \times m} = F_{m \times m}$).

**Definition 1.6.4** Given $n \times n$ matrix A and $n \times n$ matrix $A^T$, then A is symmetric if $A^T = A$.

**Definition 1.6.5** Given $n \times n$ matrix A and $n \times n$ matrix $A^T$, then A is skew-symmetric if $A^T = -A$.

**Example 1.6.2** Given the following matrix:
$$A = \begin{bmatrix} 1 & 5 & 10 \\ 5 & 3 & 7 \\ 10 & 7 & 10 \end{bmatrix}$$
Show that A is symmetric.

**Solution:** According to definition 1.6.4, A is $3 \times 3$ matrix. Thus, $A^T$ should be $3 \times 3$. By making columns of A rows, we obtain the following:
$$A^T = \begin{bmatrix} 1 & 5 & 10 \\ 5 & 3 & 7 \\ 10 & 7 & 10 \end{bmatrix}$$
Since $A^T = A = \begin{bmatrix} 1 & 5 & 10 \\ 5 & 3 & 7 \\ 10 & 7 & 10 \end{bmatrix}$, then A is symmetric.

**Example 1.6.3** Given the following matrix:
$$A = \begin{bmatrix} 0 & 2 & 5 \\ -2 & 0 & 3 \\ -5 & -3 & 0 \end{bmatrix}$$
Show that A is skew-symmetric.

**Solution:** According to definition 1.6.5, A is $3 \times 3$ matrix. Thus, $A^T$ should be $3 \times 3$. By making columns of A rows, we obtain the following:
$$A^T = \begin{bmatrix} 0 & -2 & -5 \\ 2 & 0 & -3 \\ 5 & 3 & 0 \end{bmatrix}$$

Since $A^T = -A = \begin{bmatrix} 0 & -2 & -5 \\ 2 & 0 & -3 \\ 5 & 3 & 0 \end{bmatrix}$, then A is skew-symmetric.

**Fact 1.6.1** Given $n \times n$ matrix A, if A is skew-symmetric, then all numbers on the main diagonal of A are always zeros.

**Fact 1.6.2** Given $n \times m$ matrix A, then $(A^T)^T = A$.

**Fact 1.6.3** Given $n \times m$ matrix A and $m \times k$ matrix B, then $(AB^T)^T = B^T A^T$. (Warning: $(AB^T)^T \neq A^T B^T$).

**Fact 1.6.4** Given $n \times m$ matrices A and B, then $(A \pm B)^T = A^T \pm B^T$.

**Result 1.6.1** Given matrix A and constant $\alpha$, if A is symmetric, then $\alpha A$ is symmetric such that $(\alpha A)^T = \alpha A^T$.

**Result 1.6.2** Given matrix A and constant $\alpha$, if A is skew-symmetric, then $\alpha A$ is skew-symmetric such that $(\alpha A)^T = \alpha A^T$.

**Result 1.6.3** Let A be $n \times n$ matrix, there exists a symmetric matrix B and a skew-symmetric matrix C such that A is a linear combination of B and C. This means that those should be numbers $\alpha_1$ and $\alpha_2$ such that $A = \alpha_1 B + \alpha_2 C$.

**Proof of Result 1.6.3** We will show that A is a linear combination of B and C.

We assume that B is symmetric such that $B = A + A^T$.

$B^T = (A + A^T)^T = A^T + (A^T)^T = A^T + A = B$.

Now, we assume that C is skew-symmetric such that $C = A - A^T$.

$C^T = (A - A^T)^T = A^T - (A^T)^T = A^T - A = -(A - A^T) = -C$

By using algebra, we do the following:

$$\frac{1}{2}B + \frac{1}{2}C = \frac{1}{2}(A + A^T) + \frac{1}{2}(A - A^T)$$

$$= \frac{1}{2}A + \frac{1}{2}A^T + \frac{1}{2}A - \frac{1}{2}A^T = A.$$

Thus, A is a linear combination of B and C. ∎

**Example 1.6.4** Given the following matrix:

$$A = \begin{bmatrix} 2 & 1 & 4 \\ 3 & 0 & 1 \\ 5 & 6 & 7 \end{bmatrix}$$

Find symmetric matrix B and skew-symmetric matrix C such that $A = \alpha_1 B + \alpha_2 C$ for some numbers $\alpha_1$ and $\alpha_2$.

**Solution:** As we did in the above proof, we do the following:

$$B = A + A^T = \begin{bmatrix} 2 & 1 & 4 \\ 3 & 0 & 1 \\ 5 & 6 & 7 \end{bmatrix} + \begin{bmatrix} 2 & 3 & 5 \\ 1 & 0 & 6 \\ 4 & 1 & 7 \end{bmatrix} = \begin{bmatrix} 4 & 4 & 9 \\ 4 & 0 & 7 \\ 9 & 7 & 14 \end{bmatrix}$$

$$C = A - A^T = \begin{bmatrix} 2 & 1 & 4 \\ 3 & 0 & 1 \\ 5 & 6 & 7 \end{bmatrix} - \begin{bmatrix} 2 & 3 & 5 \\ 1 & 0 & 6 \\ 4 & 1 & 7 \end{bmatrix} = \begin{bmatrix} 0 & -2 & -1 \\ 2 & 0 & -5 \\ 1 & 5 & 0 \end{bmatrix}$$

Let $\alpha_1 = \alpha_2 = \frac{1}{2}$

Thus, $A = \frac{1}{2}B + \frac{1}{2}C = \frac{1}{2}\begin{bmatrix} 4 & 4 & 9 \\ 4 & 0 & 7 \\ 9 & 7 & 14 \end{bmatrix} + \frac{1}{2}\begin{bmatrix} 0 & -2 & -1 \\ 2 & 0 & -5 \\ 1 & 5 & 0 \end{bmatrix}$

## 1.7 Determinants

In this section, we introduce step by step for finding determinant of a certain matrix. In addition, we discuss some important properties such as invertible and non-invertible. In addition, we talk about the effect of row-operations on determinants.

**Definition 1.7.1** Determinant is a square matrix. Given $M_2(\mathbb{R}) = \mathbb{R}^{2\times 2} = \mathbb{R}_{2\times 2}$, let $A \in M_2(\mathbb{R})$ where A is $2 \times 2$ matrix, $A = \begin{bmatrix} a_{11} & a_{12} \\ a_{21} & a_{22} \end{bmatrix}$. The determinant of A is represented by det(A) or |A|.

Hence, $\det(A) = |A| = a_{11}a_{22} - a_{12}a_{21} \in \mathbb{R}$. (Warning: this definition works only for $2 \times 2$ matrices).

**Example 1.7.1** Given the following matrix:
$$A = \begin{bmatrix} 3 & 2 \\ 5 & 7 \end{bmatrix}$$
Find the determinant of A.

**Solution:** Using definition 1.7.1, we do the following:
$\det(A) = |A| = (3)(7) - (2)(5) = 21 - 10 = 11$.
Thus, the determinant of A is 11.

**Example 1.7.2** Given the following matrix:
$$A = \begin{bmatrix} 1 & 0 & 2 \\ 3 & 1 & -1 \\ 1 & 2 & 4 \end{bmatrix}$$
Find the determinant of A.

**Solution:** Since A is 3 × 3 matrix such that
$A \in M_3(\mathbb{R}) = \mathbb{R}^{3 \times 3}$, then we cannot use definition 1.7.1 because it is valid only for 2 × 2 matrices. Thus, we need to use the following method to find the determinant of A.

**Step 1:** Choose any row or any column. It is recommended to choose the one that has more zeros. In this example, we prefer to choose the second column or the first row. Let's choose the second column as follows:

$$A = \begin{bmatrix} 1 & 0 & 2 \\ 3 & 1 & -1 \\ 1 & 2 & 4 \end{bmatrix}$$

$a_{12} = 0, a_{22} = 1$ and $a_{32} = 2$.

**Step 2:** To find the determinant of A, we do the following: For $a_{12}$, since $a_{12}$ is in the first row and second column, then we virtually remove the first row and second column.

$$A = \begin{bmatrix} 1 & 0 & 2 \\ 3 & 1 & -1 \\ 1 & 2 & 4 \end{bmatrix}$$

$(-1)^{1+2} a_{12} \det \begin{bmatrix} 3 & -1 \\ 1 & 4 \end{bmatrix}$

For $a_{22}$, since $a_{22}$ is in the second row and second column, then we virtually remove the second row and second column.

$$A = \begin{bmatrix} 1 & 0 & 2 \\ 3 & 1 & -1 \\ 1 & 2 & 4 \end{bmatrix}$$

$(-1)^{2+2} a_{22} \det \begin{bmatrix} 1 & 2 \\ 1 & 4 \end{bmatrix}$

For $a_{32}$, since $a_{32}$ is in the third row and second column, then we virtually remove the third row and second column.

$$A = \begin{bmatrix} 1 & 0 & 2 \\ 3 & 1 & -1 \\ 1 & 2 & 4 \end{bmatrix}$$

$(-1)^{3+2} a_{32} \det \begin{bmatrix} 1 & 2 \\ 3 & -1 \end{bmatrix}$

**Step 3:** Add all of them together as follows:

$\det(A) = (-1)^{1+2} a_{12} \det \begin{bmatrix} 3 & -1 \\ 1 & 4 \end{bmatrix} + (-1)^{2+2} a_{22} \det \begin{bmatrix} 1 & 2 \\ 1 & 4 \end{bmatrix}$

$\quad + (-1)^{3+2} a_{32} \det \begin{bmatrix} 1 & 2 \\ 3 & -1 \end{bmatrix}$

$\det(A) = (-1)^{3}(0) \det \begin{bmatrix} 3 & -1 \\ 1 & 4 \end{bmatrix} + (-1)^{4}(1) \det \begin{bmatrix} 1 & 2 \\ 1 & 4 \end{bmatrix}$

$\quad + (-1)^{5}(2) \det \begin{bmatrix} 1 & 2 \\ 3 & -1 \end{bmatrix}$

$\det(A) = (-1)(0) \det \begin{bmatrix} 3 & -1 \\ 1 & 4 \end{bmatrix} + (1)(1) \det \begin{bmatrix} 1 & 2 \\ 1 & 4 \end{bmatrix}$

$\quad + (-1)(2) \det \begin{bmatrix} 1 & 2 \\ 3 & -1 \end{bmatrix}$

$\det(A) = (-1)(0)(12 - -1) + (1)(1)(4 - 2) + (-1)(2)(-1 - 6)$

$\det(A) = 0 + 2 + 14 = 16.$

Thus, the determinant of A is 16.

**Result 1.7.1** Let $A \in M_n(\mathbb{R})$. Then, A is invertible (non-singular) if and only if $\det(A) \neq 0$.

The above result means that if $\det(A) \neq 0$, then A is invertible (non-singular), and if A is invertible (non-singular), then $\det(A) \neq 0$.

**Example 1.7.3** Given the following matrix:
$$A = \begin{bmatrix} 2 & 3 \\ 4 & 6 \end{bmatrix}$$

Is A invertible (non-singular)?

Solution: Using result 1.7.1, we do the following:
det(A) = |A| = (2)(6) − (3)(4) = 12 − 12 = 0.
Since the determinant of A is 0, then A is non-invertible (singular).
Thus, the answer is No because A is non-invertible (singular).

**Definition 1.7.2** Given $A = \begin{bmatrix} a_{11} & a_{12} \\ a_{21} & a_{22} \end{bmatrix}$. Assume that det(A) ≠ 0 such that det(A) = $a_{11}a_{22} - a_{12}a_{21}$. To find $A^{-1}$ (the inverse of A), we use the following format that applies only for 2 × 2 matrices:

$$A^{-1} = \frac{1}{\det(A)} \begin{bmatrix} a_{22} & -a_{12} \\ -a_{21} & a_{11} \end{bmatrix}$$

$$A^{-1} = \frac{1}{a_{11}a_{22} - a_{12}a_{21}} \begin{bmatrix} a_{22} & -a_{12} \\ -a_{21} & a_{11} \end{bmatrix}$$

**Example 1.7.4** Given the following matrix:
$$A = \begin{bmatrix} 3 & 2 \\ -4 & 5 \end{bmatrix}$$
Is A invertible (non-singular)? If Yes, Find $A^{-1}$.

Solution: Using result 1.7.1, we do the following:
det(A) = |A| = (3)(5) − (2)(−4) = 15 + 8 = 23 ≠ 0.
Since the determinant of A is not 0, then A is invertible (non-singular).
Thus, the answer is Yes, there exists $A^{-1}$ according to definition 1.7.2 as follows:

$$A^{-1} = \frac{1}{\det(A)} \begin{bmatrix} 5 & -2 \\ 4 & 3 \end{bmatrix} = \frac{1}{23} \begin{bmatrix} 5 & -2 \\ 4 & 3 \end{bmatrix} = \begin{bmatrix} \frac{5}{23} & -\frac{2}{23} \\ \frac{4}{23} & \frac{3}{23} \end{bmatrix}$$

**Result 1.7.2** Let $A \in M_n(\mathbb{R})$ be a triangular matrix. Then, det(A) = multiplication of the numbers on the main diagonal of A.

There are three types of triangular matrix:

a) Upper Triangular Matrix: it has all zeros on the left side of the diagonal of $n \times n$ matrix.

(i.e. $A = \begin{bmatrix} 1 & 7 & 3 \\ 0 & 2 & 5 \\ 0 & 0 & 4 \end{bmatrix}$ is an Upper Triangular Matrix).

b) Diagonal Matrix: it has all zeros on both left and right sides of the diagonal of $n \times n$ matrix.

(i.e. $B = \begin{bmatrix} 1 & 0 & 0 \\ 0 & 2 & 0 \\ 0 & 0 & 4 \end{bmatrix}$ is a Diagonal Matrix).

c) Lower Triangular Matrix: it has all zeros on the right side of the diagonal of $n \times n$ matrix.

(i.e. $C = \begin{bmatrix} 1 & 0 & 0 \\ 5 & 2 & 0 \\ 1 & 9 & 4 \end{bmatrix}$ is a Diagonal Matrix).

**Fact 1.7.1** Let $A \in M_n(\mathbb{R})$. Then, $\det(A) = \det(A^T)$.

**Fact 1.7.2** Let $A \in M_n(\mathbb{R})$. If A is an invertible (non-singular) matrix, then $A^T$ is also an invertible (non-singular) matrix. (i.e. $(A^T)^{-1} = (A^{-1})^T$).

**Proof of Fact 1.7.2** We will show that $(A^T)^{-1} = (A^{-1})^T$.

We know from previous results that $AA^{-1} = I_n$.

By taking the transpose of both sides, we obtain:

$(AA^{-1})^T = (I_n)^T$

Then, $(A^{-1})^T A^T = (I_n)^T$

Since $(I_n)^T = I_n$, then $(A^{-1})^T A^T = I_n$.

Similarly, $(A^T)^{-1} A^T = (I_n)^T = I_n$.

Thus, $(A^T)^{-1} = (A^{-1})^T$. $\square$

**The effect of Row-Operations on determinants:**

Suppose $\propto$ is a non-zero constant, and *i and k* are row numbers in the augmented matrix.

* $\propto R_i$ , $\propto \neq 0$ (Multiply a row with a non-zero constant $\propto$).

i.e. $A = \begin{bmatrix} 1 & 2 & 3 \\ 0 & 4 & 1 \\ 2 & 0 & 1 \end{bmatrix}$ $3R_2 \dashrightarrow$ $\begin{bmatrix} 1 & 2 & 3 \\ 0 & 12 & 3 \\ 2 & 0 & 1 \end{bmatrix} = B$

Assume that $\det(A) = \gamma$ where $\gamma$ is known, then $\det(B) = 3\gamma$.

Similarly, if $\det(B) = \beta$ where $\beta$ is known, then $\det(A) = \frac{1}{3}\beta$.

* $\propto R_i + R_k \dashrightarrow R_k$ (Multiply a row with a non-zero constant $\propto$, and add it to another row).

i.e. $A = \begin{bmatrix} 1 & 2 & 3 \\ 0 & 4 & 1 \\ 2 & 0 & 1 \end{bmatrix}$ $\propto R_i + R_k \dashrightarrow R_k$

$\begin{bmatrix} 1 & 2 & 3 \\ 0 & 12 & 3 \\ 2 & 0 & 1 \end{bmatrix} = B$

Then, $\det(A) = \det(B)$.

* $R_i \leftrightarrow R_k$ (Interchange two rows). It has no effect on the determinants.

In general, the effect of Column-Operations on determinants is the same as for Row-Operations.

**Example 1.7.5** Given the following $4 \times 4$ matrix A with some Row-Operations:

A $2R_1$ --→ $A_1$ $3R_3$ --→ $A_2$ $-2R_4$ --→ $A_4$

If $\det(A) = 4$, then find $\det(A_3)$

**Solution:** Using what we have learned from the effect of determinants on Row-Operations:

$\det(A_1) = 2 * \det(A) = 2 * 4 = 8$ because $A_1$ has the first row of A multiplied by 2.

$\det(A_2) = 3 * \det(A_1) = 3 * 8 = 24$ because $A_2$ has the third row of $A_1$ multiplied by 3.

Similarly, $\det(A_3) = -2 * \det(A_2) = -2 * 24 = -48$ because $A_3$ has the fourth row of $A_2$ multiplied by -2.

**Result 1.7.3** Assume A is $n \times n$ matrix with a given $\det(A) = \gamma$. Let $\alpha$ be a number. Then, $\det(\alpha A) = \alpha^n * \gamma$.

**Result 1.7.4** Assume A and B are $n \times n$ matrices.

Then: a) $\det(A) = \det(A) * \det(B)$.

  b) Assume $A^{-1}$ exists and $B^{-1}$ exists.

  Then, $(AB)^{-1} = B^{-1}A^{-1}$.

  c) $\det(AB) = \det(BA)$.

  d) $\det(A) = \det(A^T)$.

e) If $A^{-1}$ exists, then $\det(A^{-1}) = \frac{1}{\det(A)}$.

**Proof of Result 1.7.4 (b)** We will show that $(AB)^{-1} = B^{-1}A^{-1}$.

If we multiply $(B^{-1}A^{-1})$ by $(AB)$, we obtain:

$B^{-1}(A^{-1}A)B = B^{-1}(I_n)B = B^{-1}B = I_n$.

Thus, $(AB)^{-1} = B^{-1}A^{-1}$. □

**Proof of Result 1.7.4 (e)** We will show that $\det(A^{-1}) = \frac{1}{\det(A)}$.

Since $AA^{-1} = I_n$, then $\det(AA^{-1}) = \det(I_n) = 1$.

$\det(AA^{-1}) = \det(A) * \det(A^{-1}) = 1$.

Thus, $\det(A^{-1}) = \frac{1}{\det(A)}$. □

## 1.8 Cramer's Rule

In this section, we discuss how to use Cramer's Rule to solve systems of linear equations.

**Definition 1.8.1** Given $n \times n$ system of linear equations. Let $CX = A$ be the matrix form of the given system:

$$C \begin{bmatrix} x_1 \\ x_2 \\ x_3 \\ \vdots \\ x_n \end{bmatrix} = \begin{bmatrix} a_1 \\ a_2 \\ a_3 \\ \vdots \\ a_n \end{bmatrix}$$

The system has a unique solution if and only if $\det(C) \neq 0$. Cramer's Rule tells us how to find $x_1, x_2, \ldots, x_n$ as follows:

Let $C = \begin{bmatrix} 1 & 3 & 4 \\ 1 & 2 & 1 \\ 7 & 4 & 3 \end{bmatrix}$ Then, the solutions for the system of linear equations are:

$$x_1 = \frac{\det \begin{bmatrix} a_1 & 3 & 4 \\ \vdots & 2 & 1 \\ a_n & 4 & 3 \end{bmatrix}}{\det(C)}$$

$$x_2 = \frac{\det \begin{bmatrix} 1 & a_1 & 4 \\ 1 & \vdots & 1 \\ 7 & a_n & 3 \end{bmatrix}}{\det(C)}$$

$$x_3 = \frac{\det \begin{bmatrix} 1 & 3 & a_1 \\ 1 & 2 & \vdots \\ 7 & 4 & a_n \end{bmatrix}}{\det(C)}$$

**Example 1.8.1** Use Cramer's Rule to solve the following system of linear equations:
$$\begin{cases} 2x_1 + 7x_2 = 13 \\ -10x_1 + 3x_2 = -4 \end{cases}$$

**Solution:** First of all, we write $2 \times 2$ system in the form $CX = A$ according to definition 1.8.1.
$$\begin{bmatrix} 2 & 7 \\ -10 & 3 \end{bmatrix} \begin{bmatrix} x_1 \\ x_2 \end{bmatrix} = \begin{bmatrix} 13 \\ -4 \end{bmatrix}$$
Since C in this form is $\begin{bmatrix} 2 & 7 \\ -10 & 3 \end{bmatrix}$, then

$\det(C) = (2 * 3) - (7 * (-10)) = 6 - (-70) = 76 \neq 0.$

The solutions for this system of linear equations are:

$$x_1 = \frac{\det\begin{bmatrix} 13 & 7 \\ -4 & 3 \end{bmatrix}}{\det(C)} = \frac{\det\begin{bmatrix} 13 & 7 \\ -4 & 3 \end{bmatrix}}{76} = \frac{67}{76}$$

$$x_2 = \frac{\det\begin{bmatrix} 2 & 13 \\ -10 & -4 \end{bmatrix}}{\det(C)} = \frac{\det\begin{bmatrix} 2 & 13 \\ -10 & -4 \end{bmatrix}}{76} = \frac{122}{76}$$

Thus, the solutions are $x_1 = \frac{67}{76}$ and $x_2 = \frac{122}{76}$.

## 1.9 Adjoint Method

In this section, we introduce a new mathematical method to find the inverse matrix. We also give some examples of using Adjoint Method to find inverse matrix and its entry elements.

**Example 1.9.1** Given the following matrix:
$$A = \begin{bmatrix} 1 & 0 & 2 \\ 2 & 1 & -2 \\ 0 & 0 & 2 \end{bmatrix}$$
Find $A^{-1}$ using Adjoint Method.

**Solution:** To find the inverse matrix of A, we need to do the following steps:

**Step 1:** Find the determinant of matrix A, and check whether $A^{-1}$ exists or not.
Using the previous discussed method for finding the $3 \times 3$ inverse matrix, we get: $\det(A) = 2 \neq 0$.
Therefore, $A^{-1}$ exists.

**Step 2:** Calculate the coefficient matrix of A such that $A = C_{3\times 3}$ because A is $3 \times 3$ matrix.

$$C = \begin{bmatrix} c_{11} & c_{12} & c_{13} \\ c_{21} & c_{22} & c_{23} \\ c_{31} & c_{32} & c_{33} \end{bmatrix}$$

In general, to find every element of coefficient matrix A, we need to use the following form:

$c_{ik} = (-1)^{i+k} \det[\text{Remove the } i^{th} \text{ row and } k^{th} \text{ column of } A]$

Now, using the above form, we can find all coefficients of matrix A:

$c_{11} = (-1)^{1+1} \det \begin{bmatrix} 1 & -2 \\ 0 & 2 \end{bmatrix} = 2$

$c_{12} = (-1)^{1+2} \det \begin{bmatrix} 0 & -2 \\ 0 & 2 \end{bmatrix} = 0$

$c_{13} = (-1)^{1+3} \det \begin{bmatrix} 0 & 1 \\ 0 & 0 \end{bmatrix} = 0$

$c_{21} = (-1)^{2+1} \det \begin{bmatrix} 0 & 2 \\ 0 & 2 \end{bmatrix} = 0$

$c_{22} = (-1)^{2+2} \det \begin{bmatrix} 1 & 2 \\ 0 & 2 \end{bmatrix} = 2$

$c_{23} = (-1)^{2+3} \det \begin{bmatrix} 1 & 0 \\ 0 & 0 \end{bmatrix} = 0$

$c_{31} = (-1)^{3+1} \det \begin{bmatrix} 0 & 2 \\ 1 & -2 \end{bmatrix} = -2$

$c_{32} = (-1)^{3+2} \det \begin{bmatrix} 1 & 2 \\ 0 & -2 \end{bmatrix} = 2$

$c_{33} = (-1)^{3+3} \det \begin{bmatrix} 1 & 0 \\ 0 & 1 \end{bmatrix} = 1$

Hence, $C = \begin{bmatrix} 2 & 0 & 0 \\ 0 & 2 & 0 \\ -2 & 2 & 1 \end{bmatrix}$

**Step 3:** Find the adjoint of matrix A as follows:
In general, adjoint of A = adj(A) = $C^T$.

Thus, adj(A) = $C^T$ = $\begin{bmatrix} 2 & 0 & -2 \\ 0 & 2 & 2 \\ 0 & 0 & 1 \end{bmatrix}$

This formula is always true such that

$A * \text{adj}(A) = \det(A) I_n$.

If $\det(A) \neq 0$, then $A \left[ \frac{1}{\det(A)} * \text{adj}(A) \right] = I_n$.

Hence, $A^{-1} = \frac{1}{\det(A)} * \text{adj}(A) = \frac{1}{\det(A)} * C^T$.

$A^{-1} = \frac{1}{\det(A)} * C^T = \frac{1}{2} * \begin{bmatrix} 2 & 0 & -2 \\ 0 & 2 & 2 \\ 0 & 0 & 1 \end{bmatrix} = \begin{bmatrix} 1 & 0 & -1 \\ 0 & 1 & 1 \\ 0 & 0 & \frac{1}{2} \end{bmatrix}$

**Example 1.9.2** Given the following matrix:

$A = \begin{bmatrix} 2 & 0 & -2 & 1 \\ -2 & 1 & 2 & 4 \\ -4 & -1 & 3 & 0 \\ 0 & 0 & 0 & 4 \end{bmatrix}$

Find (2,4) entry of $A^{-1}$.

**Solution:** To find the (2,4) entry of $A^{-1}$, we need to do the following steps:

**Step 1:** Find the determinant of matrix A, and check whether $A^{-1}$ exists or not.
Use the Row-Operation Method, we obtain the following:

$\begin{bmatrix} 2 & 0 & -2 & 1 \\ -2 & 1 & 2 & 4 \\ -4 & -1 & 3 & 0 \\ 0 & 0 & 0 & 4 \end{bmatrix} \begin{matrix} R_1 + R_2 \to R_2 \\ 2R_2 + R_3 \to R_3 \end{matrix} \begin{bmatrix} 2 & 0 & -2 & 1 \\ 0 & 1 & 0 & 5 \\ 0 & -1 & -1 & 2 \\ 0 & 0 & 0 & 4 \end{bmatrix}$

$$R_2 + R_3 \to R_3 \begin{bmatrix} 2 & 0 & -2 & 1 \\ 0 & 1 & 0 & 5 \\ 0 & 0 & -1 & 7 \\ 0 & 0 & 0 & 4 \end{bmatrix}$$

Now, using result 1.7.2 for finding the determinant of the upper triangular matrix:

Since $A \in M_n(\mathbb{R})$ is a triangular matrix. Then, det(A) = multiplication of the numbers on the main diagonal of A.

$$\begin{bmatrix} 2 & 0 & -2 & 1 \\ 0 & 1 & 0 & 5 \\ 0 & 0 & -1 & 7 \\ 0 & 0 & 0 & 4 \end{bmatrix}$$

$\det(A) = 2 * 1 * -1 * 4 = -8 \neq 0$.

Therefore, $A^{-1}$ exists.

**Step 2:** Use the following general form for (i,k) entry of $A^{-1}$:

$(i,k) - entry\ of\ A^{-1} =$

$$\frac{c_{ki}}{\det(A)} = \frac{(-1)^{i+k} \det[Remove\ the\ k^{th}\ row\ and\ i^{th}\ column\ of\ A]}{\det(A)}$$

Now, using the above form, we can find $(2,4)$-entry of matrix A:

$$\frac{c_{42}}{\det(A)} = \frac{(-1)^{2+4} \det \begin{bmatrix} 2 & -2 & 1 \\ -2 & 2 & 4 \\ -4 & 3 & 0 \end{bmatrix}}{\det(A)} = \frac{(-1)^{2+4} \det \begin{bmatrix} 2 & -2 & 1 \\ -2 & 2 & 4 \\ -4 & 3 & 0 \end{bmatrix}}{-8}$$

**Step 3:** Find the determinant of $\begin{bmatrix} 2 & -2 & 1 \\ -2 & 2 & 4 \\ -4 & 3 & 0 \end{bmatrix}$.

Let's call the above matrix F such that F =
$$\begin{bmatrix} 2 & -2 & 1 \\ -2 & 2 & 4 \\ -4 & 3 & 0 \end{bmatrix}$$
To find the determinant of F, we need to use the Row-Operation Method to reduce the matrix F.

$$F = \begin{bmatrix} 2 & -2 & 1 \\ -2 & 2 & 4 \\ -4 & 3 & 0 \end{bmatrix} \begin{array}{c} R_1 + R_2 \to R_2 \\ 2R_1 + R_3 \to R_3 \end{array} \begin{bmatrix} 2 & -2 & 1 \\ 0 & 0 & 5 \\ 0 & -1 & 2 \end{bmatrix}$$

$$R_2 \leftrightarrow R_3 \begin{bmatrix} 2 & -2 & 1 \\ 0 & -1 & 2 \\ 0 & 0 & 5 \end{bmatrix}$$

Let's call the above matrix D such that $D = \begin{bmatrix} 2 & -2 & 1 \\ 0 & -1 & 2 \\ 0 & 0 & 5 \end{bmatrix}$

Now, using result 1.7.2 for finding the determinant of this upper triangular matrix:

$$D = \begin{bmatrix} 2 & -2 & 1 \\ 0 & -1 & 2 \\ 0 & 0 & 5 \end{bmatrix}$$

$\det(D) = 2 * -1 * 5 = -10 \neq 0$.
Therefore, $\det(F) = -\det(D) = 10 \neq 0$.
Thus, The (2,4)-entry of matrix A is:

$$\frac{c_{42}}{\det(A)} = \frac{(-1)^{2+4} \det \begin{bmatrix} 2 & -2 & 1 \\ -2 & 2 & 4 \\ -4 & 3 & 0 \end{bmatrix}}{-8} = (-1)^6 * \frac{10}{-8} = \frac{10}{-8} = -\frac{5}{4}$$
$$= -1.25$$

∴ (2,4)-entry of matrix A = −1.25

## 1.10 Exercises

1. Solve the following system of linear equations:
$$\begin{cases} 2x_3 - x_4 + x_6 = 10 \\ -x_1 + 3x_2 + x_5 = 0 \\ x_1 + 2x_2 + 2x_3 - x_4 + 2x_5 + x_6 = 12 \end{cases}$$

2. Given an augmented matrix of a system: $\begin{pmatrix} 1 & -1 & c & | & 2 \\ -1 & 1 & 4 & | & 3 \\ 4 & -3 & b & | & 10 \end{pmatrix}$

   a. For what values of c and b will the system be consistent?

   b. If the system is consistent, when do you have a unique solution?.

3. Let $F = \begin{bmatrix} -1 & 2 & -1 & 3 \\ 0 & 1 & 1 & 1 \\ -2 & 0 & -1 & 1 \end{bmatrix}$ and $H = \begin{bmatrix} 2 & 1 & 1 \\ -4 & 1 & 0 \\ 1 & 1 & 1 \\ -1 & 0 & 3 \end{bmatrix}$

   a. Find the third row of HF.

   b. Find the third column of FH.

   c. Let HF = A. Find $a_{42}$.

4. Let $A = \begin{bmatrix} 1 & 2 & 4 \\ -1 & -2 & -3 \\ -2 & -3 & -8 \end{bmatrix}$. If possible find $A^{-1}$.

5. Given A is 2 × 4 matrix such that
$\mathbf{A} \underbrace{R_1 \to \mathbf{A_1}}_{} \underbrace{2R_1 + R_2 \to R_2}_{} \mathbf{A_2}$.

a. Find two elementary matrices say $E_1, E_2$ such that $E_1 E_2 A = A_2$.

b. Find two elementary matrices say $F_1, F_2$ such that $F_1 F_2 A_2 = A$.

6. Let $A = \begin{bmatrix} 1 & 2 & 4 \\ -1 & -2 & 3 \\ -2 & -3 & -7 \end{bmatrix}$. Find det(A). Is A invertible? Explain.

7. Let $A = \begin{bmatrix} 4 & -2 \\ -3 & 2 \end{bmatrix}$. Is A invertible? If yes, find $A^{-1}$.

8. Use Cramer's Rule to find the solution to $x_2$ in the system:
$$\begin{cases} 2x_1 + x_2 - x_3 = 2 \\ -2x_1 + 4x_2 + 2x_3 = 8 \\ -2x_1 - x_2 + 8x_3 = -2 \end{cases}$$

9. Let $A = \begin{bmatrix} 2 & -4 & 2 & 1 \\ -2 & 0 & 2 & -1 \\ 1 & -2 & 12 & 4 \\ -2 & 4 & -2 & 12 \end{bmatrix}$. Find (2,4) entry of $A^{-1}$.

10. Find a $2 \times 2$ matrix $A$ such that
$$\begin{bmatrix} 5 & 5 \\ 1 & 4 \end{bmatrix} A + 3I_2 = 2A + \begin{bmatrix} -2 & 2 \\ 3 & 2 \end{bmatrix}.$$

11. Given $A^{-1} = \begin{bmatrix} 2 & -2 & -2 \\ -2 & 3 & 0 \\ -2 & 2 & 3 \end{bmatrix}$ and $B = \begin{bmatrix} 1 & -2 & -2 \\ -4 & 2 & 2 \\ 0 & 1 & -2 \end{bmatrix}$

Solve the system $AX = \begin{bmatrix} 0 \\ -1 \\ 1 \end{bmatrix}$

# Chapter 2

# Vector Spaces

We start this chapter reviewing some concepts of set theory, and we discuss some important concepts of vector spaces including span and dimension. In the remaining sections we introduce the concept of linear independence. At the end of this chapter we discuss other concepts such as subspace and basis.

## 2.1 Span and Vector Spaces

In this section, we review some concepts of set theory, and we give an introduction to span and vector spaces including some examples related to these concepts.

Before reviewing the concepts of set theory, it is recommended to revisit section 1.4, and read the notations of numbers and the representation of the three sets of numbers in figure 1.4.1.

Let's explain some symbols and notations of set theory:

$3 \in \mathbb{Z}$ means that 3 is an element of $\mathbb{Z}$.

$\frac{1}{2} \notin \mathbb{Z}$ means that $\frac{1}{2}$ is not an element of $\mathbb{Z}$.

$\{\}$ means that it is a set.

{5} means that 5 is a subset of ℤ, and the set consists of exactly one element which is 5.

**Definition 2.1.1** The span of a certain set is the set of all possible linear combinations of the subset of that set.

**Example 2.1.1** Find Span{1}.

**Solution:** According to definition 2.1.1, then the span of the set {1} is the set of all possible linear combinations of the subset of {1} which is 1.

Hence, Span{1} = ℝ.

**Example 2.1.2** Find Span{(1,2),(2,3)}.

**Solution:** According to definition 2.1.1, then the span of the set {(1,2),(2,3)} is the set of all possible linear combinations of the subsets of {(1,2),(2,3)} which are (1,2) and (2,3). Thus, the following is some possible linear combinations:

$(1,2) = 1 * (1,2) + 0 * (2,3)$

$(2,3) = 0 * (1,2) + 1 * (2,3)$

$(5,8) = 1 * (1,2) + 2 * (2,3)$

Hence, $\{(1,2), (2,3), (5,8)\} \in \text{Span}\{(1,2), (2,3)\}$.

**Example 2.1.3** Find Span{0}.

**Solution:** According to definition 2.1.1, then the span of the set {0} is the set of all possible linear combinations of the subset of {0} which is 0.

Hence, Span$\{0\} = 0$.

**Example 2.1.4** Find Span$\{c\}$ where c is a non-zero integer.

**Solution:** Using definition 2.1.1, the span of the set $\{c\}$ is the set of all possible linear combinations of the subset of $\{c\}$ which is $c \neq 0$.

Thus, Span$\{c\} = \mathbb{R}$.

**Definition 2.1.2** $\mathbb{R}^n = \{(a_1, a_2, a_3, \ldots, a_n) | a_1, a_2, a_3, \ldots, a_n \in \mathbb{R}\}$ is a set of all points where each point has exactly $n$ coordinates.

**Definition 2.1.3** $(V, +, \cdot)$ is a vector space if satisfies the following:

  a. For every $v_1, v_2 \in V$, $v_1 + v_2 \in V$.
  b. For every $\alpha \in \mathbb{R}$ and $v \in V$, $\alpha v \in V$.

(i.e. Given $Span\{x, y\}$ and $set\ \{x, y\}$, then $\sqrt{10}x + 2y \in Span\{x, y\}$. Let's assume that $v \in Span\{x, y\}$, then $v = c_1 x + c_2 y$ for some numbers $c_1$ and $c_2$).

## 2.2 The Dimension of Vector Space

In this section, we discuss how to find the dimension of vector space, and how it is related to what we have learned in section 2.1.

**Definition 2.2.1** Given a vector space $V$, the dimension of $V$ is the number of minimum elements needed in $V$ so that their $Span$ is equal to $V$, and it is denoted by $\dim(V)$. (i.e. $\dim(\mathbb{R}) = 1$ and $\dim(\mathbb{R}^2) = 2$).

**Result 2.2.1** $\dim(\mathbb{R}^n) = n$.

**Proof of Result 2.2.1** We will show that $\dim(\mathbb{R}^n) = n$.

Claim: $D = Span\{(1,0),(0.1)\} = \mathbb{R}^2$

$\alpha_1(1,0) + \alpha_2(0,1) = (\alpha_1, \alpha_2) \in \mathbb{R}^2$

Thus, $D$ is a subset of $\mathbb{R}^2$ ($D \subseteq \mathbb{R}^2$).

For every $x_1, y_1 \in \mathbb{R}$, $(x_1, y_1) \in \mathbb{R}^2$.

Therefore, $(x_1, y_1) = x_1(1,0) + y_1(0,1) \in D$.

We prove the above claim, and hence

$\dim(\mathbb{R}^n) = n$. ∎

**Fact 2.2.1** $Span\{(3,4)\} \neq \mathbb{R}^2$.

**Proof of Fact 2.2.1** We will show that $Span\{(3,4)\} \neq \mathbb{R}^2$.

Claim: $F = Span\{(6,5)\} \neq \mathbb{R}^2$ where $(6,5) \in \mathbb{R}^2$.

We cannot find a number $\alpha$ such that $(6,5) = \alpha(3,4)$

We prove the above claim, and hence $Span\{(3,4)\} \neq \mathbb{R}^2$. □

**Fact 2.2.2** $Span\{(1,0),(0,1)\} = \mathbb{R}^2$.

**Fact 2.2.3** $Span\{(2,1),(1,0.5)\} \neq \mathbb{R}^2$.

## 2.3 Linear Independence

In this section, we learn how to determine whether vector spaces are linearly independent or not.

**Definition 2.3.1** Given a vector space $(V, +, \cdot)$, we say $v_1, v_2, \dots, v_n \in V$ are linearly independent if none of them is a linear combination of the remaining $v_i's$.

(i.e. $(3,4),(2,0) \in \mathbb{R}$ are linearly independent because we cannot write them as a linear combination of each other, in other words, we cannot find a number $\alpha_1, \alpha_2$ such that $(3,4) = \alpha_1(2,0)$ and $(2,0) = \alpha_2(3,4)$).

**Definition 2.3.2** Given a vector space $(V, +, \cdot)$, we say $v_1, v_2, \dots, v_n \in V$ are linearly dependent if at least one of $v_i's$ is a linear combination of the others.

**Example 2.3.1** Assume $v_1$ and $v_2$ are linearly independent. Show that $v_1$ and $3v_1 + v_2$ are linearly independent.

**Solution:** We will show that $v_1$ and $3v_1 + v_2$ are linearly independent. Using proof by contradiction, we assume that $v_1$ and $3v_1 + v_2$ are linearly dependent. For some non-zero number $c_1$, $v_1 = c_1(3v_1 + v_2)$.
Using the distribution property and algebra, we obtain:
$$v_1 = 3v_1c_1 + v_2c_1$$
$$v_1 - 3v_1c_1 = v_2c_1$$
$$v_1(1 - 3c_1) = v_2c_1$$
$$\frac{(1 - 3c_1)}{c_1}v_1 = v_2$$

Thus, none of $v_1$ and $3v_1 + v_2$ is a linear combination of the others which means that $v_1$ and $3v_1 + v_2$ are linearly independent. This is a contradiction. Therefore, our assumption that $v_1$ and $3v_1 + v_2$ were linearly dependent is false. Hence, $v_1$ and $3v_1 + v_2$ are linearly independent. □

**Example 2.3.2** Given the following vectors:
$$v_1 = (1,0,-2)$$
$$v_2 = (-2,2,1)$$
$$v_3 = (-1,0,5)$$

Are these vectors independent elements?

**Solution:** First of all, to determine whether these vectors are independent elements or not, we need to write these vectors as a matrix.

$\begin{bmatrix} 1 & 0 & -2 \\ -2 & 2 & 1 \\ -1 & 0 & 5 \end{bmatrix}$ Each point is a row-operation. We need to reduce this matrix to Semi-Reduced Matrix.

**Definition 2.3.3** Semi-Reduced Matrix is a reduced-matrix but the leader numbers can be any non-zero number.

Now, we apply the Row-Reduction Method to get the Semi-Reduced Matrix as follows:

$\begin{bmatrix} 1 & 0 & -2 \\ -2 & 2 & 1 \\ -1 & 0 & 5 \end{bmatrix} \begin{matrix} 2R_1 + R_2 \to R_2 \\ R_1 + R_3 \to R_3 \end{matrix} \begin{bmatrix} 1 & 0 & -2 \\ 0 & \boxed{2} & -3 \\ 0 & 0 & 3 \end{bmatrix}$ This is a Semi-Reduced Matrix.

Since none of the rows in the Semi-Reduced Matrix become zero-row, then the elements are independent because we cannot write at least one of them as a linear combination of the others.

**Example 2.3.3** Given the following vectors:

$$v_1 = (1, -2, 4, 6)$$
$$v_2 = (-1, 2, 0, 2)$$
$$v_3 = (1, -2, 8, 14)$$

Are these vectors independent elements?

**Solution:** First of all, to determine whether these vectors are independent elements or not, we need to write these vectors as a matrix.

$$\begin{bmatrix} 1 & -2 & 4 & 6 \\ -1 & 2 & 0 & 2 \\ 1 & -2 & 8 & 14 \end{bmatrix}$$ Each point is a row-operation. We need to reduce this matrix to Semi-Reduced Matrix.

Now, we apply the Row-Reduction Method to get the Semi-Reduced Matrix as follows:

$$\begin{bmatrix} 1 & -2 & 4 & 6 \\ -1 & 2 & 0 & 2 \\ 1 & -2 & 8 & 14 \end{bmatrix} \begin{array}{c} R_1 + R_2 \to R_2 \\ -R_1 + R_3 \to R_3 \end{array} \begin{bmatrix} 1 & -2 & 4 & 6 \\ 0 & 0 & \boxed{4} & 8 \\ 0 & 0 & 4 & 8 \end{bmatrix}$$

$-R_2 + R_3 \to R_3$ $\begin{bmatrix} 1 & -2 & 4 & 6 \\ 0 & 0 & \boxed{4} & 8 \\ \boxed{0 & 0 & 0 & 0} \end{bmatrix}$ This is a Semi-Reduced Matrix.

Since there is a zero-row in the Semi-Reduced Matrix, then the elements are dependent because we can write at least one of them as a linear combination of the others.

## 2.4 Subspace and Basis

In this section, we discuss one of the most important concepts in linear algebra that is known as subspace. In addition, we give some examples explaining how to find the basis for subspace.

**Definition 2.4.1** Subspace is a vector space but we call it a subspace because it lives inside a bigger vector space. (i.e. Given vector spaces $V$ and $D$, then according to the figure 2.4.1, $D$ is called a subspace of $V$).

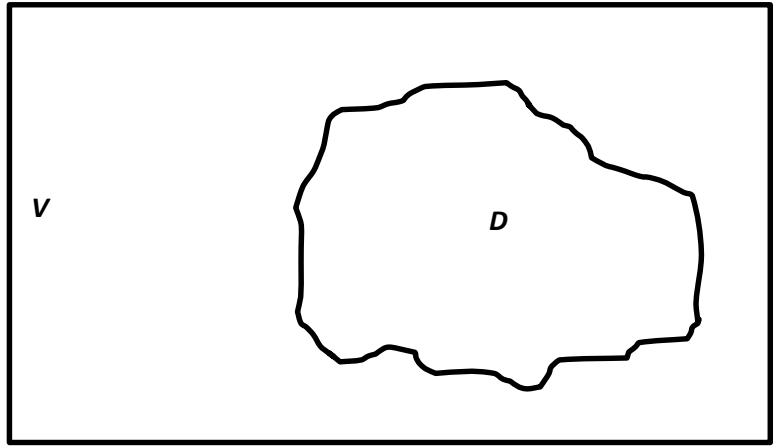

Figure 2.4.1: Subspace of $V$

**Fact 2.4.1** Every vector space is a subspace of itself.

**Example 2.4.1** Given a vector space $L = \{(c, 3c) | c \in \mathbb{R}\}$.

  a. Does $L$ live in $\mathbb{R}^2$?
  b. Does $L$ equal to $\mathbb{R}^2$?
  c. Is $L$ a subspace of $\mathbb{R}^2$?
  d. Does $L$ equal to $Span\{(0,3)\}$?
  e. Does $L$ equal to $Span\{(1,3),(2,6)\}$?

**Solution:** To answer all these questions, we need first to draw an equation from this vector space, say $y = 3x$. The following figure represents the graph of the above equation, and it passes through a point (1,3).

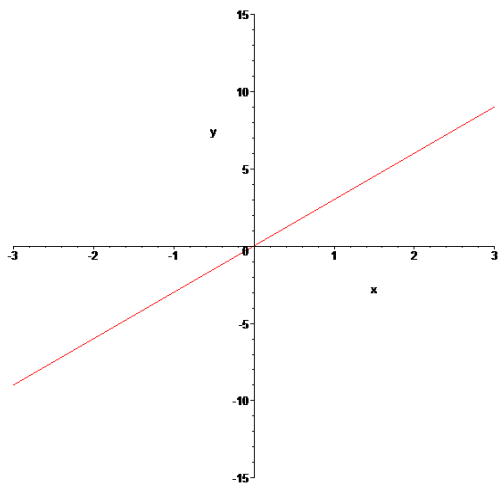

Figure 2.4.2: Graph of $y = 3x$

Now, we can answer the given questions as follows:

**Part a:** Yes; $L$ lives in $\mathbb{R}^2$.

**Part b:** No; $L$ does not equal to $\mathbb{R}^2$. To show that we prove the following claim:

Claim: $L = Span\{(5,15)\} \neq \mathbb{R}^2$ where $(5,15) \in \mathbb{R}^2$.

It is impossible to find a number $\alpha = 3$ such that

$$(20,60) = \alpha(5,15)$$

because in this case $\alpha = 4$ where $(20,60) = 4(5,15)$.

We prove the above claim, and $Span\{(5,15)\} \neq \mathbb{R}^2$.

Thus, $L$ does not equal to $\mathbb{R}^2$

**Part c:** Yes; $L$ is a subspace of $\mathbb{R}^2$ because $L$ lives inside a bigger vector space which is $\mathbb{R}^2$.

**Part d:** No; according to the graph in figure 2.4.2, (0,3) does not belong to $L$.

**Part e:** Yes; because we can write (1,3) and (2,6) as a linear combination of each other.

$\alpha_1(1,3) + \alpha_2(2,6) = \{(\alpha_1 + 2\alpha_2), (3\alpha_1 + 6\alpha_2)\}$

$\alpha_1(1,3) + \alpha_2(2,6) = \{(\alpha_1 + 2\alpha_2), 3(\alpha_1 + 2\alpha_2)\}$

Assume $c = (\alpha_1 + 2\alpha_2)$, then we obtain:

$\alpha_1(1,3) + \alpha_2(2,6) = \{(c, 3c) | c \in \mathbb{R}\} = L$.

Thus, $L = Span\{(1,3), (2,6)\}$.

**Result 2.4.1** $L$ is a subspace of $\mathbb{R}^2$ if satisfies the following:

a. $L$ lives inside $\mathbb{R}^2$.
b. $L$ has only lines through the origin (0,0).

**Example 2.4.2** Given a vector space

$D = \{(a, b, 1) | a, b \in \mathbb{R}\}$.

a. Does $D$ live in $\mathbb{R}^3$?
b. Is $D$ a subspace of $\mathbb{R}^3$?

**Solution:** Since the equation of the above vector space is a three-dimensional equation, there is no need to draw it because it is difficult to draw it exactly. Thus, we can answer the above questions immediately.

**Part a:** Yes; $D$ lives inside $\mathbb{R}^3$.

**Part b:** No; since $(0,0,0) \notin D$, then $D$ is not a subspace of $\mathbb{R}^3$.

**Fact 2.4.2** Assume $D$ lives inside $\mathbb{R}^n$. If we can write $D$ as a *Span*, then it is a subspace of $\mathbb{R}^n$.

**Fact 2.4.3** Assume $D$ lives inside $\mathbb{R}^n$. If we cannot write $D$ as a *Span*, then it is not a subspace of $\mathbb{R}^n$.

**Fact 2.4.4** Assume $D$ lives inside $\mathbb{R}^n$. If $(0,0,0,...,0)$ is in $D$, then $D$ is a subspace of $\mathbb{R}^n$.

**Fact 2.4.5** Assume $D$ lives inside $\mathbb{R}^n$. If $(0,0,0,...,0)$ is not in $D$, then $D$ is not a subspace of $\mathbb{R}^n$.

Now, we list the main results on $\mathbb{R}^n$:

**Result 2.4.2** Maximum number of independent points is $n$.

**Result 2.4.3** Choosing any $n$ independent points in $\mathbb{R}^n$, say $Q_1, Q_2, ..., Q_n$, then $\mathbb{R}^n = Span\{Q_1, Q_2, ..., Q_n\}$.

**Result 2.4.4** $\dim(\mathbb{R}^n) = n$.

Results 2.4.3 and 2.4.4 tell us the following: In order to get all $\mathbb{R}^n$, we need exactly $n$ independent points.

**Result 2.4.5** Assume $\mathbb{R}^n = Span\{Q_1, Q_2, ..., Q_k\}$, then $k \geq n$ ($n$ points of the $Q_k's$ are independents).

**Definition 2.4.2** Basis is the set of points that is needed to *Span* the vector space.

**Example 2.4.3** Let $D = Span\{(1,-1,0), (2,2,1), (0,4,1)\}$.

a. Find $\dim(D)$.

b. Find a basis for $D$.

**Solution:** First of all, we have infinite set of points, and $D$ lives inside $\mathbb{R}^3$. Let's assume the following:

$$v_1 = (1,-1,0)$$
$$v_2 = (2,2,1)$$
$$v_3 = (0,4,1)$$

**Part a:** To find $\dim(D)$, we check whether $v_1, v_2$ and $v_3$ are dependent elements or not. Using what we have learned so far from section 2.3: We need to write these vectors as a matrix.

$\begin{bmatrix} 1 & -1 & 0 \\ 2 & 2 & 1 \\ 0 & 4 & 1 \end{bmatrix}$ Each point is a row-operation. We need to reduce this matrix to Semi-Reduced Matrix.

Now, we apply the Row-Reduction Method to get the Semi-Reduced Matrix as follows:

$\begin{bmatrix} 1 & -1 & 0 \\ 2 & 2 & 1 \\ 0 & 4 & 1 \end{bmatrix} -2R_1 + R_2 \to R_2 \begin{bmatrix} 1 & -1 & 0 \\ 0 & \boxed{4} & 1 \\ 0 & 4 & 1 \end{bmatrix} -R_2 + R_3 \to R_3$

$\begin{bmatrix} 1 & -1 & 0 \\ 0 & \boxed{4} & 1 \\ \boxed{0 & 0 & 0} \end{bmatrix}$ This is a Semi-Reduced Matrix.

Since there is a zero-row in the Semi-Reduced Matrix, then these elements are dependent because we can write at least one of them as a linear combination of the others. Only two points survived in the Semi-Reduced Matrix. Thus, $\dim(D) = 2$.

**Part b:** $D$ is a plane that passes through the origin $(0,0,0)$. Since $\dim(D) = 2$, then any two independent points in $D$ will form a basis for $D$. Hence, the following are some possible bases for $D$:

Basis for $D$ is $\{(1,-1,0),(2,2,1)\}$.

Another basis for $D$ is $\{(1,-1,0),(0,4,1)\}$.

**Result 2.4.6** It is always true that $|Basis| = dim(D)$.

**Example 2.4.4** Given the following:

$M = Span\{(-1,2,0,0),(1,-2,3,0),(-2,0,3,0)\}$.

Find a basis for $M$.

**Solution:** We have infinite set of points, and $M$ lives inside $\mathbb{R}^4$. Let's assume the following:

$$v_1 = (-1,2,0,0)$$
$$v_2 = (1,-2,3,0)$$
$$v_3 = (-2,0,3,0)$$

We check if $v_1, v_2$ and $v_3$ are dependent elements. Using what we have learned so far from section 2.3 and example 2.4.3: We need to write these vectors as a matrix.

$\begin{bmatrix} -1 & 2 & 0 & 0 \\ 1 & -2 & 3 & 0 \\ -2 & 0 & 3 & 0 \end{bmatrix}$ Each point is a row-operation. We
need to reduce this matrix to Semi-Reduced Matrix.

Now, we apply the Row-Reduction Method to get the Semi-Reduced Matrix as follows:

$$\begin{bmatrix} -1 & 2 & 0 & 0 \\ 1 & -2 & 3 & 0 \\ -2 & 0 & 3 & 0 \end{bmatrix} \begin{matrix} R_1 + R_2 \to R_2 \\ -2R_1 + R_3 \to R_3 \end{matrix} \begin{bmatrix} -1 & 2 & 0 & 0 \\ 0 & 0 & \boxed{3} & 0 \\ 0 & -4 & 3 & 0 \end{bmatrix}$$

$$-R_2 + R_3 \to R_3 \begin{bmatrix} -1 & 2 & 0 & 0 \\ 0 & 0 & 3 & 0 \\ 0 & -4 & 0 & 0 \end{bmatrix}$$ This is a Semi-Reduced Matrix.

Since there is no zero-row in the Semi-Reduced Matrix, then these elements are independent. All the three points survived in the Semi-Reduced Matrix. Thus, $\dim(M) = 3$. Since $\dim(M) = 3$, then any three independent points in $M$ from the above matrices will form a basis for $M$. Hence, the following are some possible bases for $M$:

Basis for $M$ is $\{(-1,2,0,0), (0,0,3,0), (0,-4,0,0)\}$.

Another basis for $M$ is $\{(-1,2,0,0), (0,0,3,0), (0,-4,3,0)\}$.

Another basis for $M$ is $\{(-1,2,0,0), (1,-2,3,0), (-2,0,3,0)\}$.

**Example 2.4.5** Given the following:

$W = Span\{(a, -2a + b, -a) | a, b \in \mathbb{R}\}$.

a. Show that $W$ is a subspace of $\mathbb{R}^3$.

b. Find a basis for $W$.

c. Rewrite $W$ as a *Span*.

**Solution:** We have infinite set of points, and $W$ lives inside $\mathbb{R}^3$.

**Part a:** We write each coordinate of $W$ as a linear combination of the free variables $a$ and $b$.

$$a = 1 \cdot a + 0 \cdot b$$

$$-2a + b = -2 \cdot a + 1 \cdot b$$

$$-a = -1 \cdot a + 0 \cdot b$$

Since it is possible to write each coordinate of $W$ as a linear combination of the free variables $a$ and $b$, then we conclude that $W$ is a subspace of $\mathbb{R}^3$.

**Part b:** To find a basis for $W$, we first need to find $\dim(W)$. To find $\dim(W)$, let's play a game called (ON-OFF GAME) with the free variables $a$ and $b$.

| $a$ | $b$ | $Point$ |
|---|---|---|
| 1 | 0 | $(1,-2,-1)$ |
| 0 | 1 | $(0,1,0)$ |

Now, we check for independency: We already have the Semi-Reduced Matrix: $\begin{bmatrix} 1 & -2 & -1 \\ 0 & 1 & 0 \end{bmatrix}$. Thus, $\dim(W) = 2$.

Hence, the basis for $W$ is $\{(1,-2,-1),(0,1,0)\}$.

**Part b:** Since we found the basis for $W$, then it is easy to rewrite $W$ as a *Span* as follows:

$W = Span\{(1,-2,-1),(0,1,0)\}$.

**Fact 2.4.6** $\dim(W) \leq Number of\ Free - Variables$.

**Example 2.4.6** Given the following:

$H = Span\{(a^2, 3b + a, -2c, a + b + c) | a, b, c \in \mathbb{R}\}$.

Is $H$ a subspace of $\mathbb{R}^4$?

**Solution:** We have infinite set of points, and $H$ lives inside $\mathbb{R}^4$. We try write each coordinate of $H$ as a linear combination of the free variables $a, b$ and $c$.

$a^2 = $ Fixed Number $\cdot\ a + $ Fixed Number $\cdot\ b + $ Fixed Number $\cdot\ c$

$a^2$ is not a linear combination of $a, b$ and $c$.

We assume that $w = (1,1,0,1) \in H$, and $a = 1, b = c = 0$.

If $\alpha = -2$, then $-2 \cdot w = -2 \cdot (1,1,0,1) = (-2, -2, 0, -2) \notin H$.

Since it is impossible to write each coordinate of $H$ as a linear combination of the free variables $a, b$ and $c$, then we conclude that $H$ is not a subspace of $\mathbb{R}^4$.

**Example 2.4.7** Form a basis for $\mathbb{R}^4$.

**Solution:** We just need to select any random four independent points, and then we form a $4 \times 4$ matrix with four independent rows as follows:

$\begin{bmatrix} 2 & 3 & 0 & 4 \\ 0 & 5 & 1 & 1 \\ 0 & 0 & 2 & 3 \\ 0 & 0 & 0 & \pi^e \end{bmatrix}$ Note: $\pi^e$ is a number.

Let's assume the following:

$$v_1 = (2,3,0,4)$$
$$v_2 = (0,5,1,1)$$

$$v_3 = (0,0,2,3)$$
$$v_4 = (0,0,0,\pi^e)$$

Thus, the basis for $\mathbb{R}^4 = \{v_1, v_2, v_3, v_4\}$, and

$Span\{v_1, v_2, v_3, v_4\} = \mathbb{R}^4$.

**Example 2.4.8** Form a basis for $\mathbb{R}^4$ that contains the following two independent points:

$(0,2,1,4)$ and $(0,-2,3,-10)$.

**Solution:** We need to add two more points to the given one so that all four points are independent. Let's assume the following:

$$v_1 = (0,2,1,4)$$
$$v_2 = (0,-2,3,-10)$$
$$v_3 = (0,0,4,-6) \text{ This is a random point.}$$
$$v_4 = (0,0,0,1000) \text{ This is a random point.}$$

Then, we need to write these vectors as a matrix.

$\begin{bmatrix} 0 & 2 & 1 & 4 \\ 0 & -2 & 3 & -10 \\ 0 & 0 & 4 & -6 \\ 0 & 0 & 0 & 1000 \end{bmatrix}$ Each point is a row-operation. We need to reduce this matrix to Semi-Reduced Matrix.

Now, we apply the Row-Reduction Method to get the Semi-Reduced Matrix as follows:

$$\begin{bmatrix} 0 & 2 & 1 & 4 \\ 0 & -2 & 3 & -10 \\ 0 & 0 & 4 & -6 \\ 0 & 0 & 0 & 1000 \end{bmatrix} \begin{array}{c} R_1 + R_2 \to R_2 \end{array} \begin{bmatrix} 0 & 2 & 1 & 4 \\ 3 & 0 & 5 & 30 \\ 0 & 0 & 4 & -6 \\ 0 & 0 & 0 & 1000 \end{bmatrix}$$

This is a Semi-Reduced Matrix.

Thus, the basis for $\mathbb{R}^4$ is

$\{(0,2,1,4), (0,-2,3,-10), (3,0,5,30), (0,0,0,1000)\}$.

**Example 2.4.9** Given the following:

$D = Span\{(1,1,1,1), (-1,-1,0,0), (0,0,1,1)\}$

Is $(1,1,2,2) \in D$?

**Solution:** We have infinite set of points, and $D$ lives inside $\mathbb{R}^4$. There are two different to solve this example:

**The First Way:** Let's assume the following:

$$v_1 = (1,1,1,1)$$
$$v_2 = (-1,-1,0,0)$$
$$v_3 = (0,0,1,1)$$

We start asking ourselves the following question:

Question: Can we find $\alpha_1, \alpha_2$ and $\alpha_3$ such that $(1,1,2,2) = \alpha_1 \cdot v_1 + \alpha_2 \cdot v_2 + \alpha_3 \cdot v_3$?

Answer: Yes but we need to solve the following system of linear equations:

$$1 = \alpha_1 - \alpha_2 + 0 \cdot \alpha_3$$
$$1 = \alpha_1 - \alpha_2 + 0 \cdot \alpha_3$$

$$2 = \alpha_1 + \alpha_3$$

$$2 = \alpha_1 + \alpha_3$$

Using what we have learned from chapter 1 to solve the above system of linear equations, we obtain:

$$\alpha_1 = \alpha_2 = \alpha_3 = 1$$

Hence, Yes: $(1,1,2,2) \in D$.

**The Second Way (Recommended):** We first need to find $dim(D)$, and then a basis for $D$. We have to write $v_1, v_2$ and $v_3$ as a matrix.

$$\begin{bmatrix} 1 & 1 & 1 & 1 \\ -1 & -1 & 0 & 0 \\ 0 & 0 & 1 & 1 \end{bmatrix}$$ Each point is a row-operation. We need to reduce this matrix to Semi-Reduced Matrix.

Now, we apply the Row-Reduction Method to get the Semi-Reduced Matrix as follows:

$$\begin{bmatrix} 1 & 1 & 1 & 1 \\ -1 & -1 & 0 & 0 \\ 0 & 0 & 1 & 1 \end{bmatrix} R_1 + R_2 \to R_2 \begin{bmatrix} \boxed{1} & 1 & 1 & 1 \\ 0 & 0 & \boxed{1} & 1 \\ 0 & 0 & 1 & 1 \end{bmatrix}$$

$$-R_2 + R_3 \to R_3 \begin{bmatrix} 1 & 1 & 1 & 1 \\ 0 & 0 & 1 & 1 \\ 0 & 0 & 0 & 0 \end{bmatrix}$$ This is a Semi-Reduced Matrix.

Since there is a zero-row in the Semi-Reduced Matrix, then these elements are dependent. Thus, $\dim(D) = 2$.

Thus, Basis for $D$ is $\{(1,1,1,1), (0,0,1,1)\}$, and

$D = Span\{(1,1,1,1), (0,0,1,1)\}$.

Now, we ask ourselves the following question:

Question: Can we find $\alpha_1, \alpha_2$ and $\alpha_3$ such that $(1,1,2,2) = \alpha_1 \cdot (1,1,1,1) + \alpha_2 \cdot (0,0,1,1)$?

Answer: Yes:

$$1 = \alpha_1$$

$$1 = \alpha_1$$

$$2 = \alpha_1 + \alpha_2$$

$$2 = \alpha_1 + \alpha_2$$

Thus, $\alpha_1 = \alpha_2 = \alpha_3 = 1$.

Hence, Yes: $(1,1,2,2) \in D$.

## 2.5 Exercises

1. Let $M = Span\{(1,-1,1), (-1,0,1), (0,-1,2)\}$

   a. Find a basis for $M$.

   b. Is $(1,-3,5) \in M$? Why?

2. Let $D = \{(x, y, z, t) \in \mathbb{R}^4 | x, y, z, t \in \mathbb{R},$
   $x + 2z + 3t = 0, and\ y - z + t = 0\}$

   a. Show that $D$ is a subspace of $\mathbb{R}^4$.

   b. Find a basis for $D$.

   c. Write $D$ as a $Span$.

3. Let $E = \{\begin{bmatrix} m \\ m+w \\ w \end{bmatrix} | m, w \in \mathbb{R}\}$

   a. Show that $E$ is a subspace of $\mathbb{R}^3$.

   b. Find a basis for $E$.

   c. Write $E$ as a $Span$.

4. Find a basis for the subspace $F$ where
$$F = \{\begin{bmatrix} w - s + 3u & 4w + 3s - 9u & 2w \\ 8w + 2s - 6u & 5w & 0 \end{bmatrix} | w, s, u \in \mathbb{R}\}$$

5. Determine the value(s) of $x$ such that the points $(1,0,5), (1,2,4), (1,4,x)$ are dependent.

6. Let $D = Span\{\begin{bmatrix} 2 & -1 \\ 0 & 1 \end{bmatrix}, \begin{bmatrix} 1 & -1 \\ 0 & 0.5 \end{bmatrix}, \begin{bmatrix} 5 & -3 \\ 0 & 2.5 \end{bmatrix}\}$

   a. Find $dim(D)$.

   b. Find a basis $A$ for $D$.

   c. Is $C = \begin{bmatrix} -2 & 0 \\ 0 & -1 \end{bmatrix} \in D$? Why?

   d. If the answer to part c is yes, then write $C$ as a linear combination of the elements in $A$. Otherwise, write the basis $A$ as a $Span$.

7. Let $K = Span\{(1, -1, 0), (2, -1, 0), (1, 0, 0)\}$.
Find $dim(K)$.

8. Find a basis for the subspace of $\mathbb{R}^4$ spanned by $\{(2,9, -2, 53), (-3, 2, 3, -2), (8, -3, -8, 17), (0, -3, 0, 15)\}$.

9. Does the $Span\{(-2, 1, 2), (2, 1, -1), (2, 3, 0)\}$ equal to $\mathbb{R}^3$?

# Chapter 3

# Homogeneous Systems

In this chapter, we introduce the homogeneous systems, and we discuss how they are related to what we have learned in chapter 2. We start with an introduction to null space and rank. Then, we study one of the most important topics in linear algebra which is linear transformation. At the end of this chapter we discuss how to find range and kernel, and their relation to sections 3.1 and 3.2.

## 3.1 Null Space and Rank

In this section, we first give an introduction to homogeneous systems, and we discuss how to find the null space and rank of homogeneous systems. In addition, we explain how to find row space and column space.

**Definition 3.1.1** Homogeneous System is a $m \times n$ system of linear equations that has all zero constants.

(i.e. the following is an example of homogeneous system): $\begin{cases} 2x_1 + x_2 - x_3 + x_4 = 0 \\ 3x_1 + 5x_2 + 3x_3 + 4x_4 = 0 \\ -x_2 + x_3 - x_4 = 0 \end{cases}$

Imagine we have the following solution to the homogeneous system: $x_1 = x_2 = x_3 = x_4 = 0$.

Then, this solution can be viewed as a point of $\mathbb{R}^n$ (here is $\mathbb{R}^4$): $(0,0,0,0)$

**Result 3.1.1** The solution of a homogeneous system $m \times n$ can be written as

$\{(a_1, a_2, a_3, a_4, \ldots, a_n | a_1, a_2, a_3, a_4, \ldots, a_n \in \mathbb{R}\}$.

**Result 3.1.2** All solutions of a homogeneous system $m \times n$ form a subset of $\mathbb{R}^n$, and it is equal to the number of variables.

**Result 3.1.3** Given a homogeneous system $m \times n$. We write it in the matrix-form: $C \begin{bmatrix} x_1 \\ x_2 \\ x_3 \\ \vdots \\ x_n \end{bmatrix} = \begin{bmatrix} 0 \\ 0 \\ 0 \\ \vdots \\ 0 \end{bmatrix}$ where $C$ is a coefficient. Then, the set of all solutions in this system is a subspace of $\mathbb{R}^n$.

**Proof of Result 3.1.3** We assume that $M_1 = (m_1, m_2, \ldots, m_n)$ and $W_1 = (w, w_2, \ldots, w_n)$ are two solutions to the above system. We will show that $M + W$ is a solution. We write them in the matrix-form:

$C \begin{bmatrix} m_1 \\ m_2 \\ m_3 \\ \vdots \\ m_n \end{bmatrix} = \begin{bmatrix} 0 \\ 0 \\ 0 \\ \vdots \\ 0 \end{bmatrix}$ and $\begin{bmatrix} w_1 \\ w_2 \\ w_3 \\ \vdots \\ w_n \end{bmatrix} = \begin{bmatrix} 0 \\ 0 \\ 0 \\ \vdots \\ 0 \end{bmatrix}$

Now, using algebra: $M + W = C\begin{bmatrix} m_1 \\ m_2 \\ m_3 \\ \vdots \\ m_n \end{bmatrix} + C\begin{bmatrix} w_1 \\ w_2 \\ w_3 \\ \vdots \\ w_n \end{bmatrix} = \begin{bmatrix} 0 \\ 0 \\ 0 \\ \vdots \\ 0 \end{bmatrix}$

By taking $C$ as a common factor, we obtain:

$$C\left(\begin{bmatrix} m_1 \\ m_2 \\ m_3 \\ \vdots \\ m_n \end{bmatrix} + \begin{bmatrix} w_1 \\ w_2 \\ w_3 \\ \vdots \\ w_n \end{bmatrix}\right) = \begin{bmatrix} 0 \\ 0 \\ 0 \\ \vdots \\ 0 \end{bmatrix}$$

$$C\begin{bmatrix} m_1 + w_1 \\ m_2 + w_2 \\ m_3 + w_3 \\ \vdots \\ m_n + w_n \end{bmatrix} = \begin{bmatrix} 0 \\ 0 \\ 0 \\ \vdots \\ 0 \end{bmatrix}$$

Thus, $M + W$ is a solution. □

**Fact 3.1.1** If $M_1 = (m_1, m_2, \ldots, m_n)$ is a solution, and $\alpha \in \mathbb{R}$, then $\alpha M = (\alpha m_1, \alpha m_2, \ldots, \alpha m_n)$ is a solution.

**Fact 3.1.2** The only system where the solutions form a vector space is the homogeneous system.

**Definition 3.1.2** Null Space of a matrix, say $A$ is a set of all solutions to the homogeneous system, and it is denoted by $Null(A)$ or $N(A)$.

**Definition 3.1.3** Rank of a matrix, say $A$ is the number of independent rows or columns of $A$, and it is denoted by $Rank(A)$.

**Definition 3.1.4** Row Space of a matrix, say $A$ is the $Span$ of independent rows of $A$, and it is denoted by $Row(A)$.

**Definition 3.1.5** Column Space of a matrix, say $A$ is the $Span$ of independent columns of $A$, and it is denoted by $Column(A)$.

**Example 3.1.1** Given the following $3 \times 5$ matrix:
$$A = \begin{bmatrix} 1 & -1 & 2 & 0 & -1 \\ 0 & 1 & 2 & 0 & 2 \\ 0 & 0 & 0 & 1 & 0 \end{bmatrix}.$$

a. Find $Null(A)$.

b. Find $dim(Null(A))$.

c. Rewrite $Null(A)$ as $Span$.

d. Find $Rank(A)$.

e. Find $Row(A)$.

**Solution: Part a:** To find the null space of $A$, we need to find the solution of $A$ as follows:

**Step 1:** Write the above matrix as an Augmented-Matrix, and make all constants' terms zeros.

$$\begin{pmatrix} 1 & -1 & 2 & 0 & -1 & | & 0 \\ 0 & 1 & 2 & 0 & 2 & | & 0 \\ 0 & 0 & 0 & 1 & 0 & | & 0 \end{pmatrix}$$

**Step 2:** Apply what we have learned from chapter 1 to solve systems of linear equations use Row-Operation Method.

$$\begin{pmatrix} 1 & -1 & 2 & 0 & -1 & | & 0 \\ 0 & 1 & 2 & 0 & 2 & | & 0 \\ 0 & 0 & 0 & 1 & 0 & | & 0 \end{pmatrix} R_2 + R_1 \to R_1$$

$$\begin{pmatrix} 1 & 0 & 4 & 0 & 1 & | & 0 \\ 0 & 1 & 2 & 0 & 2 & | & 0 \\ 0 & 0 & 0 & 1 & 0 & | & 0 \end{pmatrix}$$ This is a Completely-Reduced Matrix.

**Step 3:** Read the solution for the above system of linear equations after using Row-Operation.

$x_1 + 4x_3 + x_5 = 0$
$x_2 + 2x_3 + 2x_5 = 0$
$\phantom{x_2 + 2x_3 + 2}x_4 = 0$

Free variables are $x_3$ and $x_5$.

Assuming that $x_3, x_5 \in \mathbb{R}$. Then, the solution of the above homogeneous system is as follows:

$x_1 = -4x_3 - x_5$
$x_2 = -2x_3 - 2x_5$
$\phantom{x_2 = -2x_3 - 2}x_4 = 0$

Thus, according to definition 3.1.2,

$Null(A) = \{(-4x_3 - x_5, -2x_3 - 2x_5, x_3, 0, x_5) | x_3, x_5 \in \mathbb{R}\}$.

**Part b:** It is always true that

$dim(Null(A)) = dim(N(A)) = The\ Number\ of\ Free\ Variables$

Here, $dim(Null(A)) = 2$.

**Definition 3.1.6** The nullity of a matrix, say $A$ is the dimension of the null space of $A$, and it is denoted by $dim(Null(A))$ or $dim(N(A))$.

**Part c:** We first need to find a basis for $Null(A)$ as follows: To find a basis for $Null(A)$, we play a game called (ON-OFF GAME) with the free variables $x_3$ and $x_5$.

| $x_3$ | $x_5$ | Point |
|---|---|---|
| 1 | 0 | $(-4,-2,1,0,0)$ |
| 0 | 1 | $(-1,-2,0,0,1)$ |

The basis for $Null(A) = \{(-4,-2,1,0,0), (-1,-2,0,0,1)\}$.

Thus, $Null(A) = Span\{(-4,-2,1,0,0), (-1,-2,0,0,1)\}$.

**Part d:** To find the rank of matrix $A$, we just need to change matrix $A$ to the Semi-Reduced Matrix. We already did that in part a. Thus, $Rank(A) = 3$.

**Part e:** To find the row space of matrix $A$, we just need to write the $Span$ of independent rows. Thus, $Row(A) = Span\{(1,-1,2,0,-1), (0,1,2,0,2), (0,0,0,1,0)\}$.
It is also a subspace of $\mathbb{R}^5$.

**Result 3.1.4** Let $A$ be $m \times n$ matrix. Then,

$Rank(A) + dim(N(A)) = n = Number\ of\ Columns\ of\ A$.

**Result 3.1.5** Let $A$ be $m \times n$ matrix. The geometric meaning of $Row(A) = Span\{Independent\ Rows\}$ "lives" inside $\mathbb{R}^n$.

**Result 3.1.6** Let $A$ be $m \times n$ matrix. The geometric meaning of $Column(A) = Span\{Independent\ Columns\}$ "lives" inside $\mathbb{R}^m$.

**Result 3.1.7** Let $A$ be $m \times n$ matrix. Then, $Rank(A) = dim(Row(A)) = dim(Column(A))$.

**Example 3.1.2** Given the following $3 \times 5$ matrix:
$$B = \begin{bmatrix} 1 & 1 & 1 & 1 & 1 \\ -1 & -1 & -1 & 0 & 2 \\ 0 & 0 & 0 & 0 & 0 \end{bmatrix}.$$

a. Find $Row(B)$.

b. Find $Column(B)$.

c. Find $Rank(B)$.

**Solution: Part a:** To find the row space of $B$, we need to change matrix $B$ to the Semi-Reduced Matrix as follows:

$$\begin{bmatrix} 1 & 1 & 1 & 1 & 1 \\ -1 & -1 & -1 & 0 & 2 \\ 0 & 0 & 0 & 0 & 0 \end{bmatrix} \begin{matrix} R_1 + R_2 \to R_2 \\ R_1 + R_3 \to R_3 \end{matrix} \begin{bmatrix} \boxed{1} & 1 & 1 & 1 & 1 \\ 0 & 0 & 0 & \boxed{1} & 3 \\ 0 & 0 & 0 & 0 & 0 \end{bmatrix}$$

This is a Semi-Reduced Matrix. To find the row space of matrix $B$, we just need to write the $Span$ of independent rows. Thus, $Row(B) = Span\{(1,1,1,1,1), (0,0,0,1,3)\}$.

**Part b:** To find the column space of $B$, we need to change matrix $B$ to the Semi-Reduced Matrix. We already did that in part a. Now, we need to locate the columns in the Semi-Reduced Matrix of $B$ that contain the leaders, and then we should locate them to the original matrix $B$.

$$\begin{bmatrix} \boxed{1} & 1 & 1 & 1 & 1 \\ 0 & 0 & 0 & \boxed{1} & 3 \\ 0 & 0 & 0 & 0 & 0 \end{bmatrix} \text{Semi-Reduced Matrix}$$

$$\begin{bmatrix} \boxed{1} & 1 & 1 & \boxed{1} & 1 \\ -1 & -1 & -1 & 0 & 2 \\ 0 & 0 & 0 & 0 & 0 \end{bmatrix} \text{Matrix } B$$

Each remaining columns is a linear combination of the first and fourth columns.

Thus, $Column(B) = Span\{(1,-1,0), (1,0,0)\}$.

**Part c:** To find the rank of matrix $B$, we just need to change matrix $A$ to the Semi-Reduced Matrix. We already did that in part a. Thus,

$Rank(A) = dim\big(Row(B)\big) = dim(Column(B)) = 2$.

# 3.2 Linear Transformation

We start this section with an introduction to polynomials, and we explain how they are similar to $\mathbb{R}^n$ as vector spaces. At the end of this section we discuss a new concept called linear transformation.

Before discussing polynomials, we need to know the following mathematical facts:

**Fact 3.2.1** $\mathbb{R}^{n \times m} = \mathbb{R}_{n \times m} = M_{n \times m}(\mathbb{R})$ is a vector space.

**Fact 3.2.2** $\mathbb{R}^{2 \times 3}$ is equivalent to $\mathbb{R}^6$ as a vector space.

(i.e. $\begin{bmatrix} 1 & 2 & 3 \\ 0 & 1 & 1 \end{bmatrix}$ is equivalent to $(1,2,3,0,1,1)$ ).

**Fact 3.2.3** $\mathbb{R}^{3 \times 2}$ is equivalent to $\mathbb{R}^6$ as a vector space.

(i.e. $\begin{bmatrix} 1 & 2 \\ 3 & 0 \\ 1 & 1 \end{bmatrix}$ is equivalent to $(1,2,3,0,1,1)$ ).

After knowing the above facts, we introduce polynomials as follows:

$P_n = Set\ of\ all\ polynomials\ of\ degree < n.$

The algebraic expression of polynomials is in the following from: $a_n x^n + a_{n-1} x^{n-1} + \cdots + a_1 x^1 + a_0$

$a_n, a_{n-1}$ and $a_1$ are coefficients.

$n$ and $n - 1$ are exponents that must be positive integers whole numbers.

$a_0$ is a constant term.

The degree of polynomial is determined by the highest power (exponent).

We list the following examples of polynomials:

- $P_2$ = Set of all polynomials of degree < 2 (i.e. $3x + 2 \in P_2$, $0 \in P_2$, $10 \in P_2$, $\sqrt{3} \in P_2$ but $\sqrt{3}\sqrt{x} \notin P_2$).
- $P_4$ = Set of all polynomials of degree < 4 (i.e. $31x^2 + 4 \in P_4$).
- If $P(x) = 3$, then $deg(P(x)) = 0$.
- $\sqrt{x} + 3$ is not a polynomial.

**Result 3.2.1** $P_n$ is a vector space.

**Fact 3.2.4** $\mathbb{R}^{2\times 3} = M_{2\times 3}(\mathbb{R})$ as a vector space same as $\mathbb{R}^6$.

**Result 3.2.2** $P_n$ is a vector space, and it is the same as $\mathbb{R}^n$. (i.e. $a_0 + a_1 x^1 + \cdots + a_{n-1} x^{n-1} \leftrightarrow (a_0, a_1, \ldots, a_{n-1})$. Note: The above form is in an ascending order.

**Result 3.2.3** $dim(P_n) = n$.

**Fact 3.2.5** $P_3 = Span\{3\ Independent\ Polynomials, and\ Each\ of\ Degree < 3\}$. (i.e. $P_3 = Span\{1, x, x^2\}$).

**Example 3.2.1** Given the following polynomials: $3x^2 - 2, -5x, 6x^2 - 10x - 4$.

a. Are these polynomials independent?
b. Let $D = Span\{3x^2 - 2, -5x, 6x^2 - 10x - 4\}$. Find a basis for $D$.

**Solution: Part a:** We know that these polynomial live in $P_3$, and as a vector space $P_3$ is the same as $\mathbb{R}^3$. According to result 3.2.2, we need to make each polynomial equivalent to $\mathbb{R}^n$ as follows:

$3x^2 - 2 = -2 + 0x + 3x^2 \leftrightarrow (-2,0,3)$

$-5x = 0 - 5x + 0x^2 \leftrightarrow (0,-5,0)$

$6x^2 - 10x - 4 = -4 - 10x + 6x^2 \leftrightarrow (-4,-10,6)$

Now, we need to write these vectors as a matrix.

$\begin{bmatrix} -2 & 0 & 3 \\ 0 & -5 & 0 \\ -4 & -10 & 6 \end{bmatrix}$ Each point is a row-operation. We need to reduce this matrix to Semi-Reduced Matrix.

Then, we apply the Row-Reduction Method to get the Semi-Reduced Matrix as follows:

$\begin{bmatrix} -2 & 0 & 3 \\ 0 & -5 & 0 \\ -4 & -10 & 6 \end{bmatrix} \; -2R_1 + R_3 \to R_3 \; \begin{bmatrix} -2 & 0 & 3 \\ 0 & -5 & 0 \\ 0 & -10 & 0 \end{bmatrix}$

$-2R_2 + R_3 \to R_3 \; \begin{bmatrix} -2 & 0 & 3 \\ 0 & -5 & 0 \\ 0 & 0 & 0 \end{bmatrix}$ This is a Semi-Reduced Matrix.

Since there is a zero-row in the Semi-Reduced Matrix, then these elements are dependent. Thus, the answer to this question is NO.

**Part b:** Since there are only 2 vectors survived after checking for dependency in part a, then the basis for $D$ $(0,-5,0) \leftrightarrow -5x$.

**Result 3.2.4** Given $v_1, v_2, \ldots, v_k$ points in $\mathbb{R}^n$ where $k < n$. Choose one particular point, say $Q$, such that $Q = c_1 v_1 + c_2 v_2 + \cdots + c_k v_k$ where $c_1, c_2, \ldots, c_k$ are constants. If $c_1, c_2, \ldots, c_k$ are unique, then $v_1, v_2, \ldots, v_k$ are independent.

Note: The word "unique" in result 3.2.4 means that there is only one value for each of $c_1, c_2, \ldots, c_k$.

**Proof of Result 3.2.4** By using proof by contradiction, we assume that $v_1 = \alpha_2 v_2 + \alpha_3 v_3 + \cdots + \alpha_k v_k$ where $\alpha_2, \alpha_3, \ldots, \alpha_k$ are constants. Our assumption means that it is dependent. Using algebra, we obtain:

$Q = c_1 \alpha_2 v_2 + c_1 \alpha_3 v_3 + \cdots + c_1 \alpha_k v_k + c_2 v_2 + \cdots + c_k v_k$.

$Q = (c_1 \alpha_2 + c_2) v_2 + (c_1 \alpha_3 + c_3) v_3 + \cdots + (c_1 \alpha_k + c_k) v_k + 0 v_1$. Thus, none of them is a linear combination of the others which means that they are linearly independent. This is a contradiction. Therefore, our assumption that $v_1, v_2, \ldots,$ and $v_k$ were linearly dependent is false. Hence, $v_1, v_2, \ldots,$ and $v_k$ are linearly independent. ◻

**Result 3.2.5** Assume $v_1, v_2, \ldots, v_k$ are independent and $Q \in Span\{v_1, v_2, \ldots, v_k\}$. Then, there exists unique number $c_1, c_2, \ldots, c_k$ such that $Q = c_1 v_1 + c_2 v_2 + \cdots + c_k v_k$.

**Linear Transformation:**
**Definition 3.2.1** $T: V \to W$ where $V$ is a domain and $W$ is a co-domain. $T$ is a linear transformation if for every $v_1, v_2 \in V$ and $\alpha \in \mathbb{R}$, we have the following:
$T(\alpha v_1 + v_2) = \alpha T(v_1) + T(v_2)$.

**Example 3.2.2** Given $T: \mathbb{R}^2 \to \mathbb{R}^3$ where $\mathbb{R}^2$ is a domain and $\mathbb{R}^3$ is a co-domain. $T((a_1, a_2)) = (3a_1 + a_2, a_2, -a_1)$.

a. Find $T((1,1))$.
b. Find $T((1,0))$.
c. Show that $T$ is a linear transformation.

**Solution: Part a:** Since $T((a_1, a_2)) = (3a_1 + a_2, a_2, -a_1)$, then $a_1 = a_2 = 1$. Thus, $T((1,1)) = (3(1) + 1, 1, -1) = (4, 1, -1)$.

**Part b:** Since $T((a_1, a_2)) = (3a_1 + a_2, a_2, -a_1)$, then $a_1 = 1$ and $a_2 = 0$. Thus, $T((1,0)) = (3(1) + 0, 0, -1) = (3, 0, -1)$.

**Part c:** Proof: We assume that $v_1 = (a_1, a_2)$, $v_2 = (b_1, b_2)$, and $\alpha \in \mathbb{R}$. We will show that $T$ is a linear transformation. Using algebra, we start from the Left-Hand-Side (LHS):

$\alpha v_1 + v_2 = (\alpha a_1 + b_1, \alpha a_2 + b_2)$
$T(\alpha v_1 + v_2) = T((\alpha a_1 + b_1, \alpha a_2 + b_2))$
$T(\alpha v_1 + v_2) = (3\alpha a_1 + 3b_1 + \alpha a_2 + b_2, \alpha a_2 + b_2, -\alpha a_1 - b_1)$

Now, we start from the Right-Hand-Side (RHS):

$\alpha T(v_1) + T(v_2) = \alpha T(a_1, a_2) + T(b_1, b_2)$
$\alpha T(v_1) + T(v_2) = \alpha(3a_1 + a_2, a_2, -a_1) + (3b_1 + b_2, b_2, -b_1)$
$= (3\alpha a_1 + \alpha a_2, \alpha a_2, -\alpha a_1) + (3b_1 + b_2, b_2, -b_1)$
$= (3\alpha a_1 + \alpha a_2 + 3b_1 + b_2, \alpha a_2 + b_2, -\alpha a_1 - b_1)$

Thus, $T$ is a linear transformation. $\square$

**Result 3.2.6** Given $T: \mathbb{R}^n \to \mathbb{R}^m$. Then, $T((a_1, a_2, a_3, \ldots, a_n)) = $ Each coordinate is a linear combination of the $a_i$'s.

**Example 3.2.3** Given $T: \mathbb{R}^3 \to \mathbb{R}^4$ where $\mathbb{R}^3$ is a domain and $\mathbb{R}^4$ is a co-domain.

a. If $T((x_1, x_2, x_3)) = (-3x_3 + 6x_1, -10x_2, 13, -x_3)$, is $T$ a linear transformation?
b. If $T((x_1, x_2, x_3)) = (-3x_3 + 6x_1, -10x_2, 0, -x_3)$, is $T$ a linear transformation?

**Solution: Part a:** Since 13 is not a linear combination of $x_1, x_2$ and $x_3$. Thus, $T$ is not a linear transformation.

**Part b:** Since 0 is a linear combination of $x_1, x_2$ and $x_3$. Thus, $T$ is a linear transformation.

**Example 3.2.4** Given $T: \mathbb{R}^2 \to \mathbb{R}^3$ where $\mathbb{R}^2$ is a domain and $\mathbb{R}^3$ is a co-domain. If $T((a_1, a_2)) = (a_1^2 + a_2, -a_2)$, is $T$ a linear transformation?

**Solution:** Since $a_1^2 + a_2$ is not a linear combination of $a_1$ and $a_2$. Hence, $T$ is not a linear transformation.

**Example 3.2.5** Given $T: \mathbb{R} \to \mathbb{R}$. If $T(x) = 10x$, is $T$ a linear transformation?

**Solution:** Since it is a linear combination of $a_1$ such that $\alpha a_1 = 10x$. Hence, $T$ is a linear transformation.

**Example 3.2.6** Find the standard basis for $\mathbb{R}^2$.

**Solution:** The standard basis for $\mathbb{R}^2$ is the rows of $I_2$.

Since $I_2 = \begin{bmatrix} 1 & 0 \\ 0 & 1 \end{bmatrix}$, then the standard basis for $\mathbb{R}^2$ is $\{(1,0),(0,1)\}$.

**Example 3.2.7** Find the standard basis for $\mathbb{R}^3$.

**Solution:** The standard basis for $\mathbb{R}^3$ is the rows of $I_3$.

Since $I_3 = \begin{bmatrix} 1 & 0 & 0 \\ 0 & 1 & 0 \\ 0 & 0 & 1 \end{bmatrix}$, then the standard basis for $\mathbb{R}^3$ is $\{(1,0,0),(0,1,0),(0,0,1)\}$.

**Example 3.2.8** Find the standard basis for $P_3$.

**Solution:** The standard basis for $P_3$ is $\{1, x, x^2\}$.

**Example 3.2.9** Find the standard basis for $P_4$.

**Solution:** The standard basis for $P_4$ is $\{1, x, x^2, x^3\}$.

**Example 3.2.10** Find the standard basis for $\mathbb{R}_{2\times 2} = M_{2\times 2}(\mathbb{R})$.

**Solution:** The standard basis for $\mathbb{R}_{2\times 2} = M_{2\times 2}(\mathbb{R})$ is $\{\begin{bmatrix} 1 & 0 \\ 0 & 0 \end{bmatrix}, \begin{bmatrix} 0 & 1 \\ 0 & 0 \end{bmatrix}, \begin{bmatrix} 0 & 0 \\ 1 & 0 \end{bmatrix}, \begin{bmatrix} 0 & 0 \\ 0 & 1 \end{bmatrix}\}$ because $\mathbb{R}_{2\times 2} = M_{2\times 2}(\mathbb{R}) = \mathbb{R}^4$ as a vector space where standard basis for $\mathbb{R}_{2\times 2} = M_{2\times 2}(\mathbb{R})$ is the rows of $I_4 = \begin{bmatrix} 1 & 0 & 0 & 0 \\ 0 & 1 & 0 & 0 \\ 0 & 0 & 1 & 0 \\ 0 & 0 & 0 & 1 \end{bmatrix}$ that are represented by $2 \times 2$ matrices.

**Example 3.2.11** Let $T: \mathbb{R}^2 \to \mathbb{R}^3$ be a linear transformation such that
$T(2,0) = (0,1,4)$

$T(-1,1) = (2,1,5)$

Find $T(3,5)$.

**Solution:** The given points are $(2,0)$ and $(-1,1)$. These two points are independent because of the following:

$$\begin{bmatrix} 2 & 0 \\ -1 & 1 \end{bmatrix} \frac{1}{2}R_1 + R_2 \to R_2 \begin{bmatrix} 2 & 0 \\ 0 & 1 \end{bmatrix}$$

Every point in $\mathbb{R}^2$ is a linear combination of $(2,0)$ and $(-1,1)$. There exists unique numbers $c_1$ and $c_2$ such that $(3,5) = c_1(2,0) + c_2(-1,1)$.
$3 = 2c_1 - c_2$
$5 = c_2$
Now, we substitute $c_2 = 5$ in $3 = 2c_1 - c_2$, we obtain:
$3 = 2c_1 - 5$
$c_1 = 4$
Hence, $(3,5) = 4(2,0) + 5(-1,1)$.
$T(3,5) = T(4(2,0) + 5(-1,1))$
$T(3,5) = 4T(2,0) + 5T(-1,1)$
$T(3,5) = 4(0,1,4) + 5(2,1,5) = (10,9,41)$
Thus, $T(3,5) = (10,9,41)$.

**Example 3.2.12** Let $T: \mathbb{R} \to \mathbb{R}$ be a linear transformation such that $T(1) = 3$. Find $T(5)$.

**Solution:** Since it is a linear transformation, then $T(5) = T(5 \cdot 1) = 5T(1) = 5(3) = 15$. If it is not a linear transformation, then it is impossible to find $T(5)$.

## 3.3 Kernel and Range

In this section, we discuss how to find the standard matrix representation, and we give examples of how to find kernel and range.

**Definition 3.3.1** Given $T: \mathbb{R}^n \to \mathbb{R}^m$ where $\mathbb{R}^n$ is a domain and $\mathbb{R}^m$ is a co-domain. Then, Standard Matrix Representation is a $m \times n$ matrix. This means that it is $dim(Co - Domain) \times dim(Domain)$ matrix.

**Definition 3.3.2** Given $T: \mathbb{R}^n \to \mathbb{R}^m$ where $\mathbb{R}^n$ is a domain and $\mathbb{R}^m$ is a co-domain. Kernel is a set of all points in the domain that have image which equals to the origin point, and it is denoted by $Ker(T)$. This means that $Ker(T) = Null\ Space\ of\ T$.

**Definition 3.3.3** Range is the column space of standard matrix representation, and it is denoted by $Range(T)$.

**Example 3.3.1** Given $T: \mathbb{R}^3 \to \mathbb{R}^4$ where $\mathbb{R}^3$ is a domain and $\mathbb{R}^4$ is a co-domain.
$T\big((x_1, x_2, x_3)\big) = (-5x_1, 2x_2 + x_3, -x_1, 0)$
   a. Find the Standard Matrix Representation.
   b. Find $T((3,2,1))$.
   c. Find $Ker(T)$.
   d. Find $Range(T)$.

**Solution: Part a:** According to definition 3.3.1, the Standard Matrix Representation, let's call it $M$, here is

$4 \times 3$. We know from section 3.2 that the standard basis for domain (here is $\mathbb{R}^3$) is $\{(1,0,0),(0,1,0),(0,0,1)\}$. We assume the following:
$$v_1 = (1,0,0)$$
$$v_2 = (0,1,0)$$
$$v_3 = (0,0,1)$$

Now, we substitute each point of the standard basis for domain in $T((x_1, x_2, x_3)) = (-5x_1, 2x_2 + x_3, -x_1, 0)$ as follows:

$T((1,0,0)) = (-5,0,-1,0)$
$T((0,1,0)) = (0,2,0,0)$
$T((0,0,1)) = (0,1,0,0)$

Our goal is to find $M$ so that $T((x_1, x_2, x_3)) = M \begin{bmatrix} x_1 \\ x_2 \\ x_3 \end{bmatrix}$.

$M = \begin{bmatrix} -5 & 0 & 0 \\ 0 & 2 & 1 \\ -1 & 0 & 0 \\ 0 & 0 & 0 \end{bmatrix}$ This is the Standard Matrix

Representation. The first, second and third columns represent $T(v_1), T(v_2)$ and $T(v_3)$.

**Part b:** Since $((x_1, x_2, x_3)) = M \begin{bmatrix} x_1 \\ x_2 \\ x_3 \end{bmatrix}$, then

$T((3,2,1)) = \begin{bmatrix} -5 & 0 & 0 \\ 0 & 2 & 1 \\ -1 & 0 & 0 \\ 0 & 0 & 0 \end{bmatrix} \begin{bmatrix} 3 \\ 2 \\ 1 \end{bmatrix}$

$$T((3,2,1)) = 3 \cdot \begin{bmatrix} -5 \\ 0 \\ -1 \\ 0 \end{bmatrix} + 2 \cdot \begin{bmatrix} 0 \\ 2 \\ 0 \\ 0 \end{bmatrix} + 1 \cdot \begin{bmatrix} 0 \\ 1 \\ 0 \\ 0 \end{bmatrix} = \begin{bmatrix} -15 \\ 5 \\ -3 \\ 0 \end{bmatrix}$$

$\begin{bmatrix} -15 \\ 5 \\ -3 \\ 0 \end{bmatrix}$ is equivalent to $(-15,5,-3,0)$. This lives in the co-domain. Thus, $T((3,2,1)) = (-15,5,-3,0)$.

**Part c:** According to definition 3.3.2, $Ker(T)$ is a set of all points in the domain that have image= $(0,0,0,0)$. Hence, $T((x_1, x_2, x_3)) = (0,0,0,0)$. This means the following: $M \begin{bmatrix} x_1 \\ x_2 \\ x_3 \end{bmatrix} = \begin{bmatrix} 0 \\ 0 \\ 0 \\ 0 \end{bmatrix}$

$$\begin{bmatrix} -5 & 0 & 0 \\ 0 & 2 & 1 \\ -1 & 0 & 0 \\ 0 & 0 & 0 \end{bmatrix} \begin{bmatrix} x_1 \\ x_2 \\ x_3 \end{bmatrix} = \begin{bmatrix} 0 \\ 0 \\ 0 \\ 0 \end{bmatrix}$$

Since $Ker(T) = Null(M)$, then we need to find $N(M)$ as follows:

$$\begin{pmatrix} -5 & 0 & 0 & | & 0 \\ 0 & 2 & 1 & | & 0 \\ -1 & 0 & 0 & | & 0 \\ 0 & 0 & 0 & | & 0 \end{pmatrix} -\frac{1}{5}R_1 \begin{pmatrix} 1 & 0 & 0 & | & 0 \\ 0 & 2 & 1 & | & 0 \\ -1 & 0 & 0 & | & 0 \\ 0 & 0 & 0 & | & 0 \end{pmatrix} R_1 + R_3 \to R_3$$

$$\begin{pmatrix} 1 & 0 & 0 & | & 0 \\ 0 & 2 & 1 & | & 0 \\ 0 & 0 & 0 & | & 0 \\ 0 & 0 & 0 & | & 0 \end{pmatrix} \frac{1}{2}R_2 \begin{pmatrix} 1 & 0 & 0 & | & 0 \\ 0 & 1 & 0.5 & | & 0 \\ 0 & 0 & 0 & | & 0 \\ 0 & 0 & 0 & | & 0 \end{pmatrix}$$ This is a Completely-Reduced Matrix. Now, we need to read the above matrix as follows:

$x_1 = 0$

$x_2 + \frac{1}{2}x_3 = 0$

$0 = 0$

$0 = 0$

To write the solution, we need to assume that $x_3 \in \mathbb{R}$ (*Free Variable*).

Hence, $x_1 = 0$ and $x_2 = -\frac{1}{2}x_3$.

$N(M) = \{(0, -\frac{1}{2}x_3, x_3) | x_3 \in \mathbb{R}\}$.

By letting $x_3 = 1$, we obtain:

$Nullity(M) = Number\ of\ Free\ Variables = 1$, and

$Basis = \{(0, -\frac{1}{2}, 1)\}$

Thus, $Ker(T) = N(M) = Span\{(0, -\frac{1}{2}, 1)\}$.

**Part d:** According to definition 3.3.3, $Range(T)$ is the column space of $M$. Now, we need to locate the columns in the Completely-Reduced Matrix in part c that contain the leaders, and then we should locate them to the original matrix as follows:

$$\begin{pmatrix} \boxed{1} & 0 & 0 & | & 0 \\ 0 & \boxed{1} & 0.5 & | & 0 \\ 0 & 0 & 0 & | & 0 \\ 0 & 0 & 0 & | & 0 \end{pmatrix}$$ Completely-Reduced Matrix

$$\begin{pmatrix} -5 & 0 & 0 & | & 0 \\ 0 & 2 & 1 & | & 0 \\ -1 & 0 & 0 & | & 0 \\ 0 & 0 & 0 & | & 0 \end{pmatrix}$$ Orignial Matrix

Thus, $Range(T) = Span\{(-5,0,-1,0), (0,2,0,0)\}$.

**Result 3.3.1** Given $T: \mathbb{R}^n \to \mathbb{R}^m$ where $\mathbb{R}^n$ is a domain and $\mathbb{R}^m$ is a co-domain. Let $M$ be a standard matrix representation. Then,
$Range(T) = Span\{Independent\ Columns of\ M\}$.

**Result 3.3.2** Given $T: \mathbb{R}^n \to \mathbb{R}^m$ where $\mathbb{R}^n$ is a domain and $\mathbb{R}^m$ is a co-domain. Let $M$ be a standard matrix representation. Then, $dim(Range(T)) = Rank(M) =$ Number of Independent Columns.

**Result 3.3.3** Given $T: \mathbb{R}^n \to \mathbb{R}^m$ where $\mathbb{R}^n$ is a domain and $\mathbb{R}^m$ is a co-domain. Let $M$ be a standard matrix representation. Then, $dim(Ker(T)) = Nullity(M)$.

**Result 3.3.4** Given $T: \mathbb{R}^n \to \mathbb{R}^m$ where $\mathbb{R}^n$ is a domain and $\mathbb{R}^m$ is a co-domain. Let $M$ be a standard matrix representation. Then, $dim(Range(T)) + dim(Ker(T)) = dim(Domain)$.

**Example 3.3.2** Given $T: P_2 \to \mathbb{R}$. $T(f(x)) = \int_0^1 f(x)dx$ is a linear transformation.
  a. Find $T(2x - 1)$.
  b. Find $Ker(T)$.
  c. Find $Range(T)$.

**Solution:**
**Part a:** $T(2x - 1) = \int_0^1 (2x - 1)dx = x^2 - x \Big|_{x=0}^{x=1} = 0$.
**Part b:** To find $Ker(T)$, we set equation of $T = 0$, and $f(x) = a_0 + a_1 x \in P_2$.

Thus, $T(f(x)) = \int_0^1 (a_0 + a_1 x)dx = a_0 x + \frac{a_1}{2}x^2 \Big|_{x=0}^{x=1} = 0$

$a_0 + \frac{a_1}{2} - 0 = 0$

$a_0 = -\frac{a_1}{2}$

Hence, $Ker(T) = \{-\frac{a_1}{2} + a_1 x | a_1 \in \mathbb{R}\}$. We also know that $dim(Ker(T)) = 1$ because there is one free variable. In addition, we can also find basis by letting $a_1$ be any real number not equal to zero, say $a_1 = 1$, as follows:

$Basis = \{-\frac{1}{2} + x\}$

Thus, $Ker(T) = Span\{-\frac{1}{2} + x\}$.

**Part c:** It is very easy to find range here. $Range(T) = \mathbb{R}$ because we linearly transform from a second degree polynomial to a real number. For example, if we linearly transform from a third degree polynomial to a second degree polynomial, then the range will be $P_2$.

## 3.4 Exercises

1. Given $T: \mathbb{R}^3 \to \mathbb{R}_{2\times 2}$ such that

   $T((x_1, x_2, x_3)) = \begin{bmatrix} x_1 & x_1 \\ x_3 & x_2 \end{bmatrix}$ is a linear transformation.

   a. Find the standard matrix representation of $T$.

   b. Find $Ker(T)$.

   c. Find a basis for $Range(T)$ and write $Range(T)$ as a $Span$.

2. Let $T: \mathbb{R}^2 \to \mathbb{R}$ be a linear transformation such that $T(1,0) = 4$, $T(2, -2) = 2$. Find the standard matrix representation of $T$.

3. Given $T: P_3 \to \mathbb{R}$ such that $T(a + bx + cx^2) = \int_0^1 (a + bx + cx^2)dx$. Find $Ker(T)$.

4. Given $T: \mathbb{R}^4 \to \mathbb{R}^3$ is a linear transformation such that $T(x, y, z, w) = (x + y + z - 2w, -2w, w)$.

   a. Find the standard matrix representation of $T$.
   b. Find $dim(Ker(T))$.
   c. Find $Range(T)$.

5. Given $T: P_4 \to \mathbb{R}_{2\times 2}$ such that
$T(g(x)) = \begin{bmatrix} g(-1) & g(0) \\ g(-1) & g(0) \end{bmatrix}$ is a linear transformation.

   a. Find the standard matrix representation of $T$.
   b. Find $Ker(T)$ and write $Ker(T)$ as a $Span$.
   c. Find a basis for $Range(T)$ and write $Range(T)$ as a $Span$.

6. Given $T: P_3 \to \mathbb{R}$ is a linear transformation such that $T(1) = 6, T(x^2 + x) = -5$, and $T(x^2 + 2x + 1) = 4$.

   a. Find $T(x), T(x^2)$ and $T(5x^2 + 3x + 8)$.
   b. Find the standard matrix representation of $T$.
   c. Find $Ker(T)$ and write $Ker(T)$ as a $Span$.

# Chapter 4
# Characteristic Equation of Matrix

In this chapter, we discuss how to find eigenvalues and eigenvectors of a matrix. Then, we introduce the diagonalizable matrix, and we explain how to determine whether a matrix is diagonalizable or not. At the end of this chapter we discuss how to find diagonal matrix and invertible matrix.

## 4.1 Eigenvalues and Eigenvectors

In this section, we give an example explaining the steps for finding eigenvalues and eigenvectors.

**Example 4.1.1** Given $n \times n$ matrix, say $A$, Let $\alpha \in \mathbb{R}$. Can we find a point $Q$ in $\mathbb{R}^n$,
$Q = (a_1, a_2, \dots, a_n)$, such that $Q \neq (0,0,0,0, \dots, 0)$, and
$$A_{n \times n} Q = A \begin{bmatrix} a_1 \\ a_2 \\ a_3 \\ \vdots \\ a_n \end{bmatrix} = \alpha \begin{bmatrix} a_1 \\ a_2 \\ a_3 \\ \vdots \\ a_n \end{bmatrix}?$$

**Solution:** If such $\alpha$ and such $Q$ exist, we say $\alpha$ is an eigenvalue of $A$, and we say $Q$ is an eigenvector of $A$ corresponds to the eigenvalue $\alpha$.

$$A \begin{bmatrix} a_1 \\ a_2 \\ a_3 \\ \vdots \\ a_n \end{bmatrix} = \alpha \begin{bmatrix} a_1 \\ a_2 \\ a_3 \\ \vdots \\ a_n \end{bmatrix}$$

$$A \begin{bmatrix} a_1 \\ a_2 \\ a_3 \\ \vdots \\ a_n \end{bmatrix} - \alpha \begin{bmatrix} a_1 \\ a_2 \\ a_3 \\ \vdots \\ a_n \end{bmatrix} = \begin{bmatrix} 0 \\ 0 \\ 0 \\ \vdots \\ 0 \end{bmatrix}$$

$$[A - \alpha I_n] \begin{bmatrix} a_1 \\ a_2 \\ a_3 \\ \vdots \\ a_n \end{bmatrix} = \begin{bmatrix} 0 \\ 0 \\ 0 \\ \vdots \\ 0 \end{bmatrix}$$

We conclude that such $\alpha$ and ($Q \neq$ origin )exist if and only if $det(A - \alpha I_n) = 0$.

Note: $det(A - \alpha I_n)$ is called characteristic polynomial of $A$.

**Example 4.1.2** Given the following $3 \times 3$ matrix:
$$A = \begin{bmatrix} 2 & 0 & 1 \\ 0 & 1 & -2 \\ 0 & 0 & -1 \end{bmatrix}$$
a. Find all eigenvalues of $A$.
b. For each eigenvalue, find the corresponding eigenspace.

**Solution: Part a:** To find all eigenvalues, we set characteristic polynomial of $A = 0$. This means that $det(A - \alpha I_n) = 0$. Then, we do the following:

$det(A - \alpha I_n) = 0$

$det\left(\begin{bmatrix} 2 & 0 & 1 \\ 0 & 1 & -2 \\ 0 & 0 & -1 \end{bmatrix} - \alpha \begin{bmatrix} 1 & 0 & 0 \\ 0 & 1 & 0 \\ 0 & 0 & 1 \end{bmatrix}\right) = 0$

$det\left(\begin{bmatrix} 2 & 0 & 1 \\ 0 & 1 & -2 \\ 0 & 0 & -1 \end{bmatrix} - \begin{bmatrix} \alpha & 0 & 0 \\ 0 & \alpha & 0 \\ 0 & 0 & \alpha \end{bmatrix}\right) = 0$

$det\left(\begin{bmatrix} 2-\alpha & 0 & 1 \\ 0 & 1-\alpha & -2 \\ 0 & 0 & -1-\alpha \end{bmatrix}\right) = 0$

$det\left(\begin{bmatrix} 2-\alpha & 0 & 1 \\ 0 & 1-\alpha & -2 \\ 0 & 0 & -1-\alpha \end{bmatrix}\right) = (2-\alpha) \cdot (1-\alpha) \cdot (-1-\alpha) = 0$

Thus, the eigenvalues of $A$ are: $\begin{array}{l} \alpha = 2 \\ \alpha = 1 \\ \alpha = -1 \end{array}$

**Part b:** Since the first eigenvalue is 2, then the eigenspace of $A$ that corresponds to the eigenvalue 2 equals to the set of all points in $\mathbb{R}_3$ such that $(point\ A) = \alpha \cdot point =$ set of all eigenvectors of $A$ that correspond to $\alpha = (2 + origin\ point)$. To find eigenspace that corresponds to $\alpha = 2$, say $E_2$, we need to find $E_2 = N(A - 2I_3)$. Note: $N$ here represents $Null$. Now, we do the following:

$E_2 = N(A - 2I_3) = N\left(\begin{bmatrix} 2 & 0 & 1 \\ 0 & 1 & -2 \\ 0 & 0 & -1 \end{bmatrix} - 2\begin{bmatrix} 1 & 0 & 0 \\ 0 & 1 & 0 \\ 0 & 0 & 1 \end{bmatrix}\right)$

$E_2 = N(A - 2I_3) = N\left(\begin{bmatrix} 2-2 & 0 & 1 \\ 0 & 1-2 & -2 \\ 0 & 0 & -1-2 \end{bmatrix}\right)$

$E_2 = N(A - 2I_3) = N\begin{bmatrix} 0 & 0 & 1 \\ 0 & -1 & -2 \\ 0 & 0 & -3 \end{bmatrix}$

Thus, we form the augmented-matrix for the homogeneous system as follows:
$$\begin{pmatrix} 0 & 0 & 1 & | & 0 \\ 0 & -1 & -2 & | & 0 \\ 0 & 0 & -3 & | & 0 \end{pmatrix}$$

It is very easy to solve it: $\begin{matrix} a_3 = 0 \\ a_2 = -2a_3 = 0 \\ a_3 = 0 \end{matrix}$

$a_1 \in \mathbb{R}$ (Free-Variable), and $a_2 = a_3 = 0$.
$E_2 = \{(a_1, 0, 0) | a_1 \in \mathbb{R}\}$. Now, we can select any value for $a_1 \neq 0$. Let's select $a_1 = 1$.
We can also find $dim(E_2) = 1$.
Hence, $E_2 = Span\{(1,0,0)\}$. Similarly, we can find the eigenspaces that correspond to the other eigenvalues: $\alpha = 1$ and $\alpha = -1$. To find eigenspace that corresponds to $\alpha = 1$, say $E_1$, we need to find $E_1 = N(A - 1I_3)$. Note: $N$ here represents $Null$. Now, we do the following:

$$E_1 = N(A - 1I_3) = N\left(\begin{bmatrix} 2 & 0 & 1 \\ 0 & 1 & -2 \\ 0 & 0 & -1 \end{bmatrix} - 1\begin{bmatrix} 1 & 0 & 0 \\ 0 & 1 & 0 \\ 0 & 0 & 1 \end{bmatrix}\right)$$

$$E_1 = N(A - 1I_3) = N\left(\begin{bmatrix} 2-1 & 0 & 1 \\ 0 & 1-1 & -2 \\ 0 & 0 & -1-1 \end{bmatrix}\right)$$

$$E_1 = N(A - 1I_3) = N\begin{bmatrix} 1 & 0 & 1 \\ 0 & 0 & -2 \\ 0 & 0 & -2 \end{bmatrix}$$

Thus, as we did previously, we form the augmented-matrix for the homogeneous system as follows:
$$\begin{pmatrix} 1 & 0 & 1 & | & 0 \\ 0 & 0 & -2 & | & 0 \\ 0 & 0 & -2 & | & 0 \end{pmatrix}$$

It is very easy to solve it: $\begin{matrix} a_1 = -a_3 \\ a_3 = 0 \end{matrix}$

$a_2 \in \mathbb{R}$ (Free-Variable), and $a_1 = a_3 = 0$.
$E_1 = \{(0, a_2, 0) | a_2 \in \mathbb{R}\}$. Now, we can select any value for $a_2 \neq 0$. Let's select $a_2 = 1$.
We can also find $dim(E_1) = 1$.
Hence, $E_1 = Span\{(0,1,0)\}$. Finally, to find eigenspace that corresponds to $\alpha = -1$, say $E_{-1}$, we need to find $E_{-1} = N(A + 1I_3)$. Note: $N$ here represents $Null$. Now, we do the following:

$$E_{-1} = N(A + 1I_3) = N\left(\begin{bmatrix} 2 & 0 & 1 \\ 0 & 1 & -2 \\ 0 & 0 & -1 \end{bmatrix} + 1\begin{bmatrix} 1 & 0 & 0 \\ 0 & 1 & 0 \\ 0 & 0 & 1 \end{bmatrix}\right)$$

$$E_{-1} = N(A + 1I_3) = N\left(\begin{bmatrix} 2+1 & 0 & 1 \\ 0 & 1+1 & -2 \\ 0 & 0 & -1+1 \end{bmatrix}\right)$$

$$E_{-1} = N(A + 1I_3) = N\begin{bmatrix} 3 & 0 & 1 \\ 0 & 2 & -2 \\ 0 & 0 & 0 \end{bmatrix}$$

Thus, as we did previously, we form the augmented-matrix for the homogeneous system as follows:
$$\left(\begin{array}{ccc|c} 3 & 0 & 1 & 0 \\ 0 & 2 & -2 & 0 \\ 0 & 0 & 0 & 0 \end{array}\right)$$

It is very easy to solve it: $\begin{array}{l} a_1 = -\frac{1}{3}a_3 \\ a_2 = a_3 \\ 0 = 0 \end{array}$

$a_3 \in \mathbb{R}$ (Free-Variable).
$E_{-1} = \{(-\frac{1}{3}a_3, a_3, a_3) | a_3 \in \mathbb{R}\}$. Now, we can select any value for $a_3 \neq 0$. Let's select $a_3 = 1$.
We can also find $dim(E_{-1}) = 1$.
Hence, $E_{-1} = Span\{(-\frac{1}{3}, 1, 1)\}$.

## 4.2 Diagonalizable Matrix

In this section, we explain the concept of diagonalization, and we give some examples explaining it, and how to find the diagonal and invertible matrices.

**Definition 4.2.1** $A$ is diagonalizable if there exists an invertible matrix $L$, and a diagonal matrix $D$ such that $A = LDL^{-1}$.

**Result 4.2.1** $A$ is $n \times n$ diagonalizable matrix if and only if $det(A - \alpha I_n)$ is written as multiplication of linear equation, say $det(A - \alpha I_n) =$ some constants $(c_1 - \alpha) \cdot (c_2 - \alpha) \cdot \ldots \cdot (c_k - \alpha)$, and $dim(E_{c_i}) = n_i$ (Multiplicity of the Eigenvalues) for $i \leq i \leq k$.

**Example 4.2.1** Assume $A$ is $5 \times 5$ matrix, and $det(A - \alpha I_5) = (3 - \alpha)^2 \cdot (-2 - \alpha) \cdot (4 - \alpha)$, and the eigenvalues of $A$ are $\begin{matrix} \alpha = 3 \\ \alpha = -2 \\ \alpha = 4 \end{matrix}$

Given $dim(E_3) = 1, dim(E_{-2}) = 1$ and $dim(E_4) = 2$.

Is $A$ diagonalizable matrix?

**Solution:** It is not a diagonalizable matrix because $dim(E_3) = 1$, and it must be equal to 2 instead of 1.

**Example 4.2.2** Assume $A$ is $4 \times 4$ matrix, and $det(A - \alpha I_4) = (2 - \alpha)^3 \cdot (3 - \alpha)$, and given $dim(E_2) = 3$ and $dim(E_3) = 1$.

Is $A$ diagonalizable matrix?

**Solution:** According result 4.2.1, it is a diagonalizable matrix.

**Example 4.2.3** Given the following $3 \times 3$ matrix:
$A = \begin{bmatrix} 2 & 0 & 1 \\ 0 & 1 & -2 \\ 0 & 0 & -1 \end{bmatrix}$. Use example 4.1.2 from section 4.1 to answer the following questions:

a. Is $A$ diagonalizable matrix? If yes, find a diagonal matrix $D$, and invertible matrix $L$ such that $A = LDL^{-1}$.
b. Find $A^6$ such that $A = LDL^{-1}$.
c. Find $A^{300}$ such that $A = LDL^{-1}$.

**Solution: Part a:** From example 4.1.2, we found the following:
$E_2 = Span\{(1,0,0)\}$
$E_1 = Span\{(0,1,0)\}$
$E_{-1} = Span\{(-\frac{1}{3}, 1, 1)\}$

According to result 4.2.1, $A$ is a diagonalizable matrix. Now, we need to find a diagonal matrix $D$, and invertible matrix $L$ such that $A = LDL^{-1}$. To find a diagonal matrix $D$, we create a $3 \times 3$ matrix, and we put eigenvalues on the main diagonal with repetition if there is a repetition, and all other elements are zeros. Hence, the diagonal matrix $D$ is as follows:

$$D = \begin{bmatrix} 2 & 0 & 0 \\ 0 & -1 & 0 \\ 0 & 0 & 1 \end{bmatrix}$$

To find an invertible matrix $L$, we create a $3 \times 3$ matrix like the one above but each eigenvalue above is represented by a column of the eigenspace that corresponds to that eigenvalue as follows:

$$L = \begin{bmatrix} 1 & -\frac{1}{3} & 0 \\ 0 & 1 & 1 \\ 0 & 1 & 0 \end{bmatrix}$$

Thus, $A = \begin{bmatrix} 2 & 0 & 1 \\ 0 & 1 & -2 \\ 0 & 0 & -1 \end{bmatrix} = LDL^{-1}$

$$= \begin{bmatrix} 1 & -\frac{1}{3} & 0 \\ 0 & 1 & 1 \\ 0 & 1 & 0 \end{bmatrix} \begin{bmatrix} 2 & 0 & 0 \\ 0 & -1 & 0 \\ 0 & 0 & 1 \end{bmatrix} \begin{bmatrix} 1 & -\frac{1}{3} & 0 \\ 0 & 1 & 1 \\ 0 & 1 & 0 \end{bmatrix}^{-1}$$

**Part b:** To find $A^6$ such that $A = LDL^{-1}$, we do the following steps:

$A = LDL^{-1}$
$A^6 = (LDL^{-1}) \cdot (LDL^{-1}) \cdot \ldots \cdot (LDL^{-1})$
$A^6 = (LD^2L^{-1}) \cdot \ldots \cdot (LDL^{-1})$
$A^6 = LD^6L^{-1}$

Thus, $A^6 = LD^6L^{-1} = L\begin{bmatrix} 2 & 0 & 0 \\ 0 & -1 & 0 \\ 0 & 0 & 1 \end{bmatrix}^6 L^{-1} =$

$L\begin{bmatrix} 2^6 & 0 & 0 \\ 0 & (-1)^6 & 0 \\ 0 & 0 & 1^6 \end{bmatrix} L^{-1} = L\begin{bmatrix} 64 & 0 & 0 \\ 0 & 1 & 0 \\ 0 & 0 & 1 \end{bmatrix} L^{-1}$

**Part c:** To find $A^{300}$ such that $A = LDL^{-1}$, we do the following steps:

$A = LDL^{-1}$
$A^{300} = (LDL^{-1}) \cdot (LDL^{-1}) \cdot \ldots \cdot (LDL^{-1})$
$A^{300} = (LD^2L^{-1}) \cdot \ldots \cdot (LDL^{-1})$
$A^{300} = LD^{300}L^{-1}$

Thus, $A^{300} = LD^{300}L^{-1} = L\begin{bmatrix} 2 & 0 & 0 \\ 0 & -1 & 0 \\ 0 & 0 & 1 \end{bmatrix}^{300} L^{-1} =$

$L\begin{bmatrix} 2^{300} & 0 & 0 \\ 0 & (-1)^{300} & 0 \\ 0 & 0 & 1^{300} \end{bmatrix} L^{-1} = L\begin{bmatrix} 2^{300} & 0 & 0 \\ 0 & 1 & 0 \\ 0 & 0 & 1 \end{bmatrix} L^{-1}$

## 4.3 Exercises

1. Given the following $4 \times 4$ matrix:
$C = \begin{bmatrix} 1 & 0 & 0 & 0 \\ 0 & 1 & 1 & 1 \\ 0 & 0 & -1 & 1 \\ 0 & 0 & 0 & -1 \end{bmatrix}$. Is $C$ diagonalizable matrix?

2. Assume that $A$ is $5 \times 5$ diagonalizable matrix, and given the following:
$E_3 = Span\{(2,1,0,0,1), (0,1,0,1,1), (0,0,2,2,0)\}$

$dim(E_3) = 3$
$E_2 = Span\{(0,0,0,1,1), (0,0,0,0,10)\}$

    a. Find the characteristic polynomial of $A$.
    b. Find a diagonal matrix $D$ such that $A = LDL^{-1}$.
    c. Find an invertible matrix $L$ such that $A = LDL^{-1}$.

3. Assume that $W$ is $3 \times 3$ matrix, and $W - \alpha I_3 = (1-\alpha)(2-\alpha)^2$. Given $N(W-I) = Span\{(1,2,0)\}$, and $N(W-2I) = Span\{(2,0,3)\}$. Is $W$ diagonalizable matrix? Explain.

# Chapter 5

# Matrix Dot Product

In this chapter, we discuss the dot product only in $\mathbb{R}^n$. In addition, we give some results about the dot product in $\mathbb{R}^n$. At the end of this chapter, we get introduced to a concept called "Gram-Schmidt Orthonormalization".

## 5.1 The Dot Product in $\mathbb{R}^n$

In this section, we first give an example of the dot product in $\mathbb{R}^n$, and then we give three important results related to this concept.

**Example 5.1.1** Assume that $A = (2,4,1,3)$ and $B = (0,1,2,5)$ where $A, B \in \mathbb{R}^4$. Find $A \cdot B$.

**Solution:** To find $A \cdot B$, we need to do a simple vector dot product as follows:

$A \cdot B = (2,4,1,3) \cdot (0,1,2,5)$

$A \cdot B = 2 \cdot 0 + 4 \cdot 1 + 1 \cdot 2 + 3 \cdot 5$

$A \cdot B = 0 + 4 + 2 + 15$

Thus, $A \cdot B = 21$

**Result 5.1.1** If $W_1$ and $W_2$ in $\mathbb{R}^n$ and $W_1 \neq (0,0,0,\ldots,0)$ and $W_2 \neq (0,0,0,\ldots,0)$, and $W_1 \cdot W_2 = 0$, then $W_1$ and $W_2$ are independent.

**Result 5.1.2** If $W_1$ and $W_2$ in $\mathbb{R}^n$ are independent, then may/maybe not $W_1 \cdot W_2 = 0$. (i.e. Assume that $W_1 = (1,1,1,1)$ and $W_2 = (0,1,1,1)$, then $W_1 \cdot W_2 = 3$)

**Result 5.1.3** If $W_1, W_2, W_3$ and $W_4$ in $\mathbb{R}^n$ and none of them is $(0,0,0,\ldots,0)$, then we say that $W_1, W_2, W_3$ and $W_4$ are orthogonal if they satisfy the following conditions:

$$W_1 \cdot W_2 = 0$$

$$W_1 \cdot W_3 = 0$$

$$W_1 \cdot W_4 = 0$$

$$W_2 \cdot W_3 = 0$$

$$W_2 \cdot W_4 = 0$$

$$W_3 \cdot W_4 = 0$$

**Result 5.1.4** Assume that $W = (x_1, x_2, x_3, \ldots, x_n)$, then the squared-norm of $W$ is written as follows:
$$||W||^2 = x_1^2 + x_2^2 + x_3^2 + \cdots + x_n^2 = W \cdot W.$$

## 5.2 Gram-Schmidt Orthonormalization

In this section, we give one example that explains the concept of Gram-Schmidt Orthonormalization, and how it is related to what we have learned in chapters 2 and 3.

**Example 5.2.1** Given the following: $A = Span\{(1,0,1,1), (0,1,0,1), (0,1,1,1)\}$. Find the orthogonal basis for $A$.

**Hint:** The orthogonal basis means Gram-Schmidt Orthonormalization.

**Solution:** To find the orthogonal basis for $A$, we need to do the following steps:

**Step 1:** Find a basis for $A$.

$$\begin{bmatrix} 1 & 0 & 1 & 1 \\ 0 & 1 & 0 & 1 \\ 0 & 1 & 1 & 1 \end{bmatrix} -R_2 + R_3 \to R_3 \begin{bmatrix} 1 & 0 & 1 & 1 \\ 0 & 1 & 0 & 1 \\ 0 & 0 & 1 & 0 \end{bmatrix}$$ This is the Semi-Reduced Matrix.

Since we do now have a zero-row, then $dim(A) = 3$. To write a basis for $A$, it is recommended to choose rows in the above Semi-Reduced Matrix. Thus, a basis for $A$ is $\{(1,0,1,1), (0,1,0,1), (0,0,1,0)\}$.

**Step 2:** Vector spaces assumptions.

Since the basis for $A$ is $\{(1,0,1,1), (0,1,0,1), (0,0,1,0)\}$, we assume the following:

$$v_1 = (1,0,1,1)$$
$$v_2 = (0,1,0,1)$$
$$v_3 = (0,0,1,0)$$

**Step 3:** Use Gram-Schmidt Orthonormalization method.

Let's assume that the orthogonal basis for $A$ is $B_{Orthogonal} = \{W_1, W_2, W_3\}$. Now, we need to find $W_1, W_2$ and $W_3$. To find them, we need to do the following:

$W_1 = v_1 = (1,0,1,1)$

$W_2 = v_2 - \alpha_1 \cdot (previous\ W_{i's})$ where $\alpha_1 = \left(\frac{v_2 \cdot W_1}{||W_1||^2}\right)$

Thus, $W_2 = v_2 - \left(\frac{v_2 \cdot W_1}{||W_1||^2}\right) \cdot W_1$

$$W_2 = v_2 - \left(\frac{v_2 \cdot W_1}{||W_1||^2}\right) \cdot W_1$$

$$= (0,1,0,1) - \left(\frac{(0,1,0,1) \cdot (1,0,1,1)}{||(1,0,1,1)||^2}\right) \cdot (1,0,1,1)$$

$$= (0,1,0,1) - \frac{1}{3} \cdot (1,0,1,1)$$

$$= (0,1,0,1) - \left(\frac{1}{3}, 0, \frac{1}{3}, \frac{1}{3}\right) = \left(-\frac{1}{3}, 1, -\frac{1}{3}, \frac{2}{3}\right)$$

$W_3 = v_3 - \alpha_2 \cdot W_2 - \alpha_1 \cdot W_1$ where $\alpha_1 = \left(\frac{v_3 \cdot W_1}{||W_1||^2}\right)$ and $\alpha_2 = \left(\frac{v_3 \cdot W_2}{||W_2||^2}\right)$.

Thus, $W_3 = v_3 - \left(\frac{v_3 \cdot W_2}{||W_2||^2}\right) \cdot W_2 - \left(\frac{v_3 \cdot W_1}{||W_1||^2}\right) \cdot W_1$

$W_3 = (0,0,1,0) - \left(\frac{(0,0,1,0) \cdot \left(-\frac{1}{3}, 1, -\frac{1}{3}, \frac{2}{3}\right)}{\left\|\left(-\frac{1}{3}, 1, -\frac{1}{3}, \frac{2}{3}\right)\right\|^2}\right) \cdot \left(-\frac{1}{3}, 1, -\frac{1}{3}, \frac{2}{3}\right)$

$- \left(\frac{(0,0,1,0) \cdot (1,0,1,1)}{||(1,0,1,1)||^2}\right) \cdot (1,0,1,1)$

Thus, $W_3 = \left(-\frac{4}{15}, -\frac{1}{5}, \frac{11}{15}, -\frac{7}{15}\right)$.

Hence, the orthogonal basis for $A$ (Gram-Schmidt Orthonormalization) is
$\{(1,0,1,1), \left(-\frac{1}{3}, 1, -\frac{1}{3}, \frac{2}{3}\right), \left(-\frac{4}{15}, -\frac{1}{5}, \frac{11}{15}, -\frac{7}{15}\right)\}$.

## 5.3 Exercises

**1.** Let $D = Span\{(0,0,1,1), (1,0,1,1), (1,-1,1,0)\}$. Find the orthogonal basis for $D$.

**2.** Let $E = Span\{(1,1,-1,0), (0,1,1,1), (3,5,-1,2)\}$. Find the orthogonal basis for $E$.

# Answers to Odd-Numbered Exercises

## 1.10 Exercises

1. $x_4, x_5, x_6 \in \mathbb{R}$
$x_1 = \frac{3}{2}x_4 + x_5 - \frac{3}{2}x_6 + 15$
$x_2 = -\frac{3}{4}x_4 + \frac{3}{2}x_5 + \frac{3}{4}x_6 - \frac{13}{2}$
$x_3 = \frac{1}{2}x_4 - \frac{1}{2}x_6 + 5$

3. a. $[-3 \quad 3 \quad -1 \quad 5]$
   b. $\begin{bmatrix} 7 \\ 4 \\ 0 \end{bmatrix}$
   c. $-2$

5. a. $\begin{bmatrix} 1 & 0 \\ 2 & 1 \end{bmatrix} \begin{bmatrix} 2 & 0 \\ 0 & 1 \end{bmatrix} A = A_2$
   b. $\begin{bmatrix} \frac{1}{2} & 0 \\ 0 & 1 \end{bmatrix} \begin{bmatrix} 1 & 0 \\ -2 & 1 \end{bmatrix} A_2 = A$

7. $A$ is invertible (non-singular), $A^{-1} = \begin{bmatrix} 1 & 1 \\ \frac{3}{2} & 2 \end{bmatrix}$

9. $\dfrac{C_{42}}{\det(A)} = \dfrac{(-1)^{2+4} \det \begin{bmatrix} 2 & 2 & 1 \\ -2 & 2 & -1 \\ 1 & 12 & 4 \end{bmatrix}}{-1144} = -\dfrac{28}{1144} = -\dfrac{7}{286}$

11. $x_1 = 0$
$x_2 = -3$
$x_3 = 1$

## 2.5 Exercises

1. a. Basis for $M$ is $\{(1,-1,1),(-1,0,1)\}$.

   b. $(1,-3,5) \in M$ because there is a zero row that corresponds to $(1,-3,5)$ which means that $(1,-3,5)$ is dependent, and $(1,-3,5)$ can be written as a linear combination.

3. a. Let $v_1 = \begin{bmatrix} x \\ x+y \\ y \end{bmatrix} \in E$, and $v_2 = \begin{bmatrix} a \\ a+z \\ z \end{bmatrix}$, then

$v_1 + v_2 = \begin{bmatrix} x+a \\ (x+a)+(y+z) \\ y+z \end{bmatrix} \in E.$

For $m \in \mathbb{R}$ and $v_1 = \begin{bmatrix} x \\ x+y \\ y \end{bmatrix} \in E$, then $\alpha v_1 =$

$\begin{bmatrix} \alpha x \\ \alpha(x+y) \\ \alpha y \end{bmatrix} = \begin{bmatrix} \alpha x \\ \alpha x + \alpha y \\ \alpha y \end{bmatrix} \in E.$ Thus, $E$ is a subspace of $\mathbb{R}^3$.

   b. Basis for $E$ is $\{(1,1,0),(0,1,1)\}$.

   c. $E = Span\{(1,1,0),(0,1,1)\}$.

5. $x = 3$

7. $dim(K) = 2$

9. No, since the determinant of $\begin{bmatrix} -2 & 2 & 2 \\ 1 & -1 & 3 \\ 2 & -1 & 0 \end{bmatrix}$ equals zero, then the elements $\{(-2,1,2),(2,1,-1),(2,3,0)\}$ do not $\mathbb{R}^3$.

## 3.4 Exercises

1. a. Standard Matrix Representation of $T$ is
$$\begin{bmatrix} 1 & 0 & 0 \\ 1 & 0 & 0 \\ 0 & 0 & 1 \\ 0 & 1 & 0 \end{bmatrix}$$
   b. Basis for $Ker(T)$ is $\{(1,1,0,0), (0,0,0,1), (0,0,1,0)\}$.
Thus, $Ker(T)$ is $Span\{(1,1,0,0), (0,0,0,1), (0,0,1,0)\}$.
   c. Basis for $Range(T) = \{\begin{bmatrix} 1 & 1 \\ 0 & 0 \end{bmatrix}, \begin{bmatrix} 0 & 0 \\ 0 & 1 \end{bmatrix}, \begin{bmatrix} 0 & 0 \\ 1 & 0 \end{bmatrix}\}$.
Thus, $Ker(T)$ is $Span\{\begin{bmatrix} 1 & 1 \\ 0 & 0 \end{bmatrix}, \begin{bmatrix} 0 & 0 \\ 0 & 1 \end{bmatrix}, \begin{bmatrix} 0 & 0 \\ 1 & 0 \end{bmatrix}\}$.

3. $Ker(T) = Span\{-0.5 + x, -\frac{1}{3} + x^2\}$.

5. a. Standard Matrix Representation of $T$ is
$$\begin{bmatrix} 1 & -1 & 1 & -1 \\ 1 & 0 & 0 & 0 \\ 1 & -1 & 1 & -1 \\ 1 & 0 & 0 & 0 \end{bmatrix}$$
   b. Basis for $Ker(T)$ is $\{x + x^2, -x + x^3\}$.
Thus, $Ker(T)$ is $Span\{x + x^2, -x + x^3\}$.
   c. Basis for $Range(T) = \{\begin{bmatrix} 1 & 1 \\ 1 & 1 \end{bmatrix}, \begin{bmatrix} -1 & 0 \\ -1 & 0 \end{bmatrix}\}$.
Thus, $Ker(T)$ is $Span\{\begin{bmatrix} 1 & 1 \\ 1 & 1 \end{bmatrix}, \begin{bmatrix} -1 & 0 \\ -1 & 0 \end{bmatrix}\}$.

## 4.3 Exercises

1. Since $dim(E_{-1}) \neq 2$, $C$ is not diagonalizable.

3. $W$ is not diagonalizable because $dim(E_2)$ must be 2 not 1.

# 5.3 Exercises

1. The orthogonal basis for *A* (Gram-Schmidt Orthonormalization) is
$\{(0,0,1,1), (1,0,0,0), (0,-1,\frac{1}{2},-\frac{1}{2})\}$.

# Bibliography

[1] Badawi, A.: MTH 221/Linear Algebra. http://www.ayman-badawi.com/MTH221.html (2004). Accessed 18 Aug 2014.